
\input amstex
\documentstyle{amsppt}
\magnification=1200
\catcode`\@=11
\redefine\logo@{}
\catcode`\@=13

\define \bn{\Bbb N}
\define \bz{\Bbb Z}
\define \bq{\Bbb Q}
\define \br{\Bbb R}
\define \bc{\Bbb C}
\define \bh{\Bbb H}

\define \M{{\Cal M}}
\define\Ha{{\Cal H}}
\define\La{{\Cal L}}

\define\geg{{\goth g}}
\define\rk{\text{rk}~}


\define\mult{\text{mult}}
\define\re{\text{re}}
\define\im{\text{im}}

\define\0o{{\overline 0}}
\define\1o{{\overline 1}}

\TagsOnRight

\document

\topmatter

\title
On classification of Lorentzian Kac--Moody algebras
\endtitle

\author
Valery A. Gritsenko and
Viacheslav V. Nikulin \footnote{Supported by
Russian Fund of Fundamental Research (grant 00-01-00170).
\hfill\hfill}
\endauthor

\address
University Lille 1,
UFR de Mathematiques,
F-59655 Villeneuve d'Ascq Cedex,
France
\vskip1pt
POMI, St. Petersburg, Russia
\endaddress
\email
Valery.Gritsenko\@agat.univ-lille1.fr
\endemail

\address
Deptm. of Pure Mathem. The University of Liverpool, Liverpool
L69 3BX, UK;
\vskip1pt
Steklov Mathematical Institute,
ul. Gubkina 8, Moscow 117966, GSP-1, Russia
\endaddress
\email
vnikulin\@liv.ac.uk\ \
slava\@nikulin.mian.su
\endemail

\abstract
We discuss a general theory of
Lorentzian Kac--Moody algebras which should be a hyperbolic analogy of
the classical theories of finite-dimensional semi-simple and affine
Kac--Moody algebras. First examples of Lorentzian Kac--Moody algebras
were found by Borcherds. We consider general finiteness results about
the set of Lorentzian Kac--Moody algebras of the rank $\ge 3$,
and the problem of their classification. As an example,
we give classification of Lorentzian Kac--Moody algebras of the rank
three with the hyperbolic root lattice $S_t^\ast$,
symmetry lattice $L_t^\ast$, and the symmetry group
$\widehat{O}^+(L_t)$, $t\in \bn$, where
$$
H=\left(\smallmatrix 0&-1\\-1&0\endsmallmatrix\right),\
S_t=H\oplus \langle 2t \rangle=
\left(\smallmatrix 0&0&-1\\0&2t&0\\-1&0&0\endsmallmatrix\right),\
L_t=H\oplus S_t=
\left(\smallmatrix 0&0 &0& 0&-1\\
				 0&0 &0  &-1& 0\\
         0&0 &2t & 0& 0\\
				 0&-1&0  & 0& 0\\
        -1& 0&0  & 0& 0
\endsmallmatrix\right),
$$
and
$\widehat{O}^+(L_t)=\{g\in O^+(L_t)\ |\ g
\text{\ is trivial on $L_t^\ast/L_t$}\}$
is an extended paramodular group. Perhaps, this is the first example when
a large class of Lorentzian Kac--Moody algebras was classified.
\endabstract

\rightheadtext
{Lorentzian Kac--Moody algebras}
\leftheadtext{V.A. Gritsenko and V.V. Nikulin}
\endtopmatter

\document

\head
0. Introduction
\endhead

There are two well-known and very important in Mathematics and
Physics types of Lie algebras. The first is the class of
finite-dimensional semi-simple (or just finite) Lie algebras.
The second is the class of affine Kac--Moody algebras.
Both these classes belong to the general class of Kac--Moody Lie algebras.
See \cite{47} about their theory.

After fundamental results by R. Borcherds \cite{2} --- \cite{7},
it seems, we now know a right definition of the next type of
Lie algebras which should serve as a hyperbolic analogy of finite
(or elliptic type) and affine (or parabolic type) Lie algebras.
Here we call this type of Lie algebras as {\it Lorentzian Kac--Moody
algebras,} but one may prefer a different name, e. g. {\it hyperbolic of
Borcherds type}. Algebras of this type were very important for
Borcherds solution \cite{5} of the famous
Moonshine Conjecture of Conway and Norton \cite{20}
about modular properties of representations of the sporadic
finite simple group with the name Monster.

Main features of the theory of
Lorentzian Kac--Moody algebras $\geg$ are as follows
(see details of this definition in Sect. 1.4):

1) A Lorentzian Kac--Moody algebra $\geg$ is graded
$$
\geg\ =\ \bigoplus_{\alpha\in R}{\geg_\alpha}=
\geg_0 \bigoplus \left(\bigoplus_
{\alpha\in \Delta_+}\geg_\alpha\right) \bigoplus
\left(\bigoplus_{\alpha\in \Delta_+}
\geg_{-\alpha}\right)
$$
by a {\it hyperbolic root lattice $R$} which is a free $\bz$-module of a
finite rank equipped with an integral symmetric bilinear form which
is non-degenerate and
has exactly one negative square. For finite and
affine Kac-Moody algebras the root lattice $R$ is positive definite and
semi-positive definite respectively.

2) {\it Weyl group $W\subset O(R)$} is a reflection group
in the hyperbolic space (or a Lobachevsky space) $\La(R)$ related
with $R$. For finite and affine Kac--Moody algebras the Weyl group
is a finite reflection group and a discrete reflection group
in Euclidean space respectively.

3) The real part of the root system of the algebra $\geg$
which is some subset of roots in the root lattice $R$,
is some hyperbolic analogy of finite and affine root systems.
Like for finite and affine root systems,
it is defined by a fundamental chamber $\M\subset \La(R)$ of
$W$ and the set $P(\M)$ of roots of $R$ which are orthogonal to
faces of $\M$ of highest dimension ($P(\M)$ is the set of
{\it simple real roots} of the algebra $\geg$).

4) {\it Denominator identity of $\geg$} (due to Weyl, Kac and Borcherds)
defines an automorphic form $\Phi$ on a IV type
(in classification of \'E. Cartan) Hermitian
symmetric domain related with $R$ which should agree
with the Weyl group $W$. The denominator identity has the form
$$
\split
&\Phi=\exp{(-2\pi i (\rho,z))}
\prod_{\alpha\in \Delta_+}{\Bigl(1-\exp{\left(-2\pi i (\alpha,z)\right)}
\Bigr)^{\mult(\alpha)}}=\\
=&\sum_{w\in W}{\varepsilon (w)
\Bigl(\exp{\left(-2\pi i (w(\rho),z)\right)}\ -\hskip-10pt
\bigr.\sum_{a\in S\cap \br_{++}\M}
{\bigl.m(a)\exp{\left(-2\pi i (w(\rho+a),z)\right)}\Bigr)}}.
\endsplit
\tag{0.1}
$$
Its infinite product part is the infinite product
expansion of the automorphic form $\Phi$ where multiplicities
$\mult(\alpha)$ give dimensions of the root spaces $\geg_\alpha$
of $\geg$. Its infinite sum part is the Fourier expansion of
the automorphic form $\Phi$ which defines the algebra $\geg$
by generators and defining relations similar to of Killing,
Cartan, Chevalley and Serre. The element $\rho \in R\otimes \bq$
is the {\it Weyl vector} defined by the equality
$$
(\rho,\,\alpha)=-\alpha^2/2\ \  \text{for any\ }
\alpha \in P(\M).
\tag{0.2}
$$
For finite and affine Kac--Moody algebras the denominator identity
gives a polynomial and a Jacobi modular form respectively. We want
to keep this modularity property for Lorentzian Kac--Moody algebras.

5) The automorphic form $\Phi$ should be {\it reflective} which
means that the divisor of $\Phi$ should be union of rational
quadratic divisors orthogonal to roots. One can consider this
condition as the globalization on the whole Hermitian symmetric domain
of zeros of the infinite product expansion of $\Phi$ in the
neighborhood of the cusp at infinity of the Hermitian symmetric domain
where the infinite product converges. This reflectivity
property is valid in all known interesting cases, and it seems,
it distinguishes the most interesting
Lie algebras. Moreover, this additional condition is very important
for classification.

For hyperbolic case, one cannot satisfy conditions 1) --- 4)
inside usual class of Kac--Moody algebras. Fundamental discovery of
Borcherds is that one can satisfy conditions 1) --- 5) after
an appropriate generalization of Kac--Moody algebras to so called
{\it generalized Kac--Moody algebras}.
He constructed many examples of generalized Kac--Moody algebras
satisfying the conditions 1) --- 5)
above. See \cite{2} --- \cite{7}. We give the
most multi-dimensional and interesting Borcherds example in Sect. 1.3.
We follow Borcherds: Lorentzian Kac--Moody algebras which we consider
will be {\it generalized Kac--Moody algebras or superalgebras}
satisfying conditions 1) --- 5).  See details of the definition in
Sect. 1.4.

Lorentzian Kac--Moody algebras used by Borcherds for solution of
Moonshine Conjecture are graded by the unimodular hyperbolic plane
(of the rank two) $H=\pmatrix 0&-1\\-1&0\endpmatrix$, and their
denominator identity is an automorphic form $j_g(\tau_1)-j_g(\tau_2)$,
$\tau_1,\,\tau_2\in \bh$,
on the product $\bh\times \bh$ of two upper-half planes $\bh$ where $g$
is a conjugacy class of the Monster and $j_g(\tau)$ is an appropriate
normalized Hauptmodul. There are many papers and reviews written on
the subject of Borcherds solution \cite{5} of Moonshine Conjecture.
E. g. see \cite{6}, \cite{10}, \cite{27}, \cite{70}. We don't consider that
in the paper. We also don't consider relation of Lorentzian
Kac--Moody algebras with Vertex Algebras which is very important for
solution of Moonshine Conjecture.
E.g. see \cite{2} --- \cite{6}, \cite{11},
\cite{25}, \cite{27}, \cite{28}, \cite{49} about this subject.

In this paper we consider classification problem of
Lorentzian Kac--Moody algebras of the rank $\rk R\ge 3$.
One of the main properties of finite and affine
Kac--Moody algebras is that they are classified (by Killing).
Finite Kac--Moody algebras are classified by Dynkin diagrams,
and affine Kac--Moody algebras are classified by extended Dynkin
diagrams. In this paper we consider
similar property of Lorentzian Kac--Moody algebras of the rank
$\rk R\ge 3$.  We consider results and conjectures which show that
number of Lorentzian Kac--Moody algebras of the rank $\ge 3$ is in
essential finite. Thus, Lorentzian Kac--Moody algebras of the rank
$\ge 3$ are all exceptional (like exceptional Lie algebras of the
type ${\Bbb E}_n$, $n=6,\,7,\,8$ , ${\Bbb F}_4$,
${\Bbb G}_2$) and
can be classified in principle. We mention that for the ranks
one and two analogous classification problems are much simpler,
but we expect that number of cases is in essential infinite for
this case.

In \S 1 we consider general definitions and theory of
Lorentzian Kac--Moody algebras (Sects. 1.1 --- 1.4), and general
finiteness results and conjectures about Lorentzian Kac--Moody algebras
of the rank $\ge 3$ (Sects. 1.5 --- 1.7). Main result here is that
number of data 1) --- 3) in data 1) --- 4) is finite
(or is in essential finite)
for $\rk R\ge 3$. We also give some classification
results about data 1) --- 3) for the rank three case.
All these results are closely related with the old results of
the second author and Vinberg about finiteness of the set of arithmetic
reflection groups in hyperbolic spaces.
See \cite{58} --- \cite{61}, \cite{71} --- \cite{73}.

The main property of the data 1) --- 3) in data 1) --- 4) which permits
to get these results, is that
the fundamental chamber $\M$ of the Weyl group $W$ has finite or
almost finite volume. Exactly here one uses that the denominator
identity in 4) gives an automorphic form. More exactly, let
$Sym(\M)\subset O^+(R)$ be the symmetry group of $\M$. Then
the corresponding semi-direct
product $W\rtimes Sym(\M)$ has finite index in $O^+(R)$, and there exists
a non-zero $\rho\in R\otimes \bq$ (it is called a {\it generalized
Weyl vector}) such that the orbit $Sym(\M)(\rho)$ is finite. Hyperbolic
lattices $R$ having a reflection subgroup $W\subset O(R)$ with this
property are called {\it reflective}. Since $O^+(R)$ is arithmetic and
has a fundamental domain of finite volume, it follows that $\M$ is finite
of finite volume if $\rho^2<0$ {\it (elliptic type);}
$\M$ is finite of finite volume in
any angle with the vertex at infinity $\br_{++}\rho$ if $\rho^2=0$
{\it (parabolic type);}  $\M$ is finite of finite volume
in any orthogonal cylinder over
compact set in the hyperplane orthogonal to $\rho$ if $\rho^2>0$
{\it (hyperbolic type).} Using this property of the fundamental
chamber $\M$, one can prove that the number of reflective hyperbolic
lattices $R$ is finite for $\rk R\ge 3$.
See \cite{57} --- \cite{61}, \cite{63}, \cite{64}, \cite{66},
\cite{68}, \cite{71}---\cite{73}.
Existence of the  Weyl vector $\rho$ satisfying \thetag{0.2},
additionally permits to prove finiteness or almost finiteness (i.e.
finiteness up to a very simple equivalence relation)
of data 1) --- 3) in data 1) --- 4).

In Sect. 1.6 we consider finiteness conjectures and results about
data 4) and 5). Some observations from \cite{65} and \cite{40}
give a hope that number of these data is also finite for $\rk R\ge 3$.
It seems, some general results about infinite automorphic products and
their divisors by Borcherds \cite{7}, \cite{9} and
by Bruinier \cite{14} --- \cite{16} are also related with this subject.
We don't discuss these results of Borcherds and Bruinier in this paper.

In Sect. 2 we consider a concrete example of classification of
Lorentzian Kac--Moody algebras of the rank three
where we use general ideas and results from \S 1.
We give classification of Lorentzian Kac--Moody algebras with
the root lattice $S_t^\ast$ where
$$
S_t=H\oplus \langle 2t \rangle=
\pmatrix 0&0&-1\\0&2t&0\\-1&0&0\endpmatrix,
$$
the symmetry lattice $L_t^\ast$ where
$$
L_t=H\oplus S_t=
\pmatrix 0&0 &0& 0&-1\\
				 0&0 &0  &-1& 0\\
         0&0 &2t & 0& 0\\
				 0&-1&0  & 0& 0\\
        -1& 0&0  & 0& 0
\endpmatrix,
$$
and the symmetry group (of the automorphic form $\Phi$)
$$
\widehat{O}^+(L_t)=\{g\in O^+(L_t)\ |\ g
\text{\ is trivial on $L_t^\ast/L_t$}\}.
$$
Here $t\in \bn$, and $\oplus$ denotes the orthogonal sum of lattices,
and $\ast$ denotes the dual lattice.
We prove (Theorem 2.1.1) that there are exactly {\bf 29} Lorentzian
Kac--Moody algebras $\geg$ (or data 1) --- 5))
with the root lattice $S_t^\ast$, symmetry lattice $L_t^\ast$ and
the symmetry group $\widehat{O}^+(L_t)$. It seems, this result is
the first where a large class of Lorentzian Kac--Moody algebras was
classified. \S 2 is devoted to the outline of the proof of this
result. Actually, we prove much more general classification result about
reflective automorphic forms with the root lattice $S_t^\ast$, the
symmetry lattice $L_t^\ast$, the symmetry group $\widehat{O}^+(L_t)$
and an infinite product expansion similar to \thetag{0.1}
(Theorems 2.2.3 and 2.4.1). This result gives information about
Lorentzian Kac--Moody algebras with symmetry lattices $L$ which are
equivariant sublattices $L\subset L_t$ of the same rank
$\rk L=\rk L_t=5$. Here equivariant means that $O(L)\subset O(L_t)$.
Moreover, one can expect that the corresponding
{\it infinite product = infinite sum} identities
(similar to \thetag{0.1}) for these
reflective automorphic forms could be related with some interesting
algebras similar to generalized Kac--Moody algebras and superalgebras.

In this paper we outline and present all ideas of the proof of the
classification Theorem 2.1.1. We hope to present details of
the proof and calculations in forthcoming more longer publication.

All 29 Lorentzian Kac--Moody algebras of Theorem 2.1.1
had been constructed in \cite{41}. To construct these algebras,
one need to construct the corresponding 29 automorphic forms
\thetag{0.1}. They were
constructed in \cite{41} together with their infinite products and
infinite sums expansions which give multiplicities and
generators of the algebra.
We construct all 29 automorphic forms of Theorem 2.1.1 using
some variant for Jacobi modular forms of Borcherds exponential lifting.
See \cite{41} about this variant. This construction requires
delicate calculations with appropriate Jacobi
modular forms of two variables with integral Fourier coefficients.
We use classification of these forms obtained in \cite{33}, \cite{34}.
We present some of these results in \S 3: Appendix.

To proof completeness of the list of these 29 automorphic forms, we
find (Theorem 2.3.2) all reflective hyperbolic lattices $S_t$
of the rank three. Only for them an automorphic form may exist.
In particular, $t\le 105$ for them. When the lattice $S_t$ is reflective,
we find all possible data 1) --- 3) for it and predict
the divisor of the reflective automorphic form $\Phi$.
Then we can see that one of
29 automorphic forms of Theorem 2.1.1 has the same divisor and
coincides with $\Phi$ by Koecher principle.

In general, in \cite{68} all reflective
hyperbolic lattices of the rank three were classified:
there are 122 main elliptic and 66 main hyperbolic types
(there are no main parabolic types). These results and reasonable
number of cases give a hope that all Lorentzian Kac--Moody algebras
of the rank three will be classified in a future.
Finiteness results for $\rk R\ge 3$ give a hope that the
same can be done for all ranks $\rk R\ge 3$.
We expect that number of cases drops when the rank is increasing:
There are no algebras when $\rk R$ is sufficiently high.
Borcherds example which we present in Sect. 1.3 is
the most highest dimensional known example of Lorentzian Kac--Moody
algebras. For this example $\rk R=26$.

In this paper, we don't tuch possible Physical applications of
Lorentzian Kac--Moody algebras. One can find some of them in
physical papers \cite{17}, \cite{18}, \cite{21} --- \cite{23},
\cite{28}, \cite{38}, \cite{42}, \cite{45},
\cite{51} --- \cite{53}, \cite{56}.

\smallpagebreak

This paper is written during our stay in the University of Lille 1,
the University of Liverpool, Steklov Mathematical Institute, Moscow
and St. Petersburg, Max-Planck-Institut f\"ur Mathematik, Bonn, and
Newton Institute for Mathematical Sciences, Cambridge.
We are grateful to the Institutes for hospitality.

\head
1. A Theory of Lorentzian Kac--Moody Algebras and general finiteness
results and conjectures
\endhead

We start with a variant of Theory of Lorentzian Kac--Moody algebras
which one can consider as a hyperbolic analogy of classical theories of
finite and affine Kac--Moody algebras. Here we follow
Borcherds and \cite{36}, \cite{40}, \cite{41}, \cite{64}, \cite{67}.

\subhead
1.1. Some general results on Kac--Moody algebras
\endsubhead
One can find all definitions and details of this section
in the classical book by Kac \cite{47}.

A {\it generalized Cartan matrix} $A$ is an integral square
matrix of a finite rank which has only $2$ on the diagonal and
non-positive integers out of the diagonal.
We shall consider only {\it symmetrizable}
generalized Cartan matrices $A$. It means that there exists a
diagonal matrix $D$ with positive rational diagonal coefficients
such that $B=DA$ is integral and symmetric. Then $B$ is called the
{\it symmetrization} of $A$. By definition,
$\text{sign}(A)=\text{sign}(B)$. We shall suppose that $A$ is
{\it indecomposable} which means that there does not exist
a decomposition $I=I_1\cup I_2$ of the set $I$ of indices of $A$
such that $a_{ij}=0$ if $i\in I_1$ and $j\in I_2$.

Each generalized Cartan matrix $A$ defines a Kac--Moody
Lie algebra $\geg(A)$ over $\bc$.
The {\it Kac--Moody algebra} $\geg(A)$ is defined
by the set of generators and defining relations prescribed by the
generalized Cartan matrix $A$. They are due to V. Kac and R. Moody. In
fact, they are a natural generalization of classical results by
Killing, Cartan, Weyl, Chevalley and Serre about finite-dimensional
semi-simple Lie algebras. One should introduce
the set of generators $h_i$, $e_i$, $f_i$, $i \in I$,
with defining relations
$$
\cases
&[h_i,\,h_j]=0,\ \ [e_i,\,f_i]=h_i, \ \ [e_i,f_j]=0,\ \text{if\ }i\not=j,\\
&[h_i,\,e_j]=a_{ij}e_j,\ \ [h_i,\,f_j]=-a_{ij}f_j,\\
&(ad~e_i)^{1-a_{ij}}e_j=(ad~f_i)^{1-a_{ij}}f_j=0,\ \text{if}\ i\not=j.
\endcases
\tag{1.1.1}
$$
The algebra $\geg(A)$ is simple or almost simple:
it is simple after factorization by some known central ideal.

We mention some general features of the theory of Kac--Moody
algebras $\geg(A)$.

1. The symmetrization $B$ defines a free $\bz$-module
$Q=\sum_{i\in I}\bz\alpha_i$ with generators $\alpha_i$, $i\in I$,
equipped with symmetric bilinear form $\left((\alpha_i,\,\alpha_j)\right)=B$
defined by the symmetrization $B$.
The $Q$ is called {\it root lattice}.
The algebra $\geg (A)$ is {\it graded by
the root lattice $Q$} (by definition, generators $h_i$, $e_i$, $f_i$ have
weights $0$, $\alpha_i$, $-\alpha_i$ respectively):
$$
\geg(A)=\bigoplus_{\alpha\in Q}{\geg_\alpha}=
\geg_0 \bigoplus \left(\bigoplus_
{\alpha\in \Delta_+}\geg_\alpha\right) \bigoplus
\left(\bigoplus_{\alpha\in -\Delta_+}
\geg_{\alpha}\right)
\tag{1.1.2}
$$
where $\geg_\alpha$ are finite dimensional linear spaces,
$[\geg_\alpha,\,\geg_\beta]\subset \geg_{\alpha +\beta}$,
$\geg_0\equiv Q\otimes \bc$ is commutative,
and is called {\it Cartan subalgebra}.
An element $0\not=\alpha\in Q$ is called {\it root} if $\geg_\alpha\not=0$.
The $\geg_\alpha$ is called {\it root space} corresponding to $\alpha$.
The dimension $\mult(\alpha)=\dim \geg_\alpha$ is called
{\it multiplicity} of the root $\alpha$. In \thetag{1.1.2},
$\Delta \subset Q$ is the set of all roots. It is divided in
the set of {\it positive} $\Delta_+\subset \sum_{i\in I}{\bz_+\alpha_i}$
and {\it negative} $-\Delta_+$ roots.
A root $\alpha \in \Delta$ is called
{\it real} if $(\alpha,\alpha)>0$.
Otherwise (if $(\alpha,\alpha)\le 0$), it is called {\it imaginary}.
Every real root $\alpha$ defines {\it a reflection}
$s_\alpha:x\mapsto x-\left(2(x,\alpha)/(\alpha,\alpha)\right)\alpha$,
$x\in Q$. All reflections $s_\alpha$ in real roots generate
{\it Weyl group} $W\subset O(S)$. The set of roots $\Delta$
and multiplicities of roots are $W$-invariant.

2. One has {\it Weyl---Kac
denominator identity} which permits to
calculate multiplicities of roots:
$$
e(-\rho)\prod_{\alpha\in \Delta_+}{(1-e(-\alpha))^{\mult (\alpha)}} =
\sum_{w\in W}{\det(w)e(-w(\rho))}.
\tag{1.1.3}
$$
Here $e(\,\cdot\, )\in \bz[Q]$ are formal exponents where
$\bz[Q]$ is the group ring of the root lattice $Q$. The
$\rho$ is called {\it Weyl vector} and is defined by
the condition $(\rho,\,\alpha_i)=-(\alpha_i,\,\alpha_i)/2$
for any $i\in I$.

The identity \thetag{1.1.3} is combinatorial, and direct
formulae for multiplicities \linebreak
$\mult (\alpha)$ are unknown in general.
One approach to solve this problem is to replace the
formal function \thetag{1.1.3} by non-formal
one (e. g. replacing formal exponents by non-formal ones)
to get a function with ``good'' properties. These good
properties may help to find the formulae for multiplicities.

\subhead
1.2. Finite and affine cases
\endsubhead
There are two cases when we have very clear picture (or Theory)
of Kac--Moody algebras.

{\it Finite case:} The generalized Cartan matrix $A$ is positive definite,
$A>0$. Then $\geg(A)$ is finite-dimensional, and we get the classical
theory of {\it finite-dimensional semi-simple Lie algebras.}

{\it Affine case:} The generalized Cartan matrix $A$ is semi-positive
definite, $A\ge 0$. Then $\geg(A)$ is called {\it affine}.

For both these cases we have three very nice properties:

(I) There exists classification of all possible generalized
Cartan matrices $A$ and the corresponding algebras
$\geg(A)$: They are classified by Dynkin (for finite case)
and by extended Dynkin (for affine case) diagrams.

(II) In the denominator identity \thetag{1.1.3},  formal exponents may
be replaces by non-formal ones to give a function with nice properties:
For finite case this gives a polynomial. For affine case this gives a
Jacobi modular form. Using these properties (or directly), one
can find all multiplicities.

(III) Both these cases have extraordinary importance in
Mathematics and Phy\-sics.

\smallpagebreak

We want to construct similar Theory
for {\it Lorentzian (or hyperbolic) case}
when the generalized Cartan matrix $A$ is {\it hyperbolic:}
it has exactly one negative square, all its other squares are
either positive or zero. There are plenty of hyperbolic
generalized Cartan matrices, it is impossible to
find all of them and classify. On the other hand, probably
not all of them give interesting Kac--Moody algebras,
and one has to impose natural conditions on these matrices.

\smallpagebreak

\subhead
1.3. Lorentzian case. Borcherds example
\endsubhead
We have the following key example due to R. Borcherds
\cite{3}---\cite{6}.

For Borcherds example, the root lattice $Q=S$ where $S$ is
a hyperbolic even unimodular lattice $S$ of signature
$(25,1)$. Here ``even'' means that $(x,x)$ is even for any $x\in S$.
``Unimodular'' means that the dual lattice $S^\ast$ coincides with  $S$,
equivalently, for a basis $e_1,\dots, e_{26}$ of $S$
the determinant of the Gram matrix $\big((e_i,\,e_j)\big)$
is equal to $\pm 1$. A lattice $S$ with these properties
is unique up to isomorphism.
For Borcherds example, the Weyl group $W$ is generated by reflections
$s_\alpha:x\mapsto x-(x,\alpha)\alpha$, $x\in S$, in all elements
$\alpha \in S$ with $\alpha^2=2$.
The group $W$ is discrete in the hyperbolic space
$\La (S)=V^+(S)/\br_{++}$ where $\br_+$ and $\br_{++}$ denote the sets
of non-negative and positive real numbers respectively.
Here $V^+(S)$ is the positive cone, i. e.
a half of the cone $V(S)=\{x\in S\otimes \br \ |\ x^2<0\}$ of the
hyperbolic lattice $S$.
The $\La (S)$ is the set of rays in $V^+(S)$.

A fundamental chamber $\M\subset \La (S)$ for $W$ is defined by the set $P$
of elements $\alpha\in S$ with $\alpha^2=2$ which are {\it orthogonal to}
$\M$. It has the following description
due to Conway \cite{19}.
There exists an orthogonal
decomposition $S=[\rho,\,e]\oplus L$ where the Gram matrix
of elements $\rho,\,e$ is equal to
$H=\left(\matrix0&-1\\-1&0\endmatrix\right)$
(in particular, $(\rho,\rho)=0$),
and $L$ is the Leech lattice, i. e. positive definite even unimodular
lattice of the rank $24$ without elements with square $2$.
The set $P$ of roots
which are orthogonal to the fundamental chamber $\M$ (or the
set of {\it simple real roots}) of $W$ is equal to
$$
P=\{\alpha \in S\ |\ (\alpha,\,\alpha)=2\ \and\ (\rho,\alpha)=-1 \}.
\tag{1.3.1}
$$
It means that the fundamental chamber $\M\subset \La(S)$ is equal to
$$
\M=\{\br_{++}x\in\La(S)\ |\ (x,\ P)\le 0\}
\tag{1.3.2}
$$
and $P$ is a minimal set with this property. We mention that the
fundamental chamber $\M$ has {\it ``almost finite'' volume}. It means
that $\M$ is finite in any angle with the center at infinity
$\br_{++} \rho$ of the hyperbolic space $\La(S)$.

The matrix
$$
A=\bigl((\alpha,\,\alpha^\prime )\bigr),\ \ \alpha,\,\alpha^\prime \in P
\tag{1.3.3}
$$
is a generalized Cartan matrix and $\rho$ is the Weyl vector:
$$
(\rho,\alpha)=-(\alpha,\alpha)/2, \ \forall \, \alpha \in P.
\tag{1.3.4}
$$
Thus, $A$ defines the Kac--Moody algebra $\geg (A)$
graded by the hyperbolic lattice $S$. But the algebra
$\geg (A)$ is not the algebra which is considered for Borcherds
example. One has to {\it ``correct''}  the algebra $\geg (A)$.

We have the classical $SL_2(\bz)$-modular cusp
form $\Delta$ of the weight $12$ on the upper-half plane $\text{Im}\,\tau>0$:
$$
\Delta=q\prod_{n=1}^{\infty}(1-q^n)^{24}=
\sum_{m\ge 0}\tau(m)q^m,
\tag{1.3.5}
$$
where $q=exp(2\pi i \tau)$. We have
$$
\Delta^{-1}=\sum_{n\ge 0}p_{24}(n)q^{n-1}
\tag{1.3.6}
$$
where $p_{24}(n)$ are positive integers.
Borcherds \cite{4} proved the {\it identity}
$$
\Phi(z)=\exp{(-2\pi i (\rho,z))}\prod_{\alpha\in \Delta_+}
{(1-\exp{(-2\pi i (\alpha,z))})^{p_{24}(1-(\alpha,\alpha)/2)}}=
$$
$$
\sum_{w\in W}{\det(w)\sum_{m>0}{\tau (m)\exp{(-2\pi i (w(m\rho),z))}}}.
\tag{1.3.7}
$$
Here $\Delta_+=\{\alpha \in S\,|\,\alpha^2=2\ \and\ (\alpha,\,\rho)<0\}\cup
(S\cap \overline{V^+(S)}-\{0\})$. The variable
$z$ runs through the {\it complexified positive cone}
$\Omega(V^+(S))=S\otimes\br\,+\,i V^+(S )$.
Moreover, Borcherds \cite{6}, \cite{7} proved that the
function $\Phi(z)$ {\it is an automorphic form of weight $12$}
with respect to the
group $O^+(T)$ where $T=H\oplus S$ is the extended lattice of
the signature $(26,2)$ (we denote as $\oplus$ the orthogonal sum).
The group $O^+(T)$ naturally acts in the Hermitian
symmetric domain of type IV
$$
\Omega(T)=\{\bc\omega \subset T\otimes \bc\ |\ (\omega,\,\omega)=0\ \and\
(\omega,\,\overline{\omega})<0 \}_0,
\tag{1.3.8}
$$
which has canonical identification with $\Omega(V^+(S))$ as follows:
$z\in\Omega (V^+(S))$ defines the element
$\bc\omega_z\in \Omega(T)$
where $\omega_z =\big((z,z)/2\big)e_1+
e_2\oplus z \in T\otimes \bc$
and $e_1$, $e_2$ is the basis of the lattice $H$ with
the Gram matrix $H$ above. Here ``automorphic of
the weight 12'' means that the function
$\widetilde{\Phi}(\lambda \omega_z)=\lambda^{-12}\Phi(z)$,
$\lambda \in \bc^\ast$,
is homogeneous of the degree $-12$ (it is obvious)
in the homogeneous cone
$\widetilde{\Omega(T)}$ over $\Omega(T)$,
and $\widetilde{\Phi}(g\omega)=\det(g)\widetilde{\Phi}(\omega)$
for any $\omega \in \widetilde{\Omega(T)}$ and any
$g\in O^+(T)$ where $O^+(T)$ is the subgroup of index $2$ of
the group $O(T)$ which keeps the connected component \thetag{1.3.8}
(marked by $0$).

The identity \thetag{1.3.7} looks very familiar to
the form \thetag{1.1.3} of the denominator identity for Kac--Moody
algebras, but it has some difference.

To interpret \thetag{1.3.7} as a denominator identity of a Lie
algebra, Borcherds introduced \cite{3}
{\it generalized Kac--Moody algebras}
$\geg(A^\prime )$ which correspond to more general matrices $A^\prime$ than
generalized Cartan matrices. Here we shall call them as
{\it generalized Cartan--Borcherds matrices}.
Difference is that a
generalized Cartan--Borcherds matrix $A^\prime$ may also have
non-positive real elements
$a_{ij}\le 0$ on the diagonal and out of the diagonal,
but all $a_{ij}\in \bz$ if $a_{ii}=2$.
A definition of the generalized Kac--Moody
algebra $\geg (A^\prime )$ corresponding to a
generalized Cartan--Borcherds matrix $A^\prime$ is similar
to \thetag{1.1.1}. One should replace the last line of \thetag{1.1.1} by
$$
(\text{ad}~e_i)^{1-a_{ij}}e_j=(\text{ad}~f_i)^{1-a_{ij}}f_j=0\
\text{if}\ i\not=j\ \text{and}\ a_{ii}=2,
\tag{1.3.9}
$$
and add the relation
$$
[e_i,\,e_j]=[f_i,\,f_j]=0\ \text{if}\ a_{ij}=0.
\tag{1.3.10}
$$
Borcherds showed that generalized Kac--Moody algebras have similar
properties to ordinary Kac--Moody algebras.
They also have a
denominator identity which has more general form than \thetag{1.1.3}
and includes \thetag{1.3.7} as a particular case.

The identity \thetag{1.3.7} is the denominator identity for
the generalized Kac--Moody algebra $\geg(A')$ where $A'$ is the
generalized Cartan--Borcherds matrix equals to the Gram
matrix $A^\prime=\bigl((\alpha,\,\alpha^\prime )\bigr)$,
$\alpha,\,\alpha^\prime \in P^\prime$ where
$$
P^\prime=P\cup 24\rho \cup 24(2\rho) \cup \dots \cup 24(n\rho)\cup\cdots
\tag{1.3.11}
$$
is the sequence of elements of the lattice $S$.
Here $24(n\rho)$ means that we take the element $n\rho$ twenty four
times to get the Gram matrix $A^\prime$. See details in
\cite{3}, \cite{4}.

In \thetag{1.3.11}, the set $P^\prime$ defining $A^\prime$
is called the {\it set of simple roots}. It is divided in the set
${P^\prime}^{\re}=P$, described in \thetag{1.3.1}, of
{\it simple real roots} (they are orthogonal
to the fundamental chamber $\M$ of the Weyl group $W$ and have
positive square) and is the same
as for the ordinary Kac--Moody algebra $\geg(A)$ defined by the
generalized Cartan matrix $A$ in \thetag{1.3.3}. The additional sequence
$$
{P^\prime}^{\im}=24\rho \cup 24(2\rho) \cup \dots
\cup 24(n\rho)\cup\cdots
\tag{1.3.12}
$$
of $P^\prime$ (elements of ${P^\prime}^{\im}$ have zero square)
is defined by the Fourier coefficients in the sum part of
the identity \thetag{1.3.7}. For example, $24$ in
\thetag{1.3.12} is defined by the $24$ in \thetag{1.3.5}.
Together ${P^\prime}^{re}$ and ${P^\prime}^{\im}$
define the generalized Cartan--Borcherds matrix
$A^\prime$ and the generalized Kac--Moody algebra $\geg(A^\prime)$.

Borcherds example is very fundamental and beautiful. It has
important applications in Mathematics, e. g. for moduli spaces of
K3 and Enriques surface (see \cite{8},
\cite{13} and \cite{54}), and in
Physics (e. g. in String Theory): it gives the Lie algebra of
physical states of the Vertex Algebra of dimension 26 defined by
the hyperbolic lattice $S$. See \cite{2}, \cite{6}, \cite{11}.

Many other similar examples of generalized Kac--Moody algebras
and superalgebras graded by hyperbolic lattices and related with
automorphic forms
were found by Borcherds in \cite{3} ---  \cite{9},
\cite{12}. One of them when Lie algebras are graded
by the unimodular even hyperbolic plane $H$ is very important for Borcherds
proof of Moonshine Conjecture about modular properties of
representations of the Monster sporadic finite simple group.
See \cite{5}, \cite{10} and reviews \cite{27}, \cite{70} about
this subject.

\subhead
1.4. A Theory of Lorentzian Kac--Moody algebras
\endsubhead
Analyzing Borcherds example, one can suggest a general
class of Lorentzian Kac--Moody algebras (or automorphic hyperbolic
Kac--Moody algebras), see \cite{36}, \cite{40}, \cite{41},
\cite{62}, \cite{64}, \cite{67}. They are defined by
{\it data (1)---(5)} below:

(1) A {\it hyperbolic lattice} $S$ (i. e. a non-degenerate integral
symmetric bilinear form of signature $(n,1)$).
It is the {\it root lattice} for the Lie algebra we want to construct.
We follow \cite{57} in notations and terminology related with lattices.

\smallpagebreak

(2) A {\it reflection group}
$W\subset O(S)$. It is generated by reflections in some set of
roots of $S$. We remind that $\alpha\in S$ is called {\it root} if
$\alpha^2>0$ and $\alpha^2\,|\,2(\alpha,\,S)$. Any root defines
a {\it reflection}
$s_\alpha: x\mapsto x-(2(x,\,\alpha)/\alpha^2)\alpha$, $x\in S$, which
gives an automorphism of the lattice $S$. The group $W$ is the
{\it Weyl group} for the Lie algebra. We assume that the
{\it Weyl group $W$ is
non-trivial}. (For trivial $W$ the definition has to be changed,
see \cite{40}.)

\smallpagebreak

(3) {\it A set $P=P(\M)$ of orthogonal roots
to the  fundamental chamber
$\M\subset \La(S)=V^+(S)/{\br_{++}}$ of $W$} and its subdivision
$P=P_\0o\cup P_\1o$ by subsets of {\it even} and {\it odd} roots respectively.
The set $P$ of roots of $S$ should have the
property
$$
\M=\{\ \br_{++}x\in \La(S)\ |\ (x,\,P)\le 0\ \}
\tag{1.4.1}
$$
(the same as \thetag{1.3.2}) and should be minimal having this property (i. e.
each face of $\M$ of highest dimension is orthogonal to exactly one
element from $P$, and each element of $P$ is orthogonal to a face of
$\M$ of highest dimension). We additionally require that
$2\alpha$ is a root of $S$ if $\alpha\in P_\1o$ (thus
$\alpha^2\,|\,(\alpha,\,S)$).
Moreover, the set $P$ should have a
{\it Weyl vector} $\rho\in S\otimes \bq$
which means that it satisfies the condition
$$
(\rho,\,\alpha)=-\alpha^2/2,\ \ \forall \alpha \in P
\tag{1.4.2}
$$
(it is the same as \thetag{1.3.4}).\footnote{Existence of the Weyl
vector $\rho$ from $S\otimes \bq$ is a very strong condition on
the data (1)---(3). It would be better to call it as
{\it a lattice Weyl vector.} We don't do it to simplify
terminology. It is possible, that there exists a more general theory of
Lorentzian Kac--Moody algebras when a Lie algebra can be graded
by a degenerate hyperbolic lattice.  Then the Weyl vector should be
a linear function on the degenerate hyperbolic lattice, satisfying
\thetag{1.4.2}, and it always exists for an appropriate
grading. It is easy to see that one should add a one-dimensional
kernel only. On the other hand, if this theory
exists, it seems, it is so general that the classification problem has
no  sense. Thus, this theory is not so interesting for us.}
The sets $P$, $P_\0o$ and $P_\1o$ are the sets of {\it simple real roots,
even simple real roots} and {\it odd simple real roots} respectively for
the Lie algebra. The set $P$ is not empty since $W$ is non-trivial.
From \thetag{1.4.2}, the Weyl vector $\rho$ is not zero too.

\smallpagebreak

The main invariant of the data (1)---(3) is the {\it generalized
Cartan matrix}
$$
A=\left({2(\alpha, \alpha^\prime)\over (\alpha,\alpha)}\right),\ \ \
\alpha,\,\alpha^\prime\in P.
\tag{1.4.3}
$$
It defines data (1)---(3) up to some very clear equivalence.

\smallpagebreak

(4) A {\it holomorphic automorphic form} $\Phi(z)$ of some
weight $k$, ($k\in \bz/2$)
on a IV type Hermitian symmetric domain, $z\in \Omega(V^+(S))=\Omega (T)$,
with respect to a subgroup $G\subset O^+(T)$ of finite index
({\it the symmetry group} of the Lie algebra) of an extended
lattice $T=H(m)\oplus S$ ({\it the symmetry lattice} of the Lie algebra)
where $H(m)=\left(\matrix0&-m\\-m&0\endmatrix\right)$, $m\in \bn$.
(See \cite{40} for more general definition.)
Definition of the automorphic form $\Phi$ is the same as for the
Borcherds example in Sect. 1.3.
The only changes are: we identify
$\Omega(V^+(S))=\Omega(T)$ by the formula
$$
z\mapsto \bc \omega_z,\ \ z\in \Omega(V^+(S)),\ \omega_z\in
\widetilde{\Omega}(T),
\tag{1.4.4}
$$
where $\omega_z=(z,z)e_1/2+e_2/m\oplus z$ and
$e_1,\,e_2$ is the basis of the lattice $H(m)$ with the matrix
above, and we assume that $\widetilde{\Phi}(\lambda\omega_z)=
\lambda^{-k}\Phi(z)$, $\lambda \in \bc^\ast$, on the cone
$\widetilde{\Omega}(T)$ has the property:
$\widetilde{\Phi}(g\omega)=\chi(g)\widetilde{\Phi}(\omega)$
for any $g\in G$ where $\chi:G\to \bc^\ast$ is a character or a multiplier
system with the kernel of finite index in $G$.

The automorphic form $\Phi$
should have Fourier expansion of the form
of denominator identity for a generalized Kac--Moody algebra
with hyperbolic generalized Cartan--Borcherds matrix and which is
agree with previous data (1) --- (3). This form is
$$
\Phi(z)=\sum_{w\in W}{\varepsilon (w)
\Bigl(\exp{\left(-2\pi i (w(\rho),z)\right)}\ -\hskip-10pt
\bigr.\sum_{a\in S\cap \br_{++}\M}
{\bigl.m(a)\exp{\left(-2\pi i (w(\rho+a),z)\right)}\Bigr)}}
\tag{1.4.5}
$$
where $\varepsilon:W\to \{\pm 1 \}$ is a quadratic character (i. e.
a homomorphism) such that
$\varepsilon(s_\alpha)=(-1)^{1+\overline i}$ if
$\alpha\in P_{\overline i}$ and $\overline i = \0o,\ \1o$.
All Fourier coefficients $m(a)$ should be integral.

Let
$$
H=\{\ g\in O^+(S)\ |\ \Phi(g(z))=\pm \Phi(z)\ \}.
\tag{1.4.6}
$$
Since $\Phi$ is an automorphic form, this subgroup has finite index in
$O^+(S)$. We have $W\subset H$. Let
$$
Sym(\M)=\{\ g\in O^+(S)\ |\ g(\M)=\M \}
\tag{1.4.7}
$$
be the symmetry group of the fundamental chamber $\M$ and
$$
Sym(P_\1o\subset P)=\{\ g\in O^+(S)\ |\ g(P)=P\ \&\ g(P_\1o)=P_\1o\ \}
\tag{1.4.8}
$$
its subgroup which keeps invariant the set $P$ of orthogonal roots to
$\M$, and its subset $P_\1o$ as well.

(*) {\it We assume that there exists a subgroup
$A\subset Sym(P_\1o\subset P)$ such that
$A\subset H$ and the semi-direct product $W\rtimes A$ has finite
index in $H$.} It follows that $A$ is a subgroup of finite index
in $Sym(P_\1o\subset P)$ and in $Sym(\M)$, and
$W\rtimes A$ has finite index in $O^+(S)$. In particular,
subgroups $W\rtimes Sym(P_\1o\subset P)$ and
$W\rtimes Sym(\M)$ have finite index in $O(S)$.
The automorphic form $\Phi$ defines the set of
{\it simple imaginary roots} and gives {\it the denominator identity}
of the Lie algebra.
Using automorphic properties
of $\Phi(z)$, it is good to {\it calculate the infinite product
part of the denominator identity}
$$
\Phi(z)=\exp{(-2\pi i (\rho,z))}
\prod_{\alpha\in \Delta_+}{\Bigl(1-\exp{\left(-2\pi i (\alpha,z)\right)}
\Bigr)^{\mult(\alpha)}}
\tag{1.4.9}
$$
which gives multiplicities $\mult(\alpha)$ of roots $\alpha$ of the
Lie algebra. See below. There are many cases known when it is possible.

\smallpagebreak

Like for Borcherds example, already {\it data (1)---(4) define a
generalized Kac--Moody algebra or superalgebra.}
See the definition below. But it was understood that it is good
to suppose (at least, to have finiteness results)
the additional condition:

\smallpagebreak

(5) The automorphic form $\Phi$ on the domain
$\Omega (V^+(S))=\Omega (T)$
should be {\it reflective}. It means that the divisor of $\Phi$
is union of rational quadratic divisors which are orthogonal to some
roots of the extended lattice $T$. Here for a root
$\alpha \in T$ (the definition of a root of $T$ is the same as for
the lattice $S$) the {\it quadratic divisor orthogonal to} $\alpha$
is equal to
$$
D_\alpha=\{\bc\omega \in \Omega (T)\ |\ (\omega,\,\alpha)=0\}.
\tag{1.4.10}
$$

The property (5) is valid for Borcherds example above and for
the most part of known cases. Moreover, it is true in the neighbourhood
of the cusp where the infinite product \thetag{1.4.9} converges.
Thus, we want it to be true globally.

\smallpagebreak

Below we give the definition of a generalized Kac--Moody
superalgebra $\geg$ corresponding to data (1)---(4). It is given by
the sequence $P^\prime \to S$ {\it of simple roots}.
This sequence is divided in a sequence ${P^\prime}^{\re}$
of {\it simple real roots} and a sequence ${P^\prime}^{\im}$
of {\it simple imaginary roots}. Both these sequences are also
divided in the sequences of {\it even and odd roots}
marked by $\0o$ and $\1o$ respectively. We set ${P^\prime}^{\re}=P$,
${P^\prime}^{\re}_{\overline i}= P_{\overline i}$ where $P$,
$P_{\overline i}$, $\overline i =\0o,\ \1o$,
are defined in the datum (3).

For any primitive $0\not=a\in S\cap \br_{++}\M$
with $(a,\,a)=0$, one should find $\tau(na)\in \bz$, $n\in \bn$,
from the identity with a formal variable $t$:
$$
1-\sum_{k\in \bn}{m(ka)t ^k}=\prod_{n\in \bn}{(1-t^n)^{\tau(na)}}.
\tag{1.4.11}
$$
We set
$$
\split
{P^\prime}^{\im}_{\0o}&=\{m(a)a\ |\
a\in S\cap\br_{++}\M,\ (a,a)<0\ \text{and}\
m(a)>0\}\cup\\
&\cup \{\tau(a)a\ |\ a\in S\cap\br_{++}\M,\ (a,a)=0\ \text{and}\
\tau (a)>0\};\\
{P^\prime}^{im}_{\1o}&=
\{-m(a)a\ |\ a\in S\cap\br_{++}\M,\ (a,a)<0\ \text{and}\
m(a)<0\}\cup\\
&\cup \{-\tau(a)a\ |\ a\in S\cap\br_{++}\M,\ (a,a)=0\ \text{and}\
\tau (a)<0\}.\\
\endsplit
\tag{1.4.12}
$$
Here $ka$ means that we repeat the element $a$ exactly
$k$ times for the sequence. The generalized Kac--Moody superalgebra
$\geg$ is a Lie superalgebra generated by $h_r$, $e_r$, $f_r$ where
$r \in P^\prime$.
All generators  $h_r$ are even, generators $e_r$, $f_r$ are even
(respectively odd) if $r$ is even (respectively odd). They have the
defining relations 1) --- 5) below:

\smallpagebreak

1) The map $r\mapsto h_r$ for $r\in P^\prime$ gives an embedding of
$S\otimes \bc$ to $\geg$ as an Abelian subalgebra (it is even);

2) $[h_r,\,e_{r^\prime}]=(r,\,r^\prime)e_{r^\prime}$ and
$[h_r,f_{r^\prime}]=-(r,\,r^\prime)f_{r^\prime}$;

3) $[e_r,\,f_{r^\prime}]=h_r$ if $r=r^\prime$, and is $0$ if
$r\not=r^\prime$;

4) $(ad~e_r)^{1-2(r,\,r^\prime)/(r,\,r)}e_{r^\prime}=
(ad~f_r)^{1-2(r,\,r^\prime)/(r,\,r)}f_{r^\prime}=0\
\text{if $r\not=r^\prime$ and $(r,r)>0$}$\hfil\hfil
\newline
\phantom{(4)(4)} (equivalently, $r\in {P^\prime}^{re}$);

5) If $(r,\,r^\prime)=0$, then $[e_r,\,e_{r^\prime}]=
[f_r,\,f_{r^\prime}]=0$.

See \cite{3}, \cite{5}, \cite{36}, \cite{40}, \cite{69}
for details.  We remark that for Lie algebras
this definition is equivalent to the definition above using the
generalized Cartan--Borcherds matrix defined by the Gram matrix
of the sequence $P^\prime$.

The algebra $\geg$ is graded by the root lattice $S$ where the
generators $h_r$, $e_r$ and $f_r$ have the weights $0$, $r\in S$
and $-r \in S$ respectively. We have
$$
\geg\ =\ \bigoplus_{\alpha\in S}{\geg_\alpha}=
\geg_0 \bigoplus \left(\bigoplus_
{\alpha\in \Delta_+}\geg_\alpha\right) \bigoplus
\left(\bigoplus_{\alpha\in \Delta_+}
\geg_{-\alpha}\right)
\tag{1.4.13}
$$
where $\geg_0=S\otimes \bc$, and  $\Delta$ is the set of roots
(i. e. $\alpha \in S$ with $\geg_\alpha \not=0$). A root $\alpha$
is called positive ($\alpha \in \Delta_+$) if $(\alpha,\,\M)\le 0$.

For a root $\alpha \in \Delta$ the {\it multiplicity of $\alpha$}
is equal to $\mult(\alpha)=\dim \geg_{\alpha,\0o}-\dim \geg_{\alpha,\1o}$.
Multiplicities $\mult(\alpha)$ of roots and the
numbers $m(a)$ defining generators of $\geg$ are related by the
denominator identity (due to Weyl, Kac, Borcherds)
$$
\split
&\exp{(-2\pi i (\rho,z))}
\prod_{\alpha\in \Delta_+}{\Bigl(1-\exp{\left(-2\pi i (\alpha,z)\right)}
\Bigr)^{\mult(\alpha)}}\\
&=\sum_{w\in W}{\varepsilon (w)
\Bigl(\exp{\left(-2\pi i (w(\rho),z)\right)}\ -\hskip-10pt
\bigr.\sum_{a\in S\cap \br_{++}\M}
{\bigl.m(a)\exp{\left(-2\pi i (w(\rho+a),z)\right)}\Bigr)}}
\endsplit\ ,
\tag{1.4.14}
$$
which identifies multiplicities of factors in \thetag{1.4.9}
and multiplicities of roots of the algebra $\geg$.
See \cite{3}, \cite{5}, \cite{47}, \cite{48}, \cite{36}, \cite{69}.
\smallpagebreak

{\it The generalized Kac--Moody superalgebras $\geg$ above
given by the data (1)---(5) constitute the Theory
of Lorentzian Kac--Moody algebras (or automorphic Lorentzian Kac--Moody
algebras) which we consider.}

\smallpagebreak

By (4), they have similar property to the property (II) from Sect. 1.2 for
finite and affine algebras: their denominator identities give
automorphic forms. For Lorentzian case,
they are automorphic forms on IV type Hermitian symmetric domains.

What is about a similar property to the property (I) from Sect. 1.2
for finite and affine Kac--Moody algebras?
How many data (1)---(5) one may have? This is the main subject of
the paper. We consider that in the next sections.

\subhead
1.5. Finiteness results about hyperbolic root systems (data (1) --- (3))
of Lorentzian Kac--Moody algebras
\endsubhead
To classify finite and affine Kac--Moody algebras, one need to classify
appropriate finite and affine root systems (authors don't know
any other method). They are classified, and
their classification gives classification of finite and affine Lie algebras.

Data (1) --- (5) define Lorentzian Kac--Moody algebras.
One can consider the set of possible data (1) --- (3) from
data (1) --- (5) as {\it hyperbolic root systems which are appropriate
for Lorentzian Kac--Moody algebras.} Main result is that they satisfy
very restricted conditions, and it is possible to classify them, in
principle; this makes the theory of Lorentzian Kac--Moody algebras
similar to theories of finite and affine Kac--Moody algebras.
Moreover, their number is in essential finite when
$\rk S\ge 3$. It means that Lorentzian Kac--Moody algebras for
$\rk S\ge 3$ are {\it exceptional}, like exceptional simple
Lie algebras ${\Bbb E}_6$, ${\Bbb E}_7$, ${\Bbb E}_8$,
${\Bbb F}_4$,  ${\Bbb G}_2$. More generally, for these results,
we can drop the condition (5) considering the set of
possible data (1) --- (3) in data (1) --- (4).

Since we assume that Weyl group is non-trivial,  $\rk S\ge 2$.
If $\rk S=2$, then $P$ has one or two elements,
and classification of data (1) --- (3) in data (1) --- (4) is not difficult.

Further we assume that $\rk S\ge 3$.

From the condition (*) in (4), we have that
$$
W\rtimes Sym(\M)\ has\ finite\ index\ in\ O^+(S).
\tag{1.5.1}
$$
Non-trivial reflection subgroups $W\subset O(S)$
having this property are called {\it of restricted arithmetic type}
(we follow the terminology in \cite{64, Sect. 1.4}). It follows
\cite{64,\ Theorem 1.4.3}  that $W$ has
{\it arithmetic type} which means that
$$
\{\ x\in S\otimes \br\ |\ (x,\,P)\le 0\ \}\
\subset \br_{++}\M\subset \overline{V^+(S)}.
\tag{1.5.2}
$$
In particular, the Weyl vector
$\rho\in \br_{++}\M\subset \overline{V^+(S)}$ and $\rho^2\le 0$.
By \thetag{1.4.2}, the $\rho$ is non-zero.

From \thetag{1.5.2} one can deduce that the set $P$
generates $S\otimes \bq$. From \thetag{1.4.2}, it is easy to see
that the Weyl vector $\rho$ is unique and it is invariant
with respect to $Sym(P_\1o\subset P)$ which has finite index in
$Sym(\M)$. One can
get the same using Fourier expansion \thetag{1.4.5}.
Really, $exp{\left(-2\pi i (\rho,\,z)\right)}$ is one of Fourier
harmonics of the holomorphic automorphic form $\Phi(z)$. It follows
that $\rho^2\le 0$. Otherwise, the automorphic form $\Phi(z)$ has
poles.  The subgroup $H\subset O^+(S)$ of finite index,
in \thetag{1.4.6}, preserves $\pm \Phi(z)$ and the set of
Fourier harmonics
$\varepsilon(w)m(a)\exp{\left(-2\pi i (w(\rho+a), z)\right)}$,
$a\in S\cap \br_{++}\M$, $m(a)\not=0$, $w\in W$.
For Fourier harmonics
$c(x)\exp{\left(-2\pi i (x,\,z)\right)}$ and
$c(y)\exp{\left(-2\pi i (y,\,z)\right)}$, we say that
$x\ge y$, if $x-y\in \overline{V^+(S)}$.
The subgroup $H$ preserves this
ordering and permutes its minimal elements
$\varepsilon(w)exp{\left(-2\pi i (w(\rho),\,z)\right)}$, $w\in W$.
Let $H_\rho\subset H$ be the stabilizer subgroup of $\rho$.
Then $W\rtimes H_\rho=H$ has finite index in $O^+(S)$.
From definition of $\rho$, it follows that
$H_\rho\subset Sym(\M)$. It follows that
$\rho$ is invariant for the subgroup
$H_\rho\subset Sym(\M)$ of finite index.
Remark that we did not use the additional condition (*) in (4) under
this consideration.

Thus, the orbit
$$
Sym(\M)(\rho)\ \ \text{is finite.}
\tag{1.5.3}
$$
Here $\rho\in S\otimes \bq$ and $\rho\not=0$. For some $m\in \bn$,
the element  $r=m\rho\in S$, $r\not=0$ and $Sym(\M)(r)$ is also finite.

Hyperbolic lattices $S$ having the properties similar to
\thetag{1.5.1} and \thetag{1.5.3} for some of their reflection
subgroup $W$ are called reflective.  Here is the exact definition.

\definition
{Definition 1.5.1} A hyperbolic lattice $M$ is called
{\it reflective} if there exists a reflection subgroup $W\subset O(M)$
with a fundamental chamber $\M\subset V^+(M)/\br_{++}$ and
the symmetry group $Sym(\M )$ of the fundamental
chamber $\M$ such that the corresponding semi-direct product
$W\rtimes Sym(\M)$ has finite index in $O(M)$ (i. e. $W$ has
{\it restricted arithmetic type}), and there exists a non-zero $r\in M$
such that the orbit $Sym(\M)(r)$ is finite. The element $r$ is called
a {\it generalized Weyl vector} for the reflection group $W$
and its fundamental chamber $\M$.

Thus, a hyperbolic lattice $M$ is reflective if it has a reflection
group $W\subset O(M)$ having restricted arithmetic type and a
generalized Weyl vector for its fundamental chamber.
The group $W$ has {\it elliptic, parabolic and hyperbolic type},
if it has a generalized Weyl vector
$r$ respectively with $r^2<0$, with $r^2=0$ and no generalized
Weyl vectors with $r^2<0$, and with $r^2>0$ and no generalized
Weyl vectors with $r^2\le 0$.

It is easy to see that a hyperbolic lattice $M$ is reflective if and
only if its full reflection group $W(M)$ (it is generated by reflections in
all roots of $M$) and its fundamental chamber have
a generalized Weyl vector. The group $W(M)$ obviously
has restricted arithmetic type since it is normal in $O(M)$.
\enddefinition

Thus, we have
\proclaim{Proposition 1.5.2} Assume that $\rk S\ge 3$. Then any
data (1) --- (3) from data (1) --- (4) satisfy the conditions:
the Weyl vector $\rho\in S\otimes \bq$ is not zero and
$\rho^2\le 0$;
$$
W\rtimes Sym(\M )\ \text{has finite index in $O(S)$, and\ }
Sym(\M)(\rho)\ \text{is finite.}
\tag{1.5.4}
$$
In particular, the hyperbolic root lattice $S$ is reflective, the Weyl
group $W$ has restricted arithmetic type, the
Weyl vector $\rho$ is a generalized Weyl
vector for $W$ and $\M$;\,
\footnote{Remark that all the statements above are valid without the
additional condition (*) in (4), and, if necessary, one can drop
this additional condition, considering more general class of
Lorentzian Kac--Moody algebras.}
$$
Sym(P_\1o\subset P)\  has\  finite\  index\ in\ Sym(\M ),
\tag{1.5.5}
$$
and
$$
Sym(P_\1o\subset P)(\rho)=\rho.
\tag{1.5.6}
$$
\endproclaim

We have the following general crucial result which, in particular,
gives finiteness for the set of root lattices in (1).

\proclaim{Theorem 1.5.3} The set or reflective hyperbolic lattices
$M$ of $\rk M\ge 3$ is finite (up to isomorphism) if we consider
lattices up to multiplication
of their forms by positive rational numbers.

In particular, the set of hyperbolic root lattices $S$ of $\rk S\ge 3$
for data (1) in data (1) --- (4) is finite if we consider
lattices up to multiplication of their forms by positive rational numbers.
\endproclaim

If rank $\ge 3$ is fixed, see the proof for elliptic case in
\cite{59}, \cite{60}; for parabolic case in \cite{64},
\cite{68}; for hyperbolic case in \cite{66}, \cite{68}.
Boundendess of the rank follows from results in \cite{60},
\cite{63}, \cite{71}.

For fixed rank the proof of Theorem 1.5.3 follows from

\proclaim{Lemma 1.5.4 (about a narrow part of the polyhedron $\M$)}
Let $M$ be a reflective hyperbolic lattice of rank $n=\rk M\ge 3$ with
a reflection group $W$ of restricted arithmetic type and
with a fundamental chamber $\M$ having a generalized Weyl vector.
Let $P(\M)$ be the set of orthogonal vectors to $\M$ directed
outwards of $\M$.

Then there exist $\alpha_1,\,\alpha_2,\, \dots ,\ \alpha_n$ from $P(\M)$ with
the properties:

(a) $\alpha_1,\,\alpha_2,\, \dots,\ \alpha_n$ generate $M\otimes \bq$;

(b) ${4(\alpha_i,\,\alpha_j)^2\over \alpha_i^2\alpha_j^2}<100^2$ for any
$1\le i,\,j \le n$;

(c) the Gram graph of
$\alpha_1,\,\alpha_2,\, \dots ,\ \alpha_n$ is connected, i. e. one
cannot divide the set $\alpha_1,\,\alpha_2,\, \dots ,\ \alpha_n$
in two orthogonal to each other non-empty subsets.
\endproclaim

We remark that for roots $\alpha_1,\,\alpha_2,\, \dots,\ \alpha_n$ the
numbers ${4(\alpha_i,\,\alpha_j)^2\over \alpha_i^2\alpha_j^2}$
are integral. The proof of Lemma 1.5.4 is based on study of geometry of
the convex polyhedra $\M$. Actually, Lemma 1.5.4 is valid for all
convex polyhedra having properties similar to $\M$. The main property is
that the polyhedron $\M$ is convex locally finite and has
{\it almost finite volume.} Here, almost finite volume means the following.
For elliptic case $\M$ is a finite polyhedron of finite volume. For
parabolic case, $\M$ is finite and has finite volume in any angle
with the vertex $\br_{++}r$ and over a compact set on
the horosphere (one has the same property for the Borcherds example
from Sect. 1.3).
For hyperbolic case, $\M$ is finite and has finite volume in any orthogonal
cylinder over compact set in the hyperplane which is orthogonal to $r$.
Here $r$ is a generalized Weyl vector for $\M$.
For example, Lemma 1.5.4 is valid for any finite convex polyhedron
of finite volume in a hyperbolic space of dimension
$n-1\ge 2$. For parabolic and hyperbolic cases, one
has to add also some conditions of periodicity of $\M$.
The proof of Lemma 1.5.4 is not trivial, and it is the most
hard for the hyperbolic case.

\smallpagebreak

Theorem 1.5.3 requires existence of a generalized
Weyl vector $r$ only. If there exists a Weyl vector $\rho$
(satisfying \thetag{1.4.2}), we can
prove finiteness results also for data (1) --- (3) from (1) --- (4).
Two data (1) --- (3) are called isomorphic if they can be identified by
an isomorphism of hyperbolic lattices in (1).
We also accept multiplication of the forms of lattices $S$ in (1)
by positive rational numbers. We know that data (1) --- (3) from
data (1) --- (4) satisfy conditions \thetag{1.5.4} --- \thetag{1.5.6}
of Proposition 1.5.2.
{\it One can consider data (1) --- (3) satisfying additional conditions
\thetag{1.5.4} --- \thetag{1.5.6} as hyperbolic root systems which are
appropriate for Lorentzian Kac--Moody algebras theory.}
From Theorem 1.5.3 and Lemma 1.5.4, one can deduce

\proclaim{Theorem 1.5.5} Assume that $\rk S\ge 3$.

Then the set of data
(1) --- (3) satisfying conditions \thetag{1.5.4} --- \thetag{1.5.6}
is empty if $\rho^2>0$ (i.e. for hyperbolic case); finite if
$\rho^2<0$ (i.e. for elliptic case). It is also finite
if $\rho^2=0$ (i.e. for parabolic case), if one fixes a constant $C>0$
and additionally requires
$[O^+(S)_\rho:Sym(P_{\1o}\subset P)]<C$ where $O^+(S)_\rho$ is
the stabilizer subgroup of $\rho$ in $O^+(S)$.
\endproclaim

See the proof in \cite{64, Theorem 1.3.3}.
We mention that for $\rho^2=0$ the group
$O^+(S)_\rho$ is a $(\rk S-2)$-dimensional
crystallographic group acting in Euclidean affine space related with
the semi-positive lattice $(\rho)^\perp_S$. In particular,
it contains $\bz^{\rk S-2}$ as a subgroup of translations
of finite index. If the constant $C\to +\infty$, number of cases for
the parabolic case tends to infinity, see \cite{64, Example 1.3.4}.
Existence of Weyl vector $\rho$ means geometrically that
faces of $\M$ of highest dimension orthogonal to $\alpha\in P$
are touching a sphere with the center $\br_{++}\rho$ and radius
depending on $\alpha^2$ where $\alpha^2$ is bounded by a constant
depending on $S$. It makes polyhedra $\M$ very special,
and implies finiteness of their set.
\smallpagebreak

As an example, let us prove Theorem 1.5.5 for elliptic case
($\rho^2<0)$ and fixed $n=\rk S\ge 3$.

By Theorem 1.5.3, the set of root lattices $S$ is finite. We
fix one of them.

Let $\alpha$ be a root of $S$. Let
$\alpha=m\alpha_0$ where $\alpha_0$ is a primitive root of $S$ and
$m\in \bn$. Since $\alpha$ is a root,
$({2\alpha\over \alpha^2},\,S)\in \bz$. Thus,
${\alpha_0\over m\alpha_0^2}\in 2^{-1}S^\ast$ where
$S^\ast=\text{Hom}(S,\,\bz)$ is the dual lattice.
It follows that
$m\alpha_0^2\,|\,2\lambda$ where $\lambda$ is the exponent (i.e. the
maximal order of elements) of the
finite Abelian group $S^\ast/S$.
Thus, $m$, $\alpha_0^2$ and $\alpha^2=m^2\alpha_0^2$
are bounded by a constant depending on $S$.

By Lemma 1.5.4, there are roots $\alpha_1, \dots ,\alpha_n$ from $P$
which satisfy conditions (a) and (b) of the lemma. Numbers
${4(\alpha_i,\,\alpha_j)^2\over \alpha_i^2\alpha_j^2}$ are integral
and are less then $100^2$. Since $\alpha_i^2$ are bounded, it follows
that the set of possible Gram matrices
$\left((\alpha_i,\, \alpha_j)\right)$, $1\le i,j\le n$, 
is finite. We fix one of them.

Since $\alpha_1, \dots ,\alpha_n$ generate $S\otimes \bq$ and
the lattice $S$ is non-degenerate, there exists a unique
$\rho\in S\otimes \bq$ such that
$(\rho,\,\alpha_i)=-\alpha_i^2/2$, $i=1,\,\dots,\,n.$
It is the Weyl vector $\rho$ for $P$.

All elements $\alpha\in P$ satisfy the condition
$(\rho,\,\alpha)=-\alpha^2/2$ where $0<\alpha^2<K$
for some constant $K$ depending on $S$. Thus, $P$ is
a subset of the set of elements $a\in S$ satisfying conditions
$-K/2<(\rho,\,a)<0$ and $0<a^2<K$.
If $\rho^2<0$ the set is finite because
the lattice $S$ is hyperbolic
(has exactly one negative square). It follows, that $P$ is finite,
and there exists only a finite number of possibilities for $P$ and its
subset $P_{\1o}$.  It follows Theorem 1.5.5 for elliptic case.

The parabolic case (when $\rho^2=0$ and $\rho\not=0$)
requires additional simple considerations with the crystallographic
group $O(S)_\rho$ and its action on the set $P$ which is infinite
for this case. The group $Sym(P_{\1o}\subset P)$ has finite
index in $O(S)_\rho$ and has finite number of orbits in $P$.

In this proof, formally we did not use the condition (c) of
Lemma 1.5.4, but it is very important for the proof of
Theorem 1.5.3 which also has been used. As we see, for the proof
not only Theorem 1.5.3 is important, but also the method of
its proof (Lemma 1.5.4) is very important.

Let us give an example of classification of data
(1) --- (3) satisfying conditions \thetag{1.5.4} ---
\thetag{1.5.6}.

\proclaim{Theorem 1.5.6} Let us consider all data (1)---(3)
satisfying conditions  \thetag{1.5.4} ---
\thetag{1.5.6} and such that additionally
$\rk S=3$, all roots $\alpha \in P$ have equal squares
$\alpha^2$ (equivalently, the generalized Cartan matrix \thetag{1.4.3}
of $P$ is symmetric), and the Weyl vector $\rho$
has $\rho^2<0$ (elliptic type).

Then the generalized Cartan matrix of $P$ is one of the following
16 symmetric matrices $A_{i,j}$ and $B_j$ below:
$$
A_{1,0}=
\left(\smallmatrix
\hphantom{-}{2}&\hphantom{-}{0}&{-1}\cr
\hphantom{-}{0}&\hphantom{-}{2}&{-2}\cr
{-1}&{-2}&\hphantom{-}{2}\cr
\endsmallmatrix\right),\ \
A_{1,I}=
\left(\smallmatrix
\hphantom{-}{2}&{-2}&{-1}\cr
{-2}&\hphantom{-}{2}&{-1}\cr
{-1}&{-1}&\hphantom{-}{2}\cr
\endsmallmatrix\right),\ \
A_{1,II}=
\left(\smallmatrix
\hphantom{-}{2}&{-2}&{-2}\cr
{-2}&\hphantom{-}{2}&{-2}\cr
{-2}&{-2}&\hphantom{-}{2}\cr
\endsmallmatrix\right),
$$
$$
A_{1,III}=
\left(\smallmatrix
\hphantom{-}{2}&{-2}&{-6}&{-6}&{-2}\cr
{-2}&\hphantom{-}{2}&\hphantom{-}{0}&{-6}&{-7}\cr
{-6}&\hphantom{-}{0}&\hphantom{-}{2}&{-2}&{-6}\cr
{-6}&{-6}&{-2}&\hphantom{-}{2}&{0}\cr
{-2}&{-7}&{-6}&\hphantom{-}{0}&\hphantom{-}{2}\cr
\endsmallmatrix\right);\ \
A_{2,0}=
\left(\smallmatrix
\hphantom{-}{2}&{-2}&{-2}\cr
{-2}&\hphantom{-}{2}&\hphantom{-}{0}\cr
{-2}&\hphantom{-}{0}&\hphantom{-}{2}\cr
\endsmallmatrix\right),\ \
A_{2,I}=
\left(\smallmatrix
\hphantom{-}{2}&{-2}&{-4}&\hphantom{-}{0}\cr
{-2}&\hphantom{-}{2}&\hphantom{-}{0}&{-4}\cr
{-4}&\hphantom{-}{0}&\hphantom{-}{2}&{-2}\cr
\hphantom{-}{0}&{-4}&{-2}&\hphantom{-}{2}\cr
\endsmallmatrix\right),
$$
$$
A_{2,II}=
\left(\smallmatrix
\hphantom{-}{2}&{-2}&{-6}&{-2}\cr
{-2}&\hphantom{-}{2}&{-2}&{-6}\cr
{-6}&{-2}&\hphantom{-}{2}&{-2}\cr
{-2}&{-6}&{-2}&\hphantom{-}{2}\cr
\endsmallmatrix\right),\ \
A_{2,III}=
\left(\smallmatrix
\hphantom{-}{2}&{-2}&{-8}&{-16}&{-18}&{-14}&{-8}&\hphantom{-}{0}\cr
{-2}&\hphantom{-}{2}&\hphantom{-}{0}&{-8}&{-14}&{-18}&{-16}&{-8}\cr
{-8}&\hphantom{-}{0}&\hphantom{-}{2}&{-2}&{-8}&{-16}&{-18}&{-14}\cr
{-16}&{-8}&{-2}&\hphantom{-}{2}&\hphantom{-}{0}&{-8}&{-14}&{-18}\cr
{-18}&{-14}&{-8}&\hphantom{-}{0}&\hphantom{-}{2}&{-2}&{-8}&{-16}\cr
{-14}&{-18}&{-16}&{-8}&{-2}&\hphantom{-}{2}&\hphantom{-}{0}&{-8}\cr
{-8}&{-16}&{-18}&{-14}&{-8}&\hphantom{-}{0}&\hphantom{-}{2}&{-2}\cr
\hphantom{-}{0}&{-8}&{-14}&{-18}&{-16}&{-8}&{-2}&\hphantom{-}{2}\cr
\endsmallmatrix\right);
$$
$$
A_{3,0}=
\left(\smallmatrix
\hphantom{-}{2}&{-2}&{-2}\cr
{-2}&\hphantom{-}{2}&{-1}\cr
{-2}&{-1}&\hphantom{-}{2}\cr
\endsmallmatrix\right),\ \
A_{3,I}=
\left(\smallmatrix
\hphantom{-}{2}&{-2}&{-5}&{-1}\cr
{-2}&\hphantom{-}{2}&{-1}&{-5}\cr
{-5}&{-1}&\hphantom{-}{2}&{-2}\cr
{-1}&{-5}&{-2}&\hphantom{-}{2}\cr
\endsmallmatrix\right),
$$
$$
A_{3,II}=
\left(\smallmatrix
\hphantom{-}{2}&{-2}&{-10}&{-14}&{-10}&{-2}\cr
{-2}&\hphantom{-}{2}&{-2}&{-10}&{-14}&{-10}\cr
{-10}&{-2}&\hphantom{-}{2}&{-2}&{-10}&{-14}\cr
{-14}&{-10}&{-2}&\hphantom{-}{2}&{-2}&{-10}\cr
{-10}&{-14}&{-10}&{-2}&\hphantom{-}{2}&{-2}\cr
{-2}&{-10}&{-14}&{-10}&{-2}&\hphantom{-}{2}\cr
\endsmallmatrix\right),
$$
$$
A_{3,III}=
\left(\smallmatrix
\hphantom{-}{2}&{-2}&{-11}&{-25}&{-37}&{-47}&{-50}
&{-46}&{-37}&{-23}&{-11}&{-1}\cr
{-2}&\hphantom{-}{2}&{-1}&{-11}&{-23}&{-37}&
{-46}&{-50}&{-47}&{-37}&{-25}&{-11}\cr
{-11}&{-1}&\hphantom{-}{2}&{-2}&{-11}&{-25}
&{-37}&{-47}&{-50}&{-46}&{-37}&{-23}\cr
{-25}&{-11}&{-2}&\hphantom{-}{2}&{-1}&{-11}
&{-23}&{-37}&{-46}&{-50}&{-47}&{-37}\cr
{-37}&{-23}&{-11}&{-1}&\hphantom{-}{2}&{-2}
&{-11}&{-25}&{-37}&{-47}&{-50}&{-46}\cr
{-47}&{-37}&{-25}&{-11}&{-2}&\hphantom{-}{2}
&{-1}&{-11}&{-23}&{-37}&{-46}&{-50}\cr
{-50}&{-46}&{-37}&{-23}&{-11}&{-1}&\hphantom{-}{2}
&{-2}&{-11}&{-25}&{-37}&{-47}\cr
{-46}&{-50}&{-47}&{-37}&{-25}&{-11}&{-2}
&\hphantom{-}{2}&{-1}&{-11}&{-23}&{-37}\cr
{-37}&{-47}&{-50}&{-46}&{-37}&{-23}&{-11}&{-1}
&\hphantom{-}{2}&{-2}&{-11}&{-25}\cr
{-23}&{-37}&{-46}&{-50}&{-47}&{-37}
&{-25}&{-11}&{-2}&\hphantom{-}{2}&{-1}&{-11}\cr
{-11}&{-25}&{-37}&{-47}&{-50}&{-46}&{-37}
&{-23}&{-11}&{-1}&\hphantom{-}{2}&{-2}\cr
{-1}&{-11}&{-23}&{-37}&{-46}&{-50}
&{-47}&{-37}&{-25}&{-11}&{-2}&\hphantom{-}{2}\cr
\endsmallmatrix\right),
$$
$$
B_1=
\left(\smallmatrix
\hphantom{-}2 & \hphantom{-}0 &-3 & -1 \cr
\hphantom{-}0 & \hphantom{-}2 &-1 & -3 \cr
-3&-1 & \hphantom{-}2 &  \hphantom{-}0 \cr
-1&-3 & \hphantom{-}0   &  \hphantom{-}2 \cr
\endsmallmatrix\right),
\hskip30pt
B_2=
\left(\smallmatrix
\hphantom{-}2  & -1 & -4 & -1\cr
-1 &  \hphantom{-}2 & -1 & -4\cr
-4 & -1 &  \hphantom{-}2 & -1\cr
-1 & -4 & -1 &  \hphantom{-}2
\endsmallmatrix\right),
$$
$$
B_3=
\left(\smallmatrix
\hphantom{-}2 & \hphantom{-}0 & -4 & -6 & -4 &  \hphantom{-}0 \cr
\hphantom{-}0 & \hphantom{-}2 &  \hphantom{-}0 & -4 & -6 & -4 \cr
-4& \hphantom{-}0 &  \hphantom{-}2 &  \hphantom{-}0 & -4 & -6 \cr
-6&-4 &  \hphantom{-}0 &  \hphantom{-}2 &  \hphantom{-}0 & -4 \cr
-4&-6 & -4 &  \hphantom{-}0 &  \hphantom{-}2 &  \hphantom{-}0 \cr
\hphantom{-}0 &-4 & -6 & -4 &  \hphantom{-}0 &  \hphantom{-}2 \cr
\endsmallmatrix\right),
\hskip30pt
B_4=
\left(\smallmatrix
\hphantom{-}2  & -1 & -7 & -10 & -7 & -1 \cr
-1 &  \hphantom{-}2 & -1 & -7  & -10& -7 \cr
-7 & -1 &  \hphantom{-}2 & -1  & -7 & -10\cr
-10& -7 & -1 &  \hphantom{-}2  & -1 & -7 \cr
-7 & -10& -7 & -1  &  \hphantom{-}2 & -1 \cr
-1 & -7 & -10& -7  & -1 &  \hphantom{-}2 \cr
\endsmallmatrix\right).
$$
For all these cases the fundamental chamber $\M$ is a closed
polygon on the hyperbolic plane with angles respectively:
\newline
$A_{1,0}:\ \pi/2,\,0,\,\pi/3;$\ \
$A_{1,I}:\ 0,\,\pi/3,\,\pi/3;$\ \
$A_{1,II}:\ 0,\,0,\,0;$\ \
$A_{1,III}:\ 0,\,\pi/2,\,0,\,\pi/2,0;$
\newline
$A_{2,0}:\ 0,\,\pi/2,\,0;$\ \
$A_{2,I}:\ 0,\,\pi/2,\,0,\,\pi/2;$\ \
$A_{2,II}:\ 0,\,0,\,0,\,0;$\ \
\newline
$A_{2,III}:\ 0,\,\pi/2,\,0,\,\pi/2,\,0,\,\pi/2,\,0,\,\pi/2;$
\newline
$A_{3,0}:\ 0,\,\pi/3,\,0;$\ \
$A_{3,I}:\ 0,\,\pi/3,\,0,\,\pi/3;$\ \
$A_{3,II}:\ 0,\,0,\,0,\,0,\,0,\,0;$\ \
\newline
$A_{3,III}: \ 0,\,\pi/3,\,0,\,\pi/3,\,0,\,\pi/3,\,0,\,\pi/3,\,0,\,
\pi/3,\,0,\,\pi/3$.
\newline
$B_1: \ \pi/2,\,\pi/3,\,\pi/2,\,\pi/3;$\ \
$B_2: \ \pi/3,\,\pi/3,\,\pi/3,\,\pi/3;$
\newline
$B_3: \ \pi/2,\,\pi/2,\,\pi/2,\,\pi/2,\,\pi/2,\,\pi/2;$\ \
$B_4: \ \pi/3,\,\pi/3,\,\pi/3,\,\pi/3,\,\pi/3,\,\pi/3.$
\newline
All these polygons are touching a circle with the center $\br_{++}\rho$
where $\rho$ is the Weyl vector. The first 12 matrices $A_{i,j}$ give
non-compact polygons (they have at least one zero angle).
The last four matrices $B_i$ give compact polygons.
\endproclaim

See the proof in \cite{40, Theorems 1.2.1 and 1.3.1}. The proof is
based on Lemma 1.5.4 and uses computer calculations. They
follow the proof of Theorem 1.5.5 which we outlined above.

Theorem 1.5.6 gives classification of generalized Cartan matrices of
the possible sets $P$, but one can get from here description
of the corresponding data (1) --- (3) as follows.
Let $A=(a_{ij})$, $1\le i,j\le m$, be one of generalized
Cartan matrices of Theorem 1.5.6. Let us consider a free $\bz$-module
$\widetilde{M}=
\oplus_{i=1}^{m}{\bz \widetilde{\alpha}_i}$
with the symmetric bilinear form
$\left((\widetilde{\alpha}_i,\,\widetilde{\alpha}_j)\right)=A$.
The matrix $A$ has rank three,
and $\widetilde{M}$ modulo the kernel of this form defines a
hyperbolic lattice $M$ of the rank three.
It is generated by images $\alpha_i$ of $\tilde{\alpha}_i$.
Any its integral overlattice $M\subset S$ of finite
index  can be taken for the datum (1). Their number is finite since
$S\subset M^\ast$. The set $P$ is given by all elements $\alpha_i$.
Since $\alpha_i^2=2$, all elements $\alpha_i$ are roots of $S$.
Reflections in roots of $P$ generate the reflection group $W\subset O(S)$.
Its fundamental chamber is equal to
$\M=\{0\not=x\in S\otimes \br\ |\ (x,\ P)\le 0\ \}/\br_{++}\subset
\overline{\La(S)}=\overline{V^+(S)/\br_{++}}$.  As a subset
$P_\1o\subset P$ one can take any number of roots
$\alpha_i\in P$ which satisfy
the condition $\alpha_i^2=2\,|\,(\alpha_i,\,S)$.
One can check (it depends on the generalized Cartan matrix
$A$ only) that all these data satisfy conditions (1) --- (3) and
conditions \thetag{1.5.4} --- \thetag{1.5.6}.

Almost for all  matrices $A_{i,j}$ of Theorem 1.5.6 and corresponding
data (1) --- (3)  one can construct additional data
(4) satisfying (5). See \cite{36}, \cite{37}, \cite{39},
\cite{41}. We give almost all these examples in Sect. 2.6 below.

\remark{Remark 1.5.7} Formally, for the considered theory of Lorentzian
Kac--Moody algebras which are given by data (1)---(5),
it is sufficient to consider only  reflective
hyperbolic lattices $S$ with a Weyl vector $\rho$.
In particular, $\rho^2\le 0$, if $\rk S \ge 3$, and all these lattices
$S$ are reflective of elliptic or parabolic type.

There are several reasons why it is important and necessary
to consider arbitrary reflective hyperbolic lattices with
arbitrary generalized Weyl vector (having
a square of any sign, in particular).

First, finiteness results are valid for general reflective
hyperbolic lattices. For the classification of reflective hyperbolic
lattices with a Weyl vector in general, it is necessary
first to find maximal reflective hyperbolic lattices with a
generalized Weyl vector, and then to find all their
sublattices of finite index having a Weyl vector. In
practice, it is impossible to consider reflective hyperbolic lattices
with a Weyl vector and with a generalized Weyl
vector separately.

Second, it seems, there is a more general class of Lie algebras
(analogous to Lorentzian Kac--Moody algebras which we consider
here) such that for this class it is necessary to consider
reflective hyperbolic lattices and identities similar to
\thetag{1.4.14} with a generalized Weyl vector $\rho$ having
square with any sign.
In \S 2 we shall give and classify many such identities. Very
few of them are related with a Weyl vector $\rho$. It may be
similar to McDonald identities which were first discovered and
later found to be related with affine Kac--Moody algebras.
Results of \cite{26} and \cite{50} give this hope.

From our point of view, all reflective hyperbolic lattices and
identities analogous to \thetag{1.4.14} with a generalized
Weyl vector $\rho$, should be interesting for one or another
theory of Lie algebras which is analogous to the theory of
Lorentzian Kac--Moody algebras which we consider here. It is
a very interesting problem to understand their importance from
Lie algebras point of view.
\endremark

\subhead
1.6. Finiteness conjectures and results about data (4), (5)
for Lorentzian Kac--Moody algebras
\endsubhead
Here we follow \cite{40} and \cite{65}.

We expect that data (4), (5) for Lorentzian Kac--Moody algebras
are also very restricted, but we don't have here so strong results
as for hyperbolic root systems (data (1) --- (3)).

\definition{Definition 1.6.1}
A lattice $Q$ with two negative squares is called {\it reflective},
if its Hermitian symmetric domain $\Omega(Q)$ has a
meromorphic (not necessarily holomorphic) automorphic form
$\Phi$ of non-zero weight with respect to $G\subset O^+(Q)$ of finite
index such that its divisor is union of rational quadratic divisors
which are orthogonal to some roots of $Q$.
The automorphic form $\Phi$ is then also called {\it reflective} for
$Q$. (Here, for further considerations, we assume that a reflective
automorphic form $\Phi$
has non-zero weight, but might be this condition can be weakened
by the condition that $\Phi$ is not a constant.)

Obviously, reflectivity of $Q$ does not change by multiplication of
the form of the lattice $Q$ by rational numbers.
\enddefinition

We suggested in \cite{40, Conjecture 2.2.1}) and \cite{65}

\proclaim{Conjecture 1.6.2} The set of reflective lattices $Q$
with two negative squares and $\rk Q\ge 5$ is finite up to
multiplication of forms of lattices $Q$ by positive rational
numbers.
\endproclaim

We expect the statement of Conjecture 1.6.2 because of
{\it Koecher principle} (e. g. see \cite{1}):
{\it Any non-constant meromorphic automorphic form on
a Hermitian symmetric domain $\Omega$ should have
non-trivial divisor in $\Omega$ if
$\dim \Omega - \dim \Omega_\infty\ge 2$.}
Here $\Omega_\infty$ is the
set of points at infinity of $\Omega$ which is added to get Satake
compactification
$G\backslash \Omega\subset G\backslash (\Omega\cup \Omega_\infty)$
of the arithmetic quotient $G\backslash \Omega$.

For a lattice $Q$ with two negative squares and
$G\subset O(Q)$ of a finite index,
$\dim \Omega (Q)=\rk Q-2$ and $\dim \Omega(Q)_\infty=s(Q)-1$
where $s(Q)$ is the rank of a maximal isotropic sublattice of $Q$.
We have $0\le s(Q) \le 2$. In particular, the Koecher principle
is valid for $\Omega(Q)$ if $\rk Q\ge 5$ or $\rk Q\ge 3+s(Q)$.

We apply Koecher principle to restrictions $\Phi|\Omega(Q_1)$ of a
reflective automorphic form $\Phi$ on all subdomains
$\Omega(Q_1)\subset \Omega (Q)$ where $Q_1\subset Q$ is a sublattice
of $Q$ with two negative squares and $\rk Q_1\ge 3+s(Q_1)$. Here
$\Phi | \Omega(Q_1)$ is an automorphic form on the Hermitian
symmetric domain $\Omega(Q_1)$ of the same weight as
$\Phi$, and it is not a constant, if it is not zero, since its weight
is not $0$. Thus, we get

\proclaim{Proposition 1.6.3} Suppose that $Q$ is a reflective
lattice with two negative squares. Then
$$
\Omega(Q_1)\bigcap
\left(\bigcup_{\text{root\ }\alpha\in Q}{\Omega(Q)_\alpha}\right)
\not=\varnothing
\tag{1.6.1}
$$
for any sublattice $Q_1\subset Q$ such that $Q_1$ has two negative
squares and $\rk Q_1\ge 3+s(Q_1)$. Here $\Omega(Q)_\alpha$ is the rational
quadratic divisor which is orthogonal to a root
$\alpha\in Q$.
\endproclaim

Below we give an example (from \cite{65}) which shows that this
condition is very strong.

For a lattice $Q$ any element $\alpha\in Q$ with
$\alpha^2=2$ is a root. If for the definition of reflective
lattices $Q$ and reflective automorphic forms $\Phi$
we shall consider only roots with square 2, we shall get {\it 2-reflective
lattices $Q$} and {\it 2-reflective automorphic forms $\Phi$.}
This is a special case of reflective lattices and reflective
automorphic forms.

We consider lattices
$$
T_n=H\oplus H\oplus E_8\oplus E_8\oplus \langle 2n \rangle
\tag{1.6.2}
$$
where $n\in \bn$. Here $E_8$ is an even unimodular positive definite
lattice of the rank 8. A lattice $\langle A \rangle$ is the lattice
with the matrix $A$ for some basis. We want to show that
the lattices $T=T_n$ are not 2-reflective for big $n$.

This example is interesting
because arithmetic quotients $G\backslash \Omega(T_n)$ by subgroups
$G\subset O^+(T_n)$ of finite index give moduli of K3 surfaces
of degree $2n$, and
$$
{Discr}=\bigcup_{\alpha\in T_n \text{\ with\ }\alpha^2=2}
{\Omega(T_n)_\alpha}
\tag{1.6.3}
$$
is discriminant of the moduli. Points of
$Discr$ give K3 surfaces with singularities. Thus,
we want to show that the discriminant of K3 surfaces moduli cannot
be given as a divisor of an automorphic form of non-zero weight
if the degree of K3 surfaces is high enough.

Let us consider an even unimodular lattice
$L=3H\oplus E_8\oplus E_8$. Consider a primitive element $h\in L$ with
$h^2=-2n$. Using standard results about indefinite lattices, and
discriminant forms technique (see \cite{57}), one can prove that
the orthogonal complement $h^\perp$ to the element $h$ in $L$
is isomorphic to $T_n$.

We use the following general construction. Assume that $K\subset L$
is a primitive sublattice with two negative squares, $\rk K\ge 3+s(K)$
where $s(K)$ is the rank of the maximal isotropic sublattice of $K$,
and $K$ does not have elements with square $2$. Let $S=K^\perp$
be the orthogonal complement to $K$ in $L$. The lattice $S$ is
hyperbolic since $L$ has exactly three negative squares.

Let us consider the set $\Delta \subset S^\ast$ with the following
properties: if $\delta_1\in \Delta$, then (i) $\delta_1^2>0$;
(ii) there exists a $\delta_2\in K^\ast$ such that either $\delta_2=0$ or
$\delta_2^2>0$, and $\delta_1+\delta_2\in L$;
(iii) $\delta_1^2+\delta_2^2=2$.
In particular, $0<\delta_1^2\le 2$ and $0\le \delta_2^2<2$.

We have a simple

\proclaim{Lemma 1.6.4}Let $h\in S$ is primitive, $h^2=-2n$ and
the lattice $h^\perp \simeq T_n$ is 2-reflective. Then there exists
$\delta\in \Delta\subset S^\ast$ such that $h\in \delta^\perp$.
\endproclaim

\demo{Proof} Assume that $T_n=h^\perp$ is reflective. We apply
Proposition 1.6.3 to $Q=T_n$, $Q_1=K$ and roots with square two.
We get that there exists $\delta\in h^\perp$ such that $\delta^2=2$ and
$\Omega(h^\perp)_\delta\cap \Omega(K)\not=\emptyset$.
We have $\delta=\delta_1+\delta_2$ where $\delta_1\in S^\ast$
and $\delta_2\in K^\ast$. It follows that $\delta_1^2+\delta_2^2=2$.
Since $(h,\,\delta)=0$ and $h\in S$, it follows that $(h,\,\delta_1)=0$.
The lattice $S$ is hyperbolic and $h^2<0$. It follows that either
$\delta_1^2>0$ or $\delta_1=0$. The last case is impossible because
then $\delta=\delta_2\in K$ and $\delta^2=2$. But we assume that $K$
does not have elements with square 2. Thus, $\delta_1^2>0$.
If $\delta_2=0$, then $\delta_1\in \Delta$ and $(h,\,\delta_1)=0$ as
we want. Let $\delta_2\not=0$.
Let $\bc\omega \in \Omega(h^\perp)_\delta\cap \Omega(K)$.
Then $\omega \in K\otimes \bc$, $(\omega,\,\omega)=0$,
$(\omega,\overline{\omega})<0$ and $(\omega,\,\delta_2)=0$.
Writing $\omega=a+bi$ where $a,\,b \in K\otimes \br$, we then get
$a^2=b^2<0$ and $(a,\,b)=(a,\,\delta_2)=(b,\,\delta_2)=0$.
Since $K$ has exactly
two negative squares, it follows that $\delta_2^2>0$.
It proves the statement.
\enddemo

We can interpret the statement of Lemma 1.6.4 geometrically as follows.
Let $\La(S)=V^+(S)/\br_{++}$ be the hyperbolic space related with the
hyperbolic lattice $S$.
Each element $\delta \in \Delta$ defines a hyperplane
$\Ha_\delta\subset  \La(S)$ orthogonal to $\delta$.
Elements $\delta \in \Delta$ have $\delta^2<2$ and $m\delta \in S$
where $m$ is the exponent of $S^\ast/S$. It follows that
{\it the set of hyperplanes $\Ha_\delta$, $\delta\in \Delta$, is locally
finite in $\La(S)$.} The lattice $h^\perp\simeq T_n$ is not reflective
if the point $\br_{++}h\in \La(S)$ does not belong to this locally
finite set of hyperplanes (one should change $h$ by $-h$ if it is necessary).
The set of points $\br_{++}h\in \La (S)$, $h\in S$, is
everywhere dense in the hyperbolic space $\La (S)$. There are plenty of
these points which do not belong to the locally finite set of hyperplanes
$\Ha_\delta$, $\delta\in \Delta$, and define then non-reflective lattices
$h^\perp\simeq T_n$. For example, it follows that there exists an
infinite sequence of integers $n$ such that the lattice $T_n$ is not
2-reflective. Exactly that had been demonstrated in \cite{65}.

Let us take a concrete lattice $K=H(2)\oplus H(2)\oplus \langle 4 \rangle$
where $M(k)$ is a lattice which is obtained from a lattice $M$ by
multiplication
of the form of the lattice $M$ by $k\in \bq$.
Then $S=\langle -4 \rangle \oplus D_8\oplus D_8$ where $D_8$
is the root lattice of the root system ${\Bbb D}_8$.
Thus the lattice $S$ is the set of integral vectors
$h=(x,u_1, \dots ,u_8,v_1, \dots ,v_8)$ such that $u_1+\cdots +u_8\equiv 0 \mod 2$ and
$v_1+\cdots +v_8\equiv 0 \mod 2$.  The form is given by
$-4xx^\prime +u_1u_1^\prime+ \cdots + u_8u_8^\prime+
v_1v_1^\prime+ \cdots + v_8v_8^\prime$.

We consider $h_0=(1,0,0, \dots ,0)$ and
$\Delta_0=\{\delta \in \Delta\ |\ (\delta,\,h_0)=0\}$ (geometrically,
it is the set of hyperplanes $H_\delta$, $\delta \in \Delta$, which
contain the point $\br_{++}h_0$).  We consider
the set of primitive $h\in S$ such that $h^2<0$,
$\br_{++}h\not\in H_\delta$, if $\delta \in \Delta_0$, and
the distance between points $\br_{++}h_0$ and $\br_{++}h$ of the
hyperbolic space $\La(S)$ is small enough for the point $\br_{++}h$
would not belong to other hyperplanes $H_\delta$, $\delta\in \Delta$.
By Lemma 1.6.4, for these $h$ the lattice $T_n=h^\perp$ is not 2-reflective.
As the result we get

\proclaim{Theorem 1.6.5}
Let us consider integers $y$ of the form
$$
y=y_1^2+y_2^2+y_3^2+y_4^2+y_5^2+y_6^2+y_7^2+y_8^2
\tag{1.6.4}
$$
where all $y_i$ are natural,
$$
y_1+y_2+y_3+y_4+y_5+y_6+y_7+y_8\equiv 0\mod 2,
$$
and
$$
0<y_1<y_2<y_3<y_4<y_5<y_6<y_7<y_8.
$$
Let $u$ and $v$ are two numbers of the form
\thetag{1.6.4} for the corresponding vectors
$(u_1,\dots ,u_8)$, $(v_1, \dots ,v_8)$,  the
vector $(x,u_1, \dots ,u_8,v_1,\dots ,v_8)$ is primitive (it is sufficient
to suppose that it is primitive in the sublattice
$u_1+\cdots + u_8\equiv v_1+\cdots + v_8\equiv 0\mod 2$) and
$x^2>(9/4)(u+v)$. Then for
$$
2n=4x^2-u-v
\tag{1.6.5}
$$
the lattice $T_n$ is not 2-reflective.
\endproclaim

Any sufficiently large even integer $y>N$ can be represented in the form
\thetag{1.6.4} for some primitive $(y_1, \dots ,y_8)$.
It follows that any sufficiently large even $2n$
can be represented in the form \thetag{1.6.5},
and the lattice $T_n$ is not 2-reflective. More exactly, elementary
estimates show that it is true for
$$
n> \left({32\over 3}+\sqrt{128+8N}\right)^2,
$$
and lattices $T_n$ are not 2-reflective for these $n$.

\subhead
1.7. An example of classification of Lorentzian Kac--Moody algebras or
the rank 3
\endsubhead
Below and in what follows we consider hyperbolic lattices
$S_t=H\oplus \langle 2t\rangle$ and lattices with
two negative squares
$L_t=H\oplus S_t=2H\oplus\langle 2t\rangle$. We denote
$$
\widehat{O}^+(L_t)=\{g\in O^+(L_t)\ |\ g \text{\ is trivial on
$L_t^\ast/L_t$}\}.
\tag{1.7.1}
$$
The group $\widehat{O}^+(L_t)$ is called {\it extended paramodular
group}.
In \S 2 below we give classification of
Lorentzian Kac--Moody algebras with the root lattice
$S_t^\ast$, symmetry lattice $L_t^\ast$ and the
symmetry group $\widehat{O}^+(L_t)$.

The lattices
$S_t^\ast=H\oplus \langle {1\over 2t}\rangle$ and
$L_t^\ast=2H\oplus \langle {1\over 2t}\rangle$ will be integral
after multiplication of their forms by $2t$. Automorphism groups
$O(L_t)=O(L_t^\ast)$ of lattices $L_t$ and $L_t^\ast$
are naturally identified. Thus, we can consider the group
$\widehat{O}^+(L_t)$ as a subgroup of $O(L_t^\ast)$.

This case is especially interesting since the lattices $S_t$ and
$L_t$ are maximal even if $t$ is square-free. It follows that many
hyperbolic even lattices $S$ of the rank three and many
even lattices $L$ of
the rank five with two negative squares have equivariant embeddings
to $S_t$ and $L_t$ respectively. Here embedding of lattices $M_1\subset M$
of the same rank is called {\it equivariant} if it induces an embedding
$O(M_1)\subset O(M)$ of their automorphism groups. Any lattice has
an equivariant embedding to a maximal one.
Thus, studying lattices $S_t$ and $L_t$, we at the same time study
Lorentzian Kac--Moody algebras with root lattices $S^\ast$
and symmetry groups $G\subset O^+(L)$ (of finite index)
where $S$ has an equivariant embedding to $S_t$,
and $L$ has an equivariant embedding to $L_t$.
See \thetag{2.2.7} below.

More generally, the same will be true for $m$-dual lattices of
lattices $S_t$ and $L_t$. Here, for a lattice $M$ and a square-free
$m\in \bn$ the $m$-dual lattice of $M$ is
$$
M^{\ast,m}= \left(\bigcap_{p|m}{\left(M\otimes \bq\cap (M\otimes \bz_p)^\ast\right)}\right) \bigcap 
\left(\bigcap_{p\nmid m}{\left(M\otimes \bq\cap 
M\otimes \bz_p\right)}\right) .
\tag{1.7.2}
$$

Another importance of this example is its relation with
the theory of Abelian surfaces $A$ over $\bc$ with polarization of
type $(1,\,t)$. We remind that it means an algebraic integral
2-dimensional cohomology class of $A$ given by a
symplectic integral form
$$J_t=
\pmatrix
0 &0 & 1&0\\
0 &0 & 0&t\\
-1&0 & 0&0\\
0 &-t& 0&0
\endpmatrix
\tag{1.7.3}
$$
in some basis of $H_1(A,\,\bz)$.
The lattice $\langle 2t \rangle$ is the Neron-Severi lattice of
general Abelian surfaces with
polarization of the type $(1,\,t)$.
The lattice $L_t(-1)$ is the lattice of transcendental
cycles (2-dimensional) of a general
Abelian surface with polarization of the type
$(1,\,t)$. The arithmetic quotient
$\widehat{O}^+(L_t)\backslash \Omega(L_t)$
gives moduli space of Abelian surfaces with polarization of the type
$(1,\,t)$ when one identifies an Abelian surface $A$ with its dual
$\widehat{A}$. See \cite{35} and \cite{41} for details.
Thus, all automorphic forms with respect to the extended
paramodular group $\widehat{O}^+(L_t)$ have some geometric
interpretation for the moduli space of Abelian surfaces.

\head
2. Classification of Lorentzian Kac--Moody algebras with the hyperbolic
root lattice $S_t^\ast$, the symmetry lattice $L_t^\ast$,
and the symmetry group $\widehat{O}^+(L_t)$.
\endhead

Here we give classification (we follow \cite{43}) of
Lorentzian Kac--Moody algebras $\geg$ with the root lattice
$S_t^\ast$, symmetry lattice $L_t^\ast$ and the
symmetry group $\widehat{O}^+(L_t)$. Here $t$ is any natural number.
See definitions in Sect. 1.7. Perhaps, this is the first case when
a large class of Lorentzian Kac--Moody algebras is classified.
We try to remind main definitions and notations from \S 1
to make reading easier.

\smallpagebreak

\subhead
2.1. The formulation of the classification result about
Lorentzian Kac--Moody algebras
\endsubhead
By Sect. 1.4, a Lorentzian Kac--Moody algebra $\geg$ with the root
lattice $S_t^\ast$, the symmetry lattice $L_t^\ast$ and
the symmetry group $\widehat{O}^+(L_t)$ is given by
a holomorphic automorphic form $\Phi (z)$,
$z\in \Omega(V^+(S_t))=S_t\otimes \br+iV^+(S_t)$,
with respect to the
group $\widehat{O}^+(L_t)$ with the Fourier expansion
$$
\Phi(z)=\sum_{w\in W}{\varepsilon(w)
\Bigl(\exp{\left(-2\pi i (w(\rho),z)\right)}\ -\hskip-10pt
\bigr.\sum_{a\in S_t^\ast \cap \br_{++}\M}
{\bigl.m(a)\exp{\left(-2\pi i (w(\rho+a),z)\right)}\Bigr)}}
\tag{2.1.1}
$$
where all coefficients $m(a)$ should be integral;
$W\subset O^+(S_t)$ is a reflection subgroup (the Weyl group
of the algebra) generated by reflections in some roots of $S_t$; the
$\varepsilon:W\to \{\pm 1\}$ is its quadratic character;
$\M\subset \La (S_t)=V^+(S_t)/\br_{++}$ is
a fundamental chamber of $W$;
$\rho\in S_t\otimes \bq$ is the Weyl vector (see
\thetag{1.4.2}) for the set $P(\M)\subset S_t^\ast$
of orthogonal roots to $\M$
(it is the set of simple real roots of the algebra).
Additionally, the automorphic form $\Phi(z)$ should be reflective,
i. e. it should have zeros only in
rational quadratic divisors which are orthogonal to roots of
$L_t$. A minor additional condition from Sect. 1.4 is that
the semi-direct product $W\rtimes Sym(P(\M)_\1o\subset P(\M))$
should have finite index in $O(S_t)$.
Here
$P(\M)_\1o\subset P(\M)$ is the set of odd simple real roots
defined by the condition that $\varepsilon(s_\alpha)=1$ for the
reflection $s_\alpha$ with respect to a root
$\alpha\in P(\M)_\1o$. The set $P(\M)_\0o =P(\M)-P(\M)_\1o$ is
the set of even simple real roots.

The automorphic form $\Phi$ automatically has an infinite product expansion
$$
\Phi(z)=\exp{(-2\pi i (\rho,z))}
\prod_{\alpha\in \Delta_+}{\Bigl(1-\exp{\left(-2\pi i (\alpha,z)\right)}
\Bigr)^{\mult(\alpha)}}
\tag{2.1.2}
$$
where $\mult(\alpha)\in \bz$ are multiplicities of roots of $\geg$
and $\Delta_+\subset S_t^\ast $ is the set of positive roots of
the algebra $\geg$ defined by the condition
$(\Delta_+,\,\M)\le 0$.
The infinite product \thetag{2.1.2} can be also used to define the
automorphic form $\Phi$ and the Kac--Moody algebra $\geg$. The
identity \thetag{2.1.1}=\thetag{2.1.2} is called the denominator
identity.

The infinite sum part \thetag{2.1.1} of the denominator identity
defines generators and defining  relations of the Lorentzian Kac--Moody
algebra $\geg$, which is a generalized Kac--Moody (or Borcherds)
superalgebra. The algebra $\geg$ is graded by the lattice
$S_t^\ast$
$$
\geg=\bigoplus_{\alpha\in S_t^\ast}{\geg_\alpha}=
\geg_0\bigoplus \left(\bigoplus_
{\alpha\in \Delta_+}\geg_\alpha\right) \bigoplus
\left(\bigoplus_{\alpha\in \Delta_+}
\geg_{-\alpha}\right),\ \ \ \geg_0=S_t\otimes \bc\,,
\tag{2.1.3}
$$
and the product part \thetag{2.1.2} of the denominator identity
gives multiplicities
$$
\mult(\alpha):= \dim \geg_\alpha=\dim \geg_{\alpha,\0o}-
\dim \geg_{\alpha,\1o}
\tag{2.1.4}
$$
of root-spaces $\geg_\alpha$, $\alpha \in \Delta_+$,
of the superalgebra $\geg$. We also have
$\mult(-\alpha)=\mult(\alpha)$.

Our classification result is

\proclaim{Theorem 2.1.1} There are exactly 29 Lorentzian
Kac--Moody algebras $\geg$ with the root lattice
$S_t^\ast$, the symmetry lattice $L_t^\ast$ and the symmetry group
$\widehat{O}^+(L_t)$ (equivalently, there are exactly
29 automorphic forms \thetag{2.1.1}) for all natural $t\in \bn$.
They exist only for
$$
\split
t=&1\,(three),\ 2\,(seven),\ 3\,(seven),\ 4\,(seven),\ 8\,(one),\
9\,(one),\\
&12\,(one),\ 16\,(one),\ 36\,(one).
\endsplit
$$
where for each $t$ we show in brackets the number of forms.
All these 29 forms are given in Sect. 2.6 below.
\endproclaim

The proof of Theorem 2.1.1 is divided in two parts. First, one
should construct all 29 automorphic forms of the theorem. Second,
one should prove that there are no other 
automorphic forms which satisfy its conditions. We consider
construction of the 29 automorphic forms in Sect. 2.2. The completeness
of their list is considered in Sects 2.2 --- 2.4.

All 29 automorphic forms of Theorem 2.1.1 were found in
\cite{36}, \cite{37}, \cite{39}, \cite{41} together with
their infinite sum and product expansions.
Formally, the forms for $t=8,\,12,\,16$ are new but they
coincide with some forms for $t=2,\,3,\,4$ respectively
after appropriate change of variables.
Thus, in fact, they are not new.

We describe all 29 automorphic forms
of Theorem 2.1.1 in Sect. 2.6 below.
To construct these automorphic forms, we used in \cite{41}
a variant of the Borcherds exponential lifting \cite{7}
which we applied to Jacobi modular forms, see \cite{41, Theorem 2.1}.
It gives the infinite product
expansion \thetag{2.1.2} of the forms, and construction of these forms
using the infinite product expansion. We consider this variant
in the next Sect. 2.2.

To find Fourier expansion of
the 29 automorphic forms, we have used in
\cite{36}, \cite{37}, \cite{39}, \cite{41} several different methods,
there are no general methods for that. One of these methods is
{\it arithmetic lifting} from \cite{29} --- \cite{32}
of Jacobi forms which gives simple Fourier
expansions \thetag{2.1.1} of some of the 29 automorphic forms.
In this paper we don't discuss the methods of construction of
the Fourier expansions of the 29 automorphic forms, and only give
the corresponding formulae. We don't need these Fourier expansions for
the proof of Theorem 2.1.1, but they are very important for the
construction of the algebras $\geg$ using generators and defining
relations.

\subhead
2.2. A variant for Jacobi forms of Borcherds automorphic products
\endsubhead
We use a general result from \cite{41} which permits to construct
as products similar to \thetag{2.1.2} many automorphic forms with
respect to the extended paramodular group.

We use a basis $f_2,\,f_{-2}$ for $H$ with the Gram matrix
$\pmatrix 0&-1\\-1&0\endpmatrix$ and $f_3$ for
$\langle 2t\rangle$. Together they give a basis $f_2,\,f_3,\,f_{-2}$
for the lattice $S_t=H\oplus \langle 2t \rangle$ with the Gram matrix
$\pmatrix 0&0&-1\\0&2t&0\\-1&0&0\endpmatrix$.
The dual lattice $S_t^\ast$ has the basis
$f_2,\,\widehat{f_3}=f_3/2t,\,f_{-2}$.
We denote
$\alpha=(n,l,m):=nf_2-l\widehat{f_3}+mf_{-2} \in S_t^\ast$,
where $\alpha^2=-2nm+{l^2\over 2t}$. We denote by
$D(\alpha)=2t\alpha^2=-4tnm+l^2$ the {\it discriminant} or
{\it norm} of $\alpha$. For the dual lattice $S_t^\ast$ we usually
use this form. Thus, the discriminant or the norm gives the
integral lattice $S_t^\ast(2t)$. Let
$z=z_3f_2+z_2f_3+z_1f_{-2} \in \Omega(V^+(S_t))$. Then
$\exp{\left(-2\pi i (\alpha,\,z)\right)}=q^nr^ls^m$
where $q=exp{(2\pi i z_1)}$, $r=\exp{(2\pi i z_2)}$,
$s=\exp{(2\pi i z_3)}$.

In \cite{41} a variant of the Borcherds exponential lifting was proved.
Borcherds exponential lifting \cite{7} gives lifting of modular forms
of one variable. In \cite{41} we constructed its variant for
Jacobi modular forms. We formulate that in Theorem 2.2.1 below.
See Sect. 3 for definitions
and results about Jacobi modular forms which we need.

Let
$$
\phi_{0,t}(\tau,z)=\sum_{k,l\in \bz} f(k,l)q^kr^l\in
J_{0,t}^{nh}
\tag{2.2.1}
$$
be a nearly holomorphic Jacobi form of weight $0$
and index $t\in \bn$ (i.e. $k$ might be negative for the Fourier expansion),
where
$q=\exp{(2\pi i \tau)}$, $\text{Im}\,\tau>0$, $r=\exp{(2\pi i z)}$,
$z\in \bc$.
It is automorphic with respect to the Jacobi group
$H(\bz) \rtimes SL_2(\bz)$ where
$H(\bz)$ is the integral Heisenberg group which is the central extension
$$
0 \to \bz \to H(\bz)\to \bz\times \bz\to 0.
$$
Nearly holomorphic means that the form $\phi_{0,t}$ is holomorphic
except a possible pole of a finite order at infinity $q=0$.
The Fourier coefficient $f(k,l)$ of $\phi_{0,t}$ depends only on
the norm  $4tk-l^2$ of $(k,l)$ and $l$ mod $2t$.
Moreover, $f(k,l)=f(k,-l)$.
From the definition of nearly holomorphic forms, it follows
that the norm $4tk-l^2$ of indices  of non-zero Fourier coefficients
$f(k,l)$ are bounded from bellow.

Let
$$
\phi^{(0)}_{0,t}(z)=\sum_{l\in \bz} f(0,l)\,r^l
\tag{2.2.2}
$$
be the $q^0$-part of $\phi_{0,t}(\tau,\,z)$.
Its Fourier coefficients are especially important for the Theorem
below.

\proclaim{Theorem 2.2.1 \cite{41, Theorem 2.1}}
Assume that $t\in \bn$ and the Fourier coefficients $f(k,l)$ of a
Jacobi form $\phi_{0,t}$
from \thetag{2.2.1} are integral.
Then the infinite product
$$
B_\phi(z)=
q^{A}r^Bs^{C}
\prod\Sb n,l,m\in \bz\\
\vspace{0.5\jot}
(n,l,m)>0\endSb
 (1-q^nr^ls^m)^{f(nm,l)},
\tag{2.2.3}
$$
where
$$
A=\frac{1}{24}\sum_{l}f(0,l),\quad
B=\frac{1}{2}\sum_{l>0}lf(0,l),\quad
C=\frac{1}{4t}\sum_{l}l^2f(0,l),
\tag{2.2.4}
$$
and $(n,l,m)>0$ means that either $m>0$ or $m=0$ and $n>0$ or
$m=n=0$ and $l<0$, defines  a meromorphic automorphic form of weight
$\frac{f(0,0)}2$ with respect to $\widehat{O}^+(L_t)$ with a character
(or a multiplier system if the weight is half-integral).
All components of the divisor of $B_\phi(z)$ á
are the rational quadratic divisors orthogonal to
$\alpha=(a,b,1)$ of the discriminant $D=-4ta+b^2>0$
(up to the action of $\widehat{O}^+(L_t)$) with multiplicities
$$
m_{D,b}=\sum_{n>0}f(n^2a,nb).
\tag{2.2.5}
$$
See other details in \cite{41, Theorem 2.1}.
\endproclaim

All 29 automorphic forms of Theorem 2.1.1 are given by
some {\it automorphic products} of Theorem 2.2.1.
It is how we call all automorphic forms of Theorem 2.2.1.
Thus, to give all these 29 automorphic forms, we should give
the corresponding 29 Jacobi forms. We give all these forms in Sect. 2.6,
Table 2.

More generally, we describe all reflective meromorphic automorphic forms
which are given by the Theorem 2.2.1. Here a meromorphic automorphic
form on the domain $\Omega(L_t)$ is called {\it reflective} if its
divisor is a sum with some multiplicities of rational quadratic
divisors orthogonal to roots of $L_t$. We remind that an element
$\alpha \in L_t$ is called {\it root} if $\alpha^2>0$ and
$\alpha^2\,|\,2(\alpha,\,L_t)$.
Using description of the divisor of $B_\phi$ in Theorem 2.2.1,
it is easy to prove

\proclaim{Lemma 2.2.2}  An infinite product $B_\phi$ given by
a Jacobi form $\phi=\phi_{0,t}$ with integral Fourier coefficients
(from Theorem 2.2.1) defines a reflective automorphic form
if and only if each non-zero Fourier coefficient $f(k,l)$
of $\phi_{0,t}$ with negative norm $4tk-l^2<0$ satisfies
$$
4tk-l^2\,|\,(4t,\,2l).
\tag{2.2.6}
$$
\endproclaim

Let us denote by $RJ_t$ the space of all Jacobi
forms $\phi=\phi_{0,t}$ of the index $t$ of Theorem 2.2.1
which give reflective automorphic products $B_\phi$
(equivalently, $\phi=\phi_{0,t}$ satisfies Lemma 2.2.2).
It is natural to call these Jacobi forms {\it reflective} either.
The space $RJ_t$ of all reflective Jacobi forms
is a free $\bz$-module with respect to addition. We have

\proclaim{Main Theorem 2.2.3} For $t\in \bn$ the space $RJ_t$ of
reflective Jacobi forms of the index $t$ is not trivial
(i.e. it is not equal to zero) if and only if $t$ is equal to
$$
\split
&1(2),\,2(3),\,3(3),\,4(3),\,5(3),\,6(4),\,7(2),\,8(3),\,9(3),\,
10(3),\,11(1),\,12(4),\,13(2),\\
&14(3),\,15(2),\,16(2),\,17(1),\,18(3),\,20(3),\,21(3),
\,22(1),\,24(2),\,25(1),\,26(1),\\
&28(1),\,30(3),\,33(1),\,34(2),\,36(3),\,39(2),\,42(1),\,45(1),\,
48(1),\,63(1),\,66(1).
\endsplit
$$
where in brackets we show the rank of the corresponding
$\bz$-module $RJ_t$ of reflective Jacobi forms.

In Table 1 of Sect. 2.5 below we give a basis of the $\bz$-module
$RJ_t$ for $t=1$, $2$, $3$, $4$, $8$, $9$, $12$, $16$, $36$
when the subspace $RJ_t$ also contains a Jacobi form which gives
the denominator identity for a Lorentzian Kac--Moody algebra
(i.e. it gives some forms of Theorem 2.1.1).
\endproclaim

All automorphic forms of Theorem 2.1.1 have the
property that they have multiplicity 1 for
components of their divisors. This follows from properties of
real roots of $\geg$ and the infinite product decomposition
\thetag{2.1.2}.
If a simple real root $\alpha \in P(\M)$ of $\geg$ is even
(i. e. $\alpha \in P(\M)_\0o$), then $m\alpha$,
$m\in \bn$, is a root of $\geg$ if and only if $m=1$,
and multiplicity of the root $\alpha$ is one.
If a simple real root $\alpha \in P(\M)$
is odd, then $m\alpha$, $m\in \bn$, is a root if and only if
$m=1$ or $m=2$.
The root $\alpha$ has multiplicity $(-1)$, and the root $2\alpha$ has
multiplicity $1$. From the full list of reflective
Jacobi forms of Main Theorem 2.2.3, it is not hard to
find all Jacobi forms $\phi$ with the divisor of multiplicity one
for $B_\phi$, since Theorem 2.2.1 gives multiplicities of
divisors of its automorphic products.
See Table 1 of Sect. 2.5 for
$t=1$, $2$, $3$, $4$, $8$, $9$, $12$, $16$, $36$.
For all other $t$ the list of all reflective Jacobi forms
of Theorem 2.2.3 will be given in a forthcoming publication.
\cite{44}. It is too long to be presented here.

Potentially, Main Theorem 2.2.3 contains information about all
reflective automorphic forms with infinite product expansion
of the type of Theorem 2.2.1 for all equivariant sublattices
$L\subset L_t$ of finite index.
Here {\it equivariant} means
that $O(L)\subset O(L_t)$; in particular, every root of $L$ is
multiple to a root of $L_t$. If $L$ has a reflective automorphic
form $\Phi$ with respect to $O(L)$ with an infinite product
expansion, then its symmetrization
$$
\prod_{g\in O(L)\backslash O(L_t)}{g^\ast\Phi}
\tag{2.2.7}
$$
is a reflective automorphic form with an infinite product expansion
for the lattice $L_t$.

Thus, potentially, Main Theorem 2.2.3 contains
important information about reflective automorphic forms with
infinite products and about automorphic forms of denominator identities
of Lorentzian Kac--Moody algebras with the symmetry lattices $L^\ast$
instead of $L_t^\ast$ and the corresponding  hyperbolic root lattices
$S^\ast$ where $S=S_t\cap L$ instead of $S_t$.
Moreover, one can possibly consider more general class of
Lie algebras for which reflective forms of Main Theorem 2.2.3 may
give kind of denominator identities. Results from \cite{26} and
\cite{50} give such a hope.

\subhead
2.3. The proof of Main Theorem 2.2.3 and
reflective hyperbolic lattices
\endsubhead
To classify finite-dimensional semi-simple and
affine Kac--Moody algebras, one needs to
classify corresponding finite and affine root systems.
To prove Main Theorem 2.2.3, one needs description of appropriate
hyperbolic root systems. They are root subsystems of
reflective hyperbolic lattices which we have considered in
Sect. 1.5.

Let $S$ be a hyperbolic (i.e. of the signature $(m,1)$) lattice,
$W(S)$ its reflection group and $\M\subset V^+(S)/\br_{++}$ its
fundamental chamber, and $Sym(\M)$ is the symmetry group of the
fundamental chamber. Thus, we have the corresponding
semi-direct product $O^+(S)=W(S)\rtimes Sym(\M)$.
We remind (see Sect. 1.5) that $S$ is called
{\it reflective} if $\M$ has a {\it generalized
Weyl vector $\rho\in
S\otimes \bq$}. It means that $\rho\not=0$ and the orbit
$Sym(\M)(\rho)$ is finite. A reflective lattice is called
{\it elliptically reflective} if it has a generalized
Weyl vector $\rho$ with $\rho^2<0$. It is called
{\it parabolically reflective} if it is not
elliptically reflective but has
a generalized Well vector $\rho$ with $\rho^2=0$. It is called
{\it hyperbolically reflective} if it is not elliptically or
parabolically reflective, but it has a generalized
Weyl vector $\rho$ with $\rho^2>0$.

Suppose that $B_\phi(z)$ is an automorphic form of
Theorem 2.2.1 which is reflective.
The inequality $(n,l,m)>0$ of Theorem 2.2.1 is a variant of
choosing a fundamental chamber $\M$ of $W(S_t)$.
The vector $\rho=(A,B,C)$ is invariant with
respect to the group $\widehat{Sym}(\M)=
Sym(\M)\cap \widehat{O}^+(L_t)$ which
has finite index in $Sym(\M)$. If $\rho=(A,B,C)$ is not zero,
it then defines a generalized Weyl vector for
$Sym(\M)$. If the form $B_\phi(z)$ has a zero vector
$\rho=(A,B,C)$, one
can change $B_\phi(z)$ by other form $w^\ast (B_\phi(z))$
where $w\in W$ is an appropriate reflection, in such a way that
$w^\ast (B_\phi(z))$ will have a non-zero generalized Weyl
vector $\rho$. Thus, we get

\proclaim{Lemma 2.3.1} If the space $RJ_t$ of reflective Jacobi forms
is not zero, then the lattice $S_t$ is reflective.
\endproclaim

It is interesting that the space $RJ_t$ may really have a
Jacobi form with zero vector $\rho=(A,B,C)$. It happens for
$t=6$ and $t=12$ when $\rk RJ_t=4$.

The monograph \cite{68} was devoted to classification of
reflective hyperbolic lattices and appropriate
(for Lorentzian Kac--Moody algebras) hyperbolic root systems of
the rank three. This classification contains 122 main elliptic
types and 66 main hyperbolic types (there are no main parabolic types).
In particular, all reflective hyperbolic lattices of the rank three
with square-free determinant were classified (e. g. it gives
classification of all reflective lattices $S_t$ when $t$
is square-free, see Theorem 2.3.2 below). It follows
classification of all maximal reflective hyperbolic
lattices of rank three. For arbitrary reflective
hyperbolic lattices of rank three, effective methods of their
enumeration and effective estimates of their
invariants were obtained.

Using these results and direct calculations, we get

\proclaim{Theorem 2.3.2} The lattice $S_t=H\oplus \langle 2t \rangle$
is reflective for the
following and the only following $t\in \bn$,
where in brackets we put the
type of reflectivity of the lattice, (e) for elliptic, (p) for
parabolic and (h) for the hyperbolic type:
$$
\split
t\,=\,&1\ -\ 22\,(e),\ 23\,(h),\ 24\ -\ 26\,(e),\ 28\,(e),\
29\,(h),\ 30\,(e),\ 31\,(h),\\
&33\,(e),\ 34\,(e),\ 35\,(h),\ 36\,(e),\ 37\,(h),\ 38\,(h),\ 39\,(e),\
40\,(h),\ 42\,(e),\\
&44\,(h),\ 45\,(e),\ 46\,(h),\ 48\,(h),\ 49\,(e),\ 50\,(e),\ 52\,(e),\
55\,(e),\ 56\,(h),\\
&57\,(h),\ 60\,(h),\ 63\,(h),\ 66\,(e),\ 70\,(h),\ 72\,(h),\
78\,(h),\ 84(h),\ 90\,(h),\\
&100\,(h),\ 105\,(h).
\endsplit
$$
In particular, there are no reflective lattices $S_t$ of parabolic type. 
\endproclaim

To prove Main Theorem 2.2.3 and to find a basis of $RJ_t$,
one needs to analyze only the $t$ of Theorem 2.3.2. In particular,
$t\le 105$. To find the rank of $RJ_t$ and to find  a basis of
$RJ_t$, one can use known generators of the graded ring of
weak Jacobi forms. Over $\bq$ (with rational Fourier coefficients)
these generators were found in \cite{24} (see Sect. 3), and this is
sufficient for calculation of $\rk RJ_t$ (e. g. one can use computer).
Generators of the graded ring of weak Jacobi forms with integral
Fourier coefficients were found in \cite{33} and \cite{34}
(see Sect. 3). There results permit to find basis of $RJ_t$ for all
$t$ of Theorem 2.2.3. We give it for
$t=1$, $2$, $3$, $4$, $8$, $9$, $12$, $16$, $36$ in
Table 1 of Sect. 2.5 below.
To simplify these calculations, one can also use arguments which
we give below for the proof of Theorem 2.1.1.

\smallpagebreak

Classification of all reflective hyperbolic lattices of
rank three in \cite{68} gives a hope that all Lorentzian
Kac--Moody algebras of the rank three (with arbitrary hyperbolic
root lattice $S$ of the rank three) can be classified.

Finiteness results for all reflective hyperbolic lattices
(Theorems 1.5.3 and 1.5.5) and the reasonable number of
main types for reflective hyperbolic lattices of
the rank three ($122+66$) give a hope that all
Lorentzian Kac--Moody algebras of rank $\ge 3$ can be
classified in a future. We expect that number of cases drops
when the rank is increasing. We think that the rank three case
is the most reach and complicated.

\subhead
2.4. Proof of Theorem 2.1.1 and reflective hyperbolic lattices
with a Weyl vector
\endsubhead
Here we sketch the proof of Theorem 2.1.1 to emphasize importance of
reflective hyperbolic lattices and hyperbolic root systems with
a Weyl vector.

Using Theorem 2.2.3 and considerations above,
we have proved that Theorem 2.1.1 gives all automorphic forms
(satisfying its conditions) which
can be obtained by the general construction of Theorem 2.2.1.
Here we want to show that there are no other automorphic forms
satisfying conditions of Theorem 2.1.1. Thus, any of forms of
Theorem 2.1.1 can be obtained by the construction of
Theorem 2.2.1.

For Theorem 2.1.1, the fundamental polyhedron $\M$ and the set
$P(\M)$ of orthogonal roots to $\M$ have a Weyl
vector $\rho$ (satisfying \thetag{1.4.2}). The $P(\M)$ and $\rho$
are invariant with respect to the group $\widehat{Sym}(P(\M))$
which has finite index in $Sym(\M)$ (we use notation
$\widehat{G}=G\cap \widehat{O}^+(L_t)$). In
particular, the corresponding lattices $S_t$ are reflective.
For any reflective lattice $S_t$ (it belongs to the list of
Theorem 2.3.2) the fundamental chamber $\M_0$ for
the full reflection group $W(S_t)$ can be calculated and is known (for
$t=1$, $2$, $3$, $4$, $8$, $9$, $12$, $16$, $36$
these calculations are presented in Table 1 below).
Thus, the fundamental chamber $\M$ is composed from the known
polygon $\M_0$ by some reflections, and it has a Weyl
vector $\rho$. Using this information, we can find all possible
$\M$, $P(\M)$, $\rho$ and predict the divisor
of the reflective automorphic form $\Phi(z)$. Looking at the list
of 29 forms of Theorem 2.1.1, one can see that one of them 
(for the corresponding $t$) has the same divisor. By
Koecher principle (we discussed this principle in Sect. 1.6), 
the form $\Phi(z)$ is equal to that form.

Similar arguments can be used to classify all {\it reflective
meromorphic automorphic forms with infinite product}
of the type \thetag{2.1.2} and with a
generalized Weyl vector $\rho\not=0$. But we
don't require that multiplicities of the infinite product are
related with Fourier coefficients of any modular form.
Like the product \thetag{2.1.2},
this product should be related with a reflection subgroup
$W\subset W(S_t)$, its fundamental chamber $\M$, the set
$P(\M)\subset S_t^\ast$ of orthogonal roots to $\M$
(they define $\Delta_+\subset S_t^\ast$),
and a generalized Weyl vector
$\rho \in S_t\otimes \bq$ (i. e.  $\rho\not=0$, the orbit
$Sym(P(\M))(\rho)$ is finite and
$W\rtimes Sym(P(\M))$ has finite index in $O(S_t)$). The function
$\mult(\alpha)$, $\alpha \in \Delta_+$, should be integral and
invariant with respect to $\widehat{Sym}(P(\M))$.
The product should converge in a neighborhood
$im (z)^2<<0$ of the cusp at infinity.
All the definitions are similar. We have

\proclaim{Theorem 2.4.1} Any reflective meromorphic automorphic form
with the root lattice $S_t^\ast$, the symmetry lattice $L_t^\ast$ and
the symmetry group $\widehat{O}^+(L_t)$ having infinite product
of the type \thetag{2.1.2} with a non-zero generalized
Weyl vector $\rho$ is $B_\phi$ where $\phi$ belongs to the list of
Main Theorem 2.2.3.
\endproclaim

Applying to the forms of Theorem 2.4.1 reflections with respect
to roots in $S_t$, one can get some automorphic forms with zero
generalized Weyl vector and with infinite product.
They appear only for
$t=6$ and $t=12$.

\subhead
2.5. Reflective Jacobi forms from $RJ_t$ for
$t=1$, $2$, $3$, $4$, $8$, $9$, $12$, $16$, $36$
\endsubhead
For these $t$ we give the basis
$\xi_{0,t}^{(1)},\,\dots\,\xi_{0,t}^{(rk)}$ of the
$\bz$-module $RJ_t$. For each Jacobi form of the basis we show
the leading part of its Fourier expansion which defines the
Jacobi form uniquely. We give all its Fourier coefficients
with negative norm (up to equivalence); the corresponding
negative norm is shown in brackets $[\cdot]$. We also give
a formula for the form which uses basic Jacobi forms.
In these formulae $E_4=E_4(\tau)$ and $\Delta_{12}=\Delta(\tau)$
are the Eisenstein series of weight $4$ and the Ramanujan function of
weight $12$ for $SL_2(\bz )$ respectively, $E_{4,m}$ ($m=1,\,2,\,3$)
are  Eisenstein-Jacobi series of weight 4 and index $m$
(see \cite{24}), and
$\phi_{0,1}$, $\phi_{0,2}$, $\phi_{0,3}$, $\phi_{0,4}$
are the four generators from \cite{33} and \cite{41}
of the graded ring of the weak Jacobi forms
of weight zero with integral Fourier coefficients.
See Sect. 3 about these Jacobi forms.

We give the set $\overline{R}$ of primitive roots in
$S_t^\ast$ up to equivalence (up to the action of the group
$\pm \widehat{O}(L_t)$). Up to this equivalence, a root
$\alpha=(n,\,l,\,m)$ is
defined by its norm $-2t\alpha^2=-4nm+l^2$ and $\pm l\mod 2t$.
We also give the matrix
$$
\text{Mul}(\overline{R},\,\xi)=\text{mul}(\gamma_i,\,\xi_{0,t}^{(j)}),
$$
where $\text{mul}(\gamma_i,\,\xi_{0,t}^{(j)})$ is the multiplicity of
the form $\Phi_{\xi_{0,t}^{(j)}}$ in
the rational quadratic divisor which is orthogonal to the root from
the equivalence class $\gamma_i\in \overline{R}$.

We give the set $P(\M_0)$ of primitive roots in $S_t^\ast$ which is
orthogonal to the fundamental chamber $\M_0$ of the
reflection group $W(S_t)$ (this is equivalent to
the ordering $(n,l,m)>0$ used in Theorem 2.2.1), and its Gram matrix
$$
G(P(\M_0))=2t((\alpha,\,\beta)),\ \alpha,\,\beta \in P(\M_0).
$$
Thus, we identify the dual lattice $S_t^\ast$ with the integral lattice
$S_t^\ast(2t)=H(2t)\oplus \langle 1 \rangle$ to make it integral.

All these data are given in Table 1 below.

\vskip40pt

\centerline{\bf Table 1. Basis of the space $RJ_t$ of
reflective}

\centerline{\bf Jacobi forms for
$t=1$, $2$, $3$, $4$, $8$, $9$, $12$, $16$, $36$.}

\vskip20pt

{\bf Case $t=1$.} The space $RJ_1$ has the basis
$$
\split
\xi_{0,1}^{(1)}=&\phi_{0,1}=
(r[-1] + 10 + r^{-1}[-1]) + O(q);\\
\xi_{0,1}^{(2)}=&{E_4}^2{E_{4,1}}/{\Delta_{12}}-57\phi_{0,1}\\
=&q^{-1}[-4] +(r^2[-4]-r[-1]+60-r^{-1}[-1]+r^{-2}[-4]) + O(q).
\endsplit
$$
We have $\overline{R}=\overline{P(\M_0)}$ and
$$
P(\M_0)=
\pmatrix
{1}&{2}&{0}\cr
{0}&{-1}&{0}\cr
{-1}&{0}&{1}\cr
\endpmatrix \equiv
\left[\matrix
4,\,\overline{0}\cr
1,\,\overline{1}\cr
4,\,\overline{0}\cr
\endmatrix\right];
\hskip10pt
\text{Mul}(P(\M_0),\,\xi)=
\pmatrix
0&1\cr
1&0\cr
0&1\cr
\endpmatrix;
$$
$$
G(P(\M_0))=
\pmatrix
{4}&{-2}&{-2}\cr
{-2}&{1}&{0}\cr
{-2}&{0}&{4}\cr
\endpmatrix.
$$

\vskip20pt

{\bf Case $t=2$.} The space $RJ_2$ has the basis
$$
\split
\xi_{0,2}^{(1)}=&{\phi_{0,2}}=
(r[-1] + 4 + r^{-1}[-1]) + O(q);\\
\xi_{0,2}^{(2)}=&({\phi_{0,1}})^2-21{\phi_{0,2}}=
(r^2[-4]- r[-1]+18-r^{-1}[-1]+r^{-2}[-4]) + O(q);\\
\xi_{0,2}^{(3)}=&{E_4}^2{E_{4,2}}/{\Delta_{12}}-
14({\phi_{0,1}})^2+216{\phi_{0,2}}=
q^{-1}[-8] + 24 + O(q).
\endsplit
$$
We have $\overline{R}=\overline{P(\M_0)}$ where
$$
P(\M_0)=
\pmatrix
{1}&{2}&{0}\cr
{0}&{-1}&{0}\cr
{-1}&{0}&{1}\cr
\endpmatrix
\equiv
\left[\matrix
4,\,\overline{2}\\
1,\,\overline{1}\\
8,\,\overline{0}\\
\endmatrix\right];
\hskip20pt
\text{Mul}(P(\M_0),\,\xi)=
\pmatrix
0&1&0\cr
1&0&0\cr
0&0&1\cr
\endpmatrix;
$$
$$
G(P(\M_0))=
\pmatrix
{4}&{-2}&{-4}\cr
{-2}&{1}&{0}\cr
{-4}&{0}&{8}\cr
\endpmatrix.
$$

\vskip20pt

{\bf Case $t=3$.} The space $RJ_3$ has the basis
$$
\split
\xi_{0,3}^{(1)}=&{\phi_{0,3}}=
(r[-1] + 2 + r^{-1}[-1]) + O(q);\\
\xi_{0,3}^{(2)}=&{\phi_{0,1}}{\phi_{0,2}}-15{\phi_{0,3}}=
(r^2[-4] -r[-1] + 12 -r^{-1}[-1]+ r^{-2}[-4]) + O(q);\\
\xi_{0,3}^{(3)}=&{E_4}^2{E_{4,3}}/{\Delta_{12}}-2({\phi_{0,1}})^3+
33{\phi_{0,1}}{\phi_{0,2}}+90{\phi_{0,3}}=
q^{-1}[-12] + 24 + O(q).
\endsplit
$$
We have $\overline{R}=\overline{P(\M_0)}$ where
$$
P(\M_0)=
\pmatrix
{1}&{2}&{0}\cr
{0}&{-1}&{0}\cr
{-1}&{0}&{1}\cr
\endpmatrix
\equiv
\left[\matrix
4,\,\overline{2}\\
1,\,\overline{1}\\
12,\,\overline{0}\\
\endmatrix\right];
\hskip20pt
\text{Mul}(P(\M_0),\,\xi)=
\pmatrix
0&1&0\cr
1&0&0\cr
0&0&1\cr
\endpmatrix;
$$
$$
G(P(\M_0))=
\pmatrix
{4}&{-2}&{-6}\cr
{-2}&{1}&{0}\cr
{-6}&{0}&{12}\cr
\endpmatrix.
$$

\vskip20pt

{\bf Case $t=4$.} The space $RJ_4$ has the basis
$$
\split
\xi_{0,4}^{(1)}=&{\phi_{0,4}}=(r[-1] + 1 + r^{-1}[-1]) + O(q);\\
\xi_{0,4}^{(2)}=&({\phi_{0,2}})^2-9{\phi_{0,4}}=
(r^2[-4] - r[-1] + 9 -r^{-1}[-1]+ r^{-2}[-4]) + O(q);\\
\xi_{0,4}^{(3)}=&{E_4}{E_{4,1}}{E_{4,3}}/{\Delta_{12}}-
2({\phi_{0,1}})^2{\phi_{0,2}}+20{\phi_{0,1}}{\phi_{0,3}}+16{\phi_{0,4}}\\
=&q^{-1}[-16] + 24 + O(q).
\endsplit
$$
We have $\overline{R}=\overline{P(\M_0)}$ where
$$P(\M_0)=
\pmatrix
{1}&{2}&{0}\cr
{0}&{-1}&{0}\cr
{-1}&{0}&{1}\cr
\endpmatrix
\equiv
\left[\matrix
4,\,\overline{2}\\
1,\,\overline{1}\\
16,\,\overline{0}\\
\endmatrix\right];
\hskip20pt
\text{Mul}(P(\M_0),\,\xi)=
\pmatrix
0&1&0\cr
1&0&0\cr
0&0&1\cr
\endpmatrix;
$$
$$
G(P(\M_0))=
\pmatrix
{4}&{-2}&{-8}\cr
{-2}&{1}&{0}\cr
{-8}&{0}&{16}\cr
\endpmatrix.
$$

\vskip20pt

{\bf Case $t=8$.} The space $RJ_8$ has the basis
$$
\split
\xi_{0,8}^{(1)}=&({\phi_{0,2}})^2{\phi_{0,4}}-
{\phi_{0,2}}({\phi_{0,3}})^2-({\phi_{0,4}})^2=
(2r[-1] - 1 + 2r^{-1}[-1]) \\
+&(-r^6[-4] - 2r^5 + 4r^4 - 4r^3 + r^2 + 6r - 8+\cdots )q+O(q^2);\\
\xi_{0,8}^{(2)}=&\phi_{0,2}(\tau,2z)=
{\phi_{0,1}}{\phi_{0,3}}{\phi_{0,4}}+{\phi_{0,2}}({\phi_{0,3}})^2-
2({\phi_{0,2}})^2{\phi_{0,4}}-2({\phi_{0,4}})^2\\
=&(r^2[-4]+ 4 + r^{-2}[-4])+(r^6[-4]-8r^4-r^2+16-\cdots)q+O(q^2);\\
\xi_{0,8}^{(3)}=&{E_4}{E_{4,3}}\Big({E_{4,2}}{\phi_{0,3}}-
{E_{4,1}}{\phi_{0,4}}\Big)/{\Delta_{12}}-
3({\phi_{0,1}})^2({\phi_{0,3}})^2+
2({\phi_{0,1}})^2{\phi_{0,2}}{\phi_{0,4}}\\
+&{\phi_{0,1}}{\phi_{0,3}}{\phi_{0,4}}-16({\phi_{0,4}})^2=
q^{-1}[-32] + 24 + (8r^6[-4] + 256r^5 + 2268r^4\\
+&9472r^3+23608r^2 + 39424r + 46812 + \cdots)q+O(q^2).
\endsplit
$$
We have $\overline{R}=
\overline{P(\M_0)}$ where
$$
P(\M_0)=\pmatrix
{1}&{2}&{0}\cr
{0}&{-1}&{0}\cr
{-1}&{0}&{1}\cr
{1}&{6}&{1}\cr
\endpmatrix
\equiv
\left[\matrix
4,\,\overline{2}\\
1,\,\overline{1}\\
32,\,\overline{0}\\
4,\,\overline{6}
\endmatrix\right];
\hskip20pt
\text{Mul}(P(\M_0),\,\xi)=
\pmatrix
0  &1  &0\cr
2  &1  &0\cr
0  &0  &1\cr
-1 &1  &8
\endpmatrix;
$$
$$
G(P(\M_0))=
\pmatrix
{4}&{-2}&{-16}&{-4}\cr
{-2}&{1}&{0}&{-6}\cr
{-16}&{0}&{32}&{0}\cr
{-4}&{-6}&{0}&{4}\cr
\endpmatrix.
$$

\vskip20pt

{\bf Case $t=9$.} The space $RJ_9$ has the basis
$$
\split
\xi_{0,9}^{(1)}=&-{\phi_{0,1}}({\phi_{0,4}})^2+
6{\phi_{0,2}}{\phi_{0,3}}{\phi_{0,4}}-5({\phi_{0,3}})^3=
(3r[-1] - 2 + 3r^{-1}[-1]) \\
+&(-4r^6 + 6r^5 - 12r^4 + 22r^3 - 30r^2 + 36r -36 +\cdots )q+
(-r^9[-9] - 6r^8 \\
+&15r^7 - 36r^6 + 72r^5 - 120r^4 + 171r^3 - 216r^2 +
255r -268 +\cdots)q^2 + O(q^3);\\
\xi_{0,9}^{(2)}=&{\phi_{0,1}}({\phi_{0,4}})^2-
5{\phi_{0,2}}{\phi_{0,3}}{\phi_{0,4}}+4({\phi_{0,3}})^3=
(r^2[-4]-r[-1]+4- \cdots )\\
+&(3r^6-8r^5+9r^4-24r^3+31r^2-32r+42-\cdots)q+
(r^9[-9]+7r^8-15r^7\\
+&33r^6-80r^5+110r^4-177r^3+219r^2-241r+286-
\cdots)q^2+O(q^3);\\
\xi_{0,9}^{(3)}=&{E_{4,2}}{E_{4,3}}\Big({E_{4,1}}{\phi_{0,3}}-
{E_4}{\phi_{0,4}}\Big)/{\Delta_{12}}-
3{\phi_{0,1}}{\phi_{0,2}}({\phi_{0,3}})^2+
2({\phi_{0,1}})^2{\phi_{0,3}}{\phi_{0,4}}\\
-&30{\phi_{0,1}}({\phi_{0,4}})^2+
27{\phi_{0,2}}{\phi_{0,3}}{\phi_{0,4}}+9({\phi_{0,3}})^3=
q^{-1}[-36] + 24 \\
+&(33r^6 + 486r^5 +3159r^4 + 10758r^3 + 24057r^2 + 37908r +
44082 +\cdots)q\\
+&(2r^9[-9] + 243r^8 + 5346r^7 +
44055r^6 + 204120r^5 + 642978r^4 +
1483416r^3 \\
+&2632905r^2 + 3679020r + 4109590 +\cdots)q^2+O(q^3).
\endsplit
$$
We have
$\overline{R}=\overline{P(\M_0)}$ where
$$
P(\M_0)=\pmatrix
{1}&{2}&{0}\cr
{0}&{-1}&{0}\cr
{-1}&{0}&{1}\cr
{2}&{9}&{1}\cr
\endpmatrix\
\equiv
\left[\matrix
4 & \overline{2}\\
1&   \overline{1}\\
36& \overline{0}\\
9&\overline{9}
\endmatrix\right];
\hskip10pt
\text{Mul}(P(\M_0),\,\xi)=
\pmatrix
0&1&0\\
3&0&0\\
0&0&1\\
-1&1&3
\endpmatrix;
$$
$$
G(P(\M_0))=
\pmatrix
{4}&{-2}&{-18}&{0}\cr
{-2}&{1}&{0}&{-9}\cr
{-18}&{0}&{36}&{-18}\cr
{0}&{-9}&{-18}&{9}\cr
\endpmatrix.
$$

\vskip20pt

{\bf Case $t=12$.} The space $RJ_{12}$ has the basis
$$
\split
\xi_{0,12}^{(1)}=&\bigl(\vartheta(\tau,z)/\eta(\tau)\bigr)^{12}=
(r^8[-16] - 8r^7[-1] + 24r^6 - 24r^5 - 36r^4 + 120r^3 \\
-&88r^2 - 88r + 198 - \cdots)q + (-4r^{10}[-4] + 24r^9 - 32r^8 - 104r^7 +
396r^6 \\
-&352r^5 - 512r^4 + 1440r^3 - 904r^2 - 1008r + 2112 - \cdots)q^2+O(q^3);
\endsplit
$$
$$
\split
\xi_{0,12}^{(2)}=&3{\phi_{0,2}}({\phi_{0,3}})^2{\phi_{0,4}}-
({\phi_{0,2}})^2({\phi_{0,4}})^2-2({\phi_{0,3}})^4-({\phi_{0,4}})^3\\
=&(r[-1] - 1 + r^{-1}[-1]) +
(-r^7[-1] + r^6 - r^5 + r^4 - r^2 +2r - 2 + \cdots)q\\
+&(-r^{10}[-4] + r^8 - 2r^7 + 3r^6 - 3r^5 +2r^4 -
2r^2 +5r -6 +  \cdots)q^2 + O(q^3);
\endsplit
$$
$$
\split
\xi_{0,12}^{(3)}=&2({\phi_{0,2}})^2({\phi_{0,4}})^2-
5{\phi_{0,2}}({\phi_{0,3}})^2{\phi_{0,4}}+3({\phi_{0,3}})^4+
({\phi_{0,4}})^3=
(r^2[-4] - r[-1] \\
+&3 - \cdots)+  (r^7[-1] - 3r^6 + r^5 - 3r^4 +3r^3 -2r +6 -\cdots)q\\
+&(2r^{10}[-4] - 3r^8 + 2r^7 - 9r^6 + 3r^5 - 6r^4 + 7r^2 - 5r +
18 -\cdots)q^2 + O(q^3);\\
\xi_{0,12}^{(4)}=&
{E_{4,3}}\Big({E_{4,1}}{E_{4,2}}({\phi_{0,3}})^2-
2{E_4}{E_{4,2}}{\phi_{0,3}}{\phi_{0,4}}+
{E_4}{E_{4,1}}({\phi_{0,4}})^2\Big)/{\Delta_{12}}\\
-&2({\phi_{0,1}})^2{\phi_{0,2}}({\phi_{0,4}})^2+
+5({\phi_{0,1}})^2({\phi_{0,3}})^2{\phi_{0,4}}-
3{\phi_{0,1}}{\phi_{0,2}}({\phi_{0,3}})^3\\
-&36{\phi_{0,1}}{\phi_{0,3}}({\phi_{0,4}})^2+
24{\phi_{0,2}}({\phi_{0,3}})^2{\phi_{0,4}}+
9({\phi_{0,3}})^4+16({\phi_{0,4}})^3\\
=&q^{-1}[-48] + 24 + (24r^7[-1] + 264r^6 +
1608r^5 + 5610r^4 +13464r^3 \\
+& 24312r^2 + 34056r  + 38208 +\cdots)q
+(12r^{10}[-4] + 440r^9 + 5544r^8 \\
+& 34104r^7 + 135388r^6 + 395808r^5
+902352r^4 + 1667360r^3 + 2550552r^2 \\
+&3276240r + 3558160 + \cdots)q^2+O(q^3)
\endsplit.
$$
We have
$$
P(\M_0)=\pmatrix
{1}&{2}&{0}\cr
{0}&{-1}&{0}\cr
{-1}&{0}&{1}\cr
{1}&{8}&{1}\cr
\endpmatrix
\equiv
\left[\matrix
4&\overline{2}\\
1&\overline{1}\\
48&\overline{0}\\
16&\overline{8}
\endmatrix\right];\
G(P(\M_0))=\pmatrix
{4}&{-2}&{-24}&{-8}\cr
{-2}&{1}&{0}&{-8}\cr
{-24}&{0}&{48}&{0}\cr
{-8}&{-8}&{0}&{16}\cr
\endpmatrix.
$$
$$
\overline{R}=
\left[\matrix
4,\overline{2}\\
1,\overline{1}\\
48,\overline{0}\\
16,\overline{8}\\
4,\overline{10}\\
1,\overline{7}\\
\endmatrix\right];
\ \ \
\text{Mul}(\overline{R},\,\xi)=
\pmatrix
0  & 0 & 1 & 0\\
0  & 1 & 0 & 0\\
0  & 0 & 0 & 1\\
1  & 0 & 0 & 0\\
-4 &-1 & 2 &12\\
-12&-2 & 3 &36
\endpmatrix;
$$

\vskip20pt

{\bf Case $t=16$.} The space $RJ_{16}$ has the basis
$$
\split
\xi_{0,16}^{(1)}=&\phi_{0,4}(\tau,\,2z)=
(r^2[-4] +1 + r^{-2}[-4]) +
(-r^8 - r^6 + r^2 + 2 + \cdots)q^2\\
+&(r^{14}[-4] -2r^{10} - 4r^8 - 4r^6 + 5r^2 + 8 +\cdots)q^3+O(q^4);
\endsplit
$$
$$
\split
\xi_{0,16}^{(2)}=&
{E_{4,3}}\Big({E_4}{E_{4,1}}({\phi_{0,3}})^4-
({E_4})^2({\phi_{0,3}})^3{\phi_{0,4}}-
2{E_{4,1}}{E_{4,2}}({\phi_{0,3}})^2{\phi_{0,4}}\\
+&{E_4}{E_{4,2}}{\phi_{0,3}}({\phi_{0,4}})^2-
{E_{4,1}}{E_{4,2}}({\phi_{0,3}})^2{\phi_{0,4}} +
2{E_4}{E_{4,2}}{\phi_{0,3}}({\phi_{0,4}})^2\\
-&{E_4}{E_{4,1}}({\phi_{0,4}})^3\Big)/{\Delta_{12}}+
2({\phi_{0,1}})^3({\phi_{0,3}})^3{\phi_{0,4}}-
3({\phi_{0,1}})^2{\phi_{0,2}}({\phi_{0,3}})^4\\
-&7({\phi_{0,1}})^2({\phi_{0,3}})^2({\phi_{0,4}})^2-
31{\phi_{0,1}}{\phi_{0,2}}({\phi_{0,3}})^3{\phi_{0,4}}+
46{\phi_{0,1}}({\phi_{0,3}})^5\\
+&72{\phi_{0,1}}{\phi_{0,3}}({\phi_{0,4}})^3+
7({\phi_{0,2}})^3({\phi_{0,3}})^2{\phi_{0,4}}-
72{\phi_{0,2}}({\phi_{0,3}})^2({\phi_{0,4}})^2\\
-&197({\phi_{0,3}})^4{\phi_{0,4}}+
2({\phi_{0,1}})^2{\phi_{0,2}}({\phi_{0,4}})^3+
21({\phi_{0,3}})^4{\phi_{0,4}}-26({\phi_{0,4}})^4\\
+&2{\phi_{0,1}}({\phi_{0,2}})^2{\phi_{0,3}}({\phi_{0,4}})^2-
({\phi_{0,2}})^2({\phi_{0,3}})^4-4({\phi_{0,2}})^2({\phi_{0,4}})^3-
2({\phi_{0,2}})^4({\phi_{0,4}})^2\\
=&q^{-1}[-64]+(8r[-1] +14 + 8r^{-1}[-1]) +
(21r^8 + 200r^7 + 1036r^6 \\
+&3360r^5 +8100r^4 + 15240r^3 + 23604r^2 + 30352r + 33058 + \cdots)q\\
+&(56r^{11} + 1008r^{10} + 7336r^9 + 32932r^8 + 108800r^7 +
283504r^6 +610344r^5 \\
+&1112832r^4 + 1750728r^3 + 2401952r^2 + 2896688r^1 +
3081400 + \cdots)q^2\\
+&(4r^{14}[-4] +560r^{13} + 8092r^{12} + 58328r^{11} +
283784r^{10} +1042328r^9 \\
+&3082176r^8 + 7616904r^7 + 16136000r^6 + 29802144r^5 +
48582612r^4 \\
+&70497736r^3 +91619124r^2 + 107054192r +
112732002 + \cdots)q^3+O(q^4)
\endsplit.
$$
We have $\overline{R}=\overline{P(\M_0)}$ where
$$
P(\M_0)=
\pmatrix
{1}&{2}&{0}\cr
{0}&{-1}&{0}\cr
{-1}&{0}&{1}\cr
{5}&{32}&{3}\cr
{3}&{14}&{1}\cr
\endpmatrix
\equiv
\left[\matrix
4&\overline{2}\\
1&\overline{1}\\
64&\overline{0}\\
64&\overline{0}\\
4&\overline{14}\\
\endmatrix\right];
\hskip10pt
\text{Mul}(P(\M_0),\,\xi)=
\pmatrix
1 & 0\\
1&8\\
0&1\\
0&1\\
1&4
\endpmatrix;
$$
$$
G(P(\M_0))=
\pmatrix
{4}&{-2}&{-32}&{-32}&{-4}\cr
{-2}&{1}&{0}&{-32}&{-14}\cr
{-32}&{0}&{64}&{-64}&{-64}\cr
{-32}&{-32}&{-64}&{64}&{0}\cr
{-4}&{-14}&{-64}&{0}&{4}\cr
\endpmatrix.
$$

\vskip20pt

{\bf Case $t=36$.} The space $RJ_{36}$ has the basis
$$
\split
\xi_{0,36}^{(1)}=&3\xi_{0,36}^{(2)}- \xi_{0,9}^{(1)}(\tau,2z)\\
=&(-3r[-1] +5 - 3r^{-1}[-1]) +
(r^{12} + 3r^{11}+\cdots )q\\
+&(r^{18}[-36] - 3r^{17}[-1] +9r^{16} +\cdots)q^2+
(6r^{20} - 3r^{19}+\cdots )q^3\\
+&(4r^{24} - 15r^{22} +\cdots )q^4+
(3r^{27}[-9] - 9r^{26}+3r^{25}+\cdots)q^5\\
+&(3r^{29} + 6r^{28} +\cdots )q^6+
(3r^{32}[-16] -25r^{30} + 9r^{29}+\cdots)q^7\\
+&(-3r^{33} + 33r^{32}+\cdots )q^8+O(q^9) ;
\endsplit
$$
$$
\split
\xi_{0,36}^{(2)}=&
(\vartheta (\tau,\, 10z)\vartheta (\tau,\,z))/
(\vartheta (\tau , \, 5z)\vartheta (\tau,\,2z))\\
=&((r^2[-4] - r[-1] + 1 - r^{-1}[-1] + r^{-2}[-4]) +
(-r^{12} + r^{11} - r^{10} + \cdots )q\\
+&(-r^{17}[-1] + r^{16} - r^{15}+\cdots )q^2+
(-r^{19} + 2r^{18} - 3r^{17} +\cdots )q^3\\
+&(-r^{21}  + 2r^{20} - 4r^{19}+\cdots )q^4+
(r^{27}[-9] -r^{26} + r^{25}+\cdots )q^5\\
+&(r^{29} - 2r^{28}  + 3r^{27} + \cdots )q^6+
(r^{32}[-16] -r^{30}+ 3r^{29}+\cdots )q^7\\
+&(r^{34}[-4] - r^{33} + r^{32} - 3r^{30}+\cdots)q^8 + O(q^9) ;
\endsplit
$$
$$
\split
{}&\xi_{0,36}^{(3)}
=\bigl[(-36 \phi_{0,3}^4 + 56 \phi_{0,4}^3) \phi_{0,4}^3 \phi_{0,3}^3 E_4^2 +
(45 \phi_{0,3}^8 - 126 \phi_{0,4}^3 \phi_{0,3}^4 +
\phi_{0,4}^6) \phi_{0,4}^2 E_4 E_{4,1}\\
{}&+
(- 10 \phi_{0,3}^8 + 126 \phi_{0,4}^3 \phi_{0,3}^4
- 8 \phi_{0,4}^6) \phi_{0,4} \phi_{0,3} E_4 E_{4,2}\\ 
{}&+
(\phi_{0,3}^8 - 84 \phi_{0,4}^3 \phi_{0,3}^4
+ 28 \phi_{0,4}^6)\phi_{0,3}^2 E_{4,1} E_{4,2}\bigr] E_{4,3}/\Delta\\
{}&+
\bigr[770 \phi_{0,4}^6 \phi_{0,3}^2- 731 \phi_{0,4}^7 \phi_{0,2}
- 731 \phi_{0,4}^6 \phi_{0,2}^3 +
2924 \phi_{0,4}^5 \phi_{0,3}^2 \phi_{0,2}^2\\
{}&
 - 3655 \phi_{0,4}^4 \phi_{0,3}^4 \phi_{0,2}
- 29 \phi_{0,4}^4 \phi_{0,3}^2 \phi_{0,2}^4
+ 133 \phi_{0,4}^3 \phi_{0,3}^4 \phi_{0,2}^3 +
+ 1472 \phi_{0,4}^3 \phi_{0,3}^6\\
{}&
- 225 \phi_{0,4}^2 \phi_{0,3}^6 \phi_{0,2}^2+ 167 \phi_{0,4} \phi_{0,3}^8
\phi_{0,2} - 46 \phi_{0,3}^{10}\bigl]D_{0,6}\\
{}&
+(72 \phi_{0,3}^4 - 
112\phi_{0,4}^3) \phi_{0,4}^3 \phi_{0,3}^3\phi_{0,1}^3 
+ \bigl[-731 \phi_{0,4}^6 \phi_{0,2}^3 + 
1462 \phi_{0,4}^5 \phi_{0,3}^2\phi_{0,2}^2\\ 
{}&+
(-126 \phi_{0,4} \phi_{0,3}^8 + 1039 \phi_{0,4}^4 \phi_{0,3}^4
- 733 \phi_{0,4}^7) \phi_{0,2}\\
{}&
+ 29  \phi_{0,3}^{10} - 1615 \phi_{0,4}^3 \phi_{0,3}^6
- 714 \phi_{0,4}^6 \phi_{0,3}^2\bigr]\phi_{0,4} \phi_{0,1}^2\\
{}&
+
\bigl[12425 \phi_{0,4}^6 \phi_{0,3} \phi_{0,2}^4 
- 50600 \phi_{0,4}^5 \phi_{0,3}^3 \phi_{0,2}^3\\
{}&
+ (67608 \phi_{0,4}^4 \phi_{0,3}^5 + 
20633 \phi_{0,4}^7 \phi_{0,3})\phi_{0,2}^2
-(3 \phi_{0,3}^{11} + 37314 \phi_{0,4}^3 \phi_{0,3}^7\\
{}&
+ 35785 \phi_{0,4}^6 \phi_{0,3}^3) \phi_{0,2}+
8144 \phi_{0,4}^2 \phi_{0,3}^9 + 17917 \phi_{0,4}^5 \phi_{0,3}^5
+ 8005 \phi_{0,4}^8 \phi_{0,3}\bigr] \phi_{0,1}\\
{}&
- 29 \phi_{0,4}^5 \phi_{0,2}^8
+ 162 \phi_{0,4}^4 \phi_{0,3}^2 \phi_{0,2}^7
-(358 \phi_{0,3}^4 +10464 \phi_{0,4}^3) \phi_{0,4}^3\phi_{0,2}^6\\
{}& 
+ (392 \phi_{0,3}^4 + 
45141 \phi_{0,4}^3)\phi_{0,4}^2 \phi_{0,3}^2 \phi_{0,2}^5
-(213 \phi_{0,3}^8 + 66918 \phi_{0,4}^3 \phi_{0,3}^4\\
{}&
+30811 \phi_{0,4}^6)\phi_{0,4} \phi_{0,2}^4 
+(46 \phi_{0,3}^{8} + 43947 \phi_{0,4}^3 \phi_{0,3}^4
+ 83053 \phi_{0,4}^6)  \phi_{0,3}^{2}\phi_{0,2}^3\\
{}&
-(14354  \phi_{0,3}^8  +64611 \phi_{0,4}^3 \phi_{0,3}^4 
+30093 \phi_{0,4}^6) \phi_{0,4}^2\phi_{0,2}^2\\
{}&
 + (3426 \phi_{0,3}^{8}  - 496 \phi_{0,4}^3 \phi_{0,3}^4
+ 37331 \phi_{0,4}^6)\phi_{0,4} \phi_{0,3}^{2} \phi_{0,2}
- 569 \phi_{0,3}^{12}\\
{}&
+ 3899 \phi_{0,4}^3 \phi_{0,3}^8
- 455 \phi_{0,4}^6 \phi_{0,3}^4  -9861 \phi_{0,4}^9 - 83 \xi_{0,36}^{(1)}\\
{}&
=q^{-1}[-144] + 24 +
(24r^{12} + 72r^{11} +\cdots)q\\
{}&+(4r^{18}[-36] + 144r^{16} + 672r^{15}+\cdots)q^2+
(144r^{20} + 1008r^{19}+\cdots)q^3\\
{}&+(24r^{24}+ 288r^{23}+\cdots)q^4+
(8r^{27}[-9] +216r^{26} + 3096r^{25}+\cdots )q^5\\
{}&+(72r^{29} + 1584r^{28} + 15720r^{27}+\cdots)q^6+
(9r^{32}[-16] + 288r^{31}\\
{}&+5304r^{30}+\cdots)q^7+
(672r^{33}+ 12096r^{32}+\cdots)q^8+O(q^9)
\endsplit
$$
where
$$
\split
D_{0,6}&=\bigl(\theta(\tau,z)/\eta(q)\bigr)^{12}=
-\phi_{0,1}^2\phi_{0,4}+9\phi_{0,1}\phi_{0,2}\phi_{0,3}
-8\phi_{0,2}^3-27\phi_{0,3}^2\\
{}& =q(r^6 - 12r^5 + 66r^4 - 220r^3 + 495r^2 - 792r + 924 -
\dots )+q^2(\dots)
\endsplit
$$
is the generator of the ideal of the weak Jacobi forms of weight $0$
without $q^{0}$-term   (see \cite{33}).
We have
$$
P(\M_0)=
\left(\smallmatrix
{1}&{2}&{0}\cr
{0}&{-1}&{0}\cr
{-1}&{0}&{1}\cr
{2}&{18}&{1}\cr
{5}&{27}&{1}\cr
{7}&{32}&{1}\cr
\endsmallmatrix\right);\ \ \
G(P(\M_0))=\left(\smallmatrix
{4}&{-2}&{-72}&{-36}&{-18}&{-8}\cr
{-2}&{1}&{0}&{-18}&{-27}&{-32}\cr
{-72}&{0}&{144}&{-72}&{-288}&{-432}\cr
{-36}&{-18}&{-72}&{36}&{-18}&{-72}\cr
{-18}&{-27}&{-288}&{-18}&{9}&{0}\cr
{-8}&{-32}&{-432}&{-72}&{0}&{16}\cr
\endsmallmatrix\right);
$$
$$
\overline{R}=
\left[\matrix
1,\ \overline{1}\\
1,\ \overline{17}\\
4,\ \overline{2}\\
4,\ \overline{34}\\
9,\ \overline{27}\\
16, \overline{32}\\
36, \overline{18}\\
144, \overline{0}
\endmatrix\right];
\hskip10pt
\text{Mul}(\overline{R},\,\xi)=
\pmatrix
-3  &0  & 0\\
-3  & 0 & 0\\
 0  & 1 & 0\\
 0  & 1 & 0\\
 4  & 1 & 12\\
 3  & 1 & 9\\
 1  & 0 & 4\\
 0  & 0 & 1
\endpmatrix .
$$

\subhead
2.6. The list of algebras of Theorem 2.1.1
\endsubhead
In Table 2 below we give the list of all Lorentzian Kac--Moody
algebras from Theorem 2.1.1. For each algebra from the list, the
infinite product part of its denominator identity is defined by the
infinite product $B_\xi$ of Theorem 2.2.1 for some Jacobi form $\xi$
from the space $RJ_t$ which is described in Table 1 of Sect. 2.5
by its bases. This product is characterized by the property that
its multiplicities are equal to 0 or 1
for any rational quadratic divisor which is orthogonal to
a root of $L_t$. Since $B_\xi$ is reflective, it is the whole
divisor of $B_\xi$. We denote the corresponding Lorentzian
Kac--Moody algebra as $\geg(\xi)$ since it is defined by
the Jacobi form $\xi$.

We also describe the Fourier expansion of
the automorphic form $B_\xi$, which gives the
infinite sum part of the denominator identity of the algebra
$\geg(\xi)$. For some of these forms only rational expressions
from known Fourier expansions are known. These results were
obtained in \cite{36}, \cite{37}, \cite{39}, \cite{41}.
We should say that these calculations are very non-trivial
because there does not exist a general method of finding
these Fourier expansions. We don't discuss these calculations
in this paper.

We describe the fundamental chamber $\M$ and
the set $P(\M)$ of orthogonal roots to $\M$ defining the Weyl
group and the set of simple real roots of the algebra
$\geg(\xi)$. We also give the subset $P(\M)_{\1o}\subset P(\M)$
of odd roots. It means that the corresponding
generators $e_\alpha$, $f_\alpha$, $\alpha \in P(\M)_{\1o}$,
should be super (odd). If we don't mention the set $P(\M)_{\1o}$,
it is empty. We also give the generalized Cartan matrix
$$
A=\left({2(\alpha_i,\alpha_j)\over \alpha_i^2}\right),\ \
\alpha_i,\,\alpha_j \in P(\M),
$$
which is the main invariant of the algebra. Many of these matrices
will be matrices of Theorem 1.5.6 or were considered in
\cite{41, Sect. 5.1}. Then we follow notations there.
We also give the Weyl vector $\rho$.

All these polygons $\M$ are composed from the fundamental polyhedron
$\M_0$ for $W(S_t)$ using some group of symmetries of the polyhedron $\M$.
We use these symmetries to describe the sets $P(\M)$ and $P(\M)_{\1o}$
using the set $P(\M_0)$.  We numerate elements 
$\alpha_1,\dots ,\alpha_k$ of $P(\M_0)$ as they are given in Table 1. 
We denote by $s_\alpha$ the reflection in the root $\alpha$. It is
given by the formula
$$
s_\alpha\,:\,x\to x-{2(x,\,\alpha )\over \alpha^2}\alpha,\ \
x\in S_t^\ast.
$$
We denote by $[g_1,\dots ,g_k]$ the group generated by
$g_1,\dots ,g_k$.

\vskip40pt

\centerline{\bf Table 2. The list of all Lorentzian Kac--Moody algebras}
\centerline{\bf with the root lattice $S_t^\ast$, symmetry lattice
$L_t^\ast$}
\centerline{\bf and the symmetry group $\widehat{O}(L_t)$}

\vskip20pt

\centerline{\bf Case $t=1$}

\vskip10pt

{\it The Algebra $\geg(\xi_{0,1}^{(1)})$.}
The fundamental chamber
$\M=[s_{\alpha_1},\,s_{\alpha_3}](\M_0)$
is the right triangle with zero angles. We have
$$
P(\M)=[s_{\alpha_1},\,s_{\alpha_3}](\alpha_2)=\{(0,-1,0),\
(1,1,0),\ (0,1,1)\}
$$
with the group of symmetries $[s_{\alpha_1},s_{\alpha_3}]$ which is
the dihedral group $D_3$ of order 6 (we use the same notation in general
for the dihedral group $D_n$). The generalized Cartan matrix is
$$
A_{1,II}=
\pmatrix
\hphantom{-}{2}&{-2}&{-2}\cr
{-2}&\hphantom{-}{2}&{-2}\cr
{-2}&{-2}&\hphantom{-}{2}\cr
\endpmatrix      .
$$
The Weyl vector
$\rho=({1\over 2},{1\over 2},{1\over 2})$.
The automorphic form $B_{\xi_{0,1}^{(1)}}$ coincides with
the classical automorphic form $\Delta_5$ of
the weight 5 which is product (divided by $64$) of
ten even theta-constants of genus 2. This automorphic form also
gives the discriminant of genus 2 algebraic curves. By Maass \cite{55},
$$
\Delta_5=\sum\Sb n,l,m \equiv 1\operatorname{mod}2\\
\vspace{0.5\jot} n,m>0\endSb \hskip10pt
\sum_{d|(n,l,m)} (-1)^{\frac{l+d+2}2}\,d^{4}
\tau_9\,(\frac{4nm-l^2}{d^2})\,
q^{n/2}\,r^{l/2}\,s^{m/2}\,,
$$
where $\eta(\tau)^9=\sum_{n\ge 1}{\tau_9(n)q^{n/24}}$ (see \S 3,
\thetag{3.4} about $\eta(\tau)$).
It gives the infinite sum part of the denominator identity of the algebra
$\geg(\xi_{0,1}^{(1)})$. Thus, the denominator identity of
the algebra $\geg(\xi_{0,1}^{(1)})$ has the form
$$
\align
&\sum\Sb n,l,m \equiv 1\operatorname{mod}2\\
\vspace{0.5\jot} n,m>0\endSb \hskip10pt
\sum_{a|(n,l,m)} (-1)^{\frac{l+a+2}2}\,a^{4}
\tau_9\,(\frac{4nm-l^2}{a^2})\,
q^{n/2}\,r^{l/2}\,s^{m/2}\\
\vspace{2\jot}
{}&=(qrs)^{1/2} \prod
\Sb n,\,l,\,m\in \Bbb Z\\
\vspace{0.5\jot}
(n,l,m)>0\,\endSb
\bigl(1-q^n r^l s^m\bigr)^{f_{1}(nm,l)}
\tag{2.6.1}
\endalign
$$
where
$$
\xi_{0,1}^{(1)}(\tau,\,z)=
\phi_{0,1}(\tau,\,z)=\sum_{k,l\in \bz} f_1(k,l)q^kr^l
$$
is Fourier expansion of the Jacobi form $\xi_{0,1}^{(1)}=\phi_{0,1}$
from Table 1, case t=1, and $\phi_{0,1}$ is described in \S 3,
\thetag{3.13}.
See \cite{36}, \cite{37}, \cite{41} for details in this case.

\smallpagebreak

Denominator identities like \thetag{2.6.1} can be similarly written
in all cases below (since Fourier expansions of the corresponding
automorphic forms $\xi$ and $B_\xi$ are known). We leave this
to a reader.

\smallpagebreak

{\it The Algebra $\geg(\xi_{0,1}^{(2)})$.}
The chamber $\M=\M_0$ is a triangle with
angles $0$, $\pi/2$, $\pi/3$. The set
$P(\M)=P(\M_0)$, $P(\M)_{\1o}=\{\alpha_2\}$. The generalized
Cartan matrix is
$$
A_{1,I,\0o}=A_{1,I,\1o}=
\pmatrix
2 & -1 & -1\\
-4&  2 &  0\\
-1&  0 &  2
\endpmatrix .
$$
The Weyl vector $\rho=({5\over 2},{1\over 2},{3\over 2})$.
The automorphic form $B_{\xi_{0,1}^{(2)}}$ coincides with
Igusa's \cite{46} modular form $\Delta_{30}=\Delta_{35}/\Delta_5$
of the weight 30. See \cite{39}, \cite{41}.
Fourier expansion of
$\Delta_{35}$ was found by Igusa in \cite{46}. Another expression for
Fourier expansion of $\Delta_{35}$ as Hecke product of $\Delta_5$
was found in \cite{39}. Let
$$
\align
[\Delta_5(z)]_{T(2)}=
{}&
\prod_{a,b,c \mod 2} 
\Delta_5 ({\tsize \frac{z_1+a}2 ,
\tsize\frac{z_2+b}2, \tsize\frac{z_3+c}2})\\
{}&
\times\prod_{a \mod 2}
\Delta_5 ({\tsize\frac{z_1+a}2,z_2,2z_3})\,
\Delta_5 ({\tsize 2z_1, z_2, \frac{z_3+a}2})\\
\vspace{2\jot}
{}&\times\Delta_5 ({\tsize 2z_1, 2z_2, 2z_3})
\prod_{b\mod 2}
\Delta_5 ({\tsize 2z_1, -z_1+z_2, \frac{z_1-2z_2+z_3+b}2}).
\endalign
$$
In \cite{39}, it was shown that
$$
\Delta_{35}(z)={[\Delta_5(z)]_{T(2)}\over \Delta_5(z)^8}\,.
$$
Thus, $\Delta_{30}(z)=[\Delta_5(z)]_{T(2)}/\Delta_5(z)^9$.
This gives Fourier expansions of $\Delta_{35}$ and $\Delta_{30}$
as finite products and quotients of known Fourier expansions.
See \cite{39} and \cite{41} for details.

{\it The Algebra $\geg(\xi_{0,1}^{(1)}+\xi_{0,1}^{(2)})$.}
The chamber $\M=\M_0$ (the same as for $\geg(\xi_{0,1}^{(2)})$) 
and
$$
P(\M)=\{\alpha_1,\,2\alpha_2,\,\alpha_3\}.
$$
The generalized Cartan matrix is
$$
A_{1,0}=
\pmatrix
2&-2&-1\\
-2&2&0\\
-1&0&2
\endpmatrix.
$$
The Weyl vector $\rho=(3,\,1,\,2)$. The automorphic form
$B_{\xi_{0,1}^{(1)}+\xi_{0,1}^{(2)}}$ coincides with
Igusa's \cite{46} modular form $\Delta_{35}$ of the weight 35 which
has been considered above.

\vskip20pt

\centerline{\bf Case $t=2$}

\vskip10pt

{\it The Algebra $\geg(\xi_{0,2}^{(1)})$.} The chamber
$\M=[s_{\alpha_1},\,s_{\alpha_3}](\M_0)$ is
the right quadrangle with zero angles; the set 
$$
P(\M)=[s_{\alpha_1},\,s_{\alpha_3}](\alpha_2)=
\{(0,-1,0),\ (1,1,0),\ (1,3,1),\ (0,1,1)\}
$$
with the group of symmetries $[s_{\alpha_1},\,s_{\alpha_3}]$
which is $D_4$. The generalized Cartan matrix is
$$
A_{2,II}=
\pmatrix
\hphantom{-}{2}&{-2}&{-6}&{-2}\cr
{-2}&\hphantom{-}{2}&{-2}&{-6}\cr
{-6}&{-2}&\hphantom{-}{2}&{-2}\cr
{-2}&{-6}&{-2}&\hphantom{-}{2}\cr
\endpmatrix.
$$
The Weyl vector $\rho=({1\over 4},\,{1\over 2},\,{1\over 4})$.
The automorphic form $B_{\xi_{0,2}^{(1)}}$ coincides with
the automorphic form $\Delta_2$ of the weight 2 which was
introduced in \cite{36} and \cite{41}.
Its Fourier expansion is
$$
\Delta_2=\sum_{N\ge 1}\
\sum\Sb
 n,\,m >0,\,l\in \bz\\
\vspace{0.5\jot} n,\,m\equiv 1\mod 4\\
\vspace{0.5\jot} 2nm-l^2=N^2
\endSb
\hskip-4pt
N\biggl(\frac {-4}{Nl}\biggr)
\sum_{a\,|\,(n,l,m)} \biggl(\frac {-4}{a}\biggr)
\, q^{n/4} r^{l/2} s^{m/4}\,.
$$
Here $\left(\frac {\,m\,}{\,n\,}\right)$ is the Jacobi symbol
(or the generalized Legendre symbol).

{\it The Algebra $\geg(\xi_{0,2}^{(2)})$.}
The chamber $\M=[s_{\alpha_3}](\M_0)$ is a triangle with
angles $0,\,0,\,\pi/2$. The sets
$$
P(\M)=\{\alpha_1,\,\alpha_2,\,s_{\alpha_3}(\alpha_1)=(0,2,1)\},\ \ \
P(\M)_{\1o}=\{\alpha_2\}
$$
with the group of symmetries $[s_{\alpha_3}]$ which is $D_1$.
The generalized Cartan matrix is
$$
A_{2,I,\0o}=
\pmatrix
2 & -1 &  0\\
-4&  2 & -4\\
0 & -1 &  2
\endpmatrix.
$$
The Weyl vector $\rho=({3\over 4},{1\over 2},{3\over 4})$.
The automorphic form $B_{\xi_{0,2}^{(2)}}$ coincides with the
automorphic form $\Delta_{9}=\Delta_{11}/\Delta_2$ of the weight
9 from \cite{41}.
There are two formulae for the Fourier expansion of $\Delta_{11}$.
This automorphic form is the lifting of its first Fourier-Jacobi
coefficient :
$$
\Delta_{11}(z)=\hbox{Lift}(\eta(z_1)^{21}\vartheta(z_1,2z_2)).
$$
It gives us a simple exact formula for the Fourier coefficients of
$\Delta_{11}$ in terms of the Fourier coefficients of the Jacobi form
$\eta(z_1)^{21}\vartheta(z_1 ,2z_2)$.
This formula is  similar to the  formula for  $\Delta_5$ above
(see \cite{41, Example 1.15}).

The second expression for $\Delta_{11}$ is given by
the multiplicative symmetrisation of $\Delta_5$ for $t=1$ above. Let
$$
\text{Ms}_2(\Delta_5)(z_1,z_2,z_3)=
\Delta_5(z_1,2z_2,4z_3)\Delta_5(z_1,z_2,z_3)\Delta_5(z_1,z_2,z_3+1).
$$
By \cite{41, \thetag{3.10}},
$$
\Delta_{11}(z)={\text{Ms}_2(\Delta_5)(z)\over \Delta_2(z)^2}.
$$
Thus, $\Delta_9=\text{Ms}_2(\Delta_5)(z)/\Delta_2(z)^3$.
It gives Fourier expansions of $\Delta_9$
as finite products and quotients of known Fourier expansions.
See \cite{39}, \cite{41} for details.

{\it The Algebra $\geg(\xi_{0,2}^{(1)}+\xi_{0,2}^{(2)})$.} The
polygon $\M=[s_{\alpha_3}](\M_0)$
is a triangle with angles $0,\,0,\,\pi/2$
(the same as for $\geg(\xi_{0,2}^{(2)})$); the set
$$
P(\M)=\{\alpha_1,\,2\alpha_2,\,s_{\alpha_3}(\alpha_1)=(0,2,1)\}
$$
with the group of symmetries $[s_{\alpha_3}]$ which is $D_1$.
The generalized Cartan matrix is
$$
A_{2,0}=
\pmatrix
2 & -2& 0\\
-2&  2&-2\\
 0&	-2&	2
\endpmatrix.
$$
The Weyl vector $\rho=(1,1,1)$. The automorphic form
$B_{\xi_{0,2}^{(1)}+\xi_{0,2}^{(2)}}$ coincides with the
$\Delta_{11}$ of the weight 11. We have discussed its Fourier
expansion above. See \cite{39}, \cite{41}.

{\it The Algebra $\geg(\xi_{0,2}^{(3)})$.} 
The chamber $\M=[s_{\alpha_1},\,s_{\alpha_2}](\M_0)$ is an
infinite polygon with angles $\pi/2$ and
which is touching a horosphere with the
center at $\br_{++}\rho$ where $\rho=(1,0,0)$ is
the Weyl vector. The set
$$
P(\M)=[s_{\alpha_1},\,s_{\alpha_2}](\alpha_3)
$$
with the group of symmetries $[s_{\alpha_1},\,s_{\alpha_2}]$
which is $D_\infty$. The generalized
Cartan matrix is the symmetric matrix
$$
A_{2,\1o}=
\left({(\alpha,\,\alpha^\prime)\over 4}\right),\ \
\alpha,\alpha^\prime \in P(\M).
$$
The automorphic form
$B_{\xi_{0,2}^{(3)}}$ is $\Psi_{12}^{(2)}$ of the
weight 12 from \cite{41}.

It was shown in \cite{41, Remark 4.4}
that $\Psi_{12}^{(2)}$ can be obtained as restriction (possibly with
some multiplicative constant) of the
Borcherds automorphic form $\Phi$ from \thetag{1.3.7}.
The Borcherds automorphic form $\Phi$ is defined on Hermitian symmetric
domain $\Omega(2H\oplus \La)$ where $\La$ is Leech lattice.
One should restrict $\Phi$
on the subdomain $\Omega(2H+\bz v)$ where $v\in \La$ is a primitive
element with $v^2=2t$ where $t=2$ for this case. Thus, we have
$$
\Psi_{12}^{(2)}=c\Phi|\Omega(2H+\bz v)
$$
where $c$ is some constant.
It gives some Fourier expansion of $\Psi_{12}^{(2)}$.
The same construction is valid for automorphic forms
$\Psi_{12}^{(t)}$ which we consider below when $t=3$ and $t=4$.

{\it The Algebra $\geg(\xi_{0,2}^{(1)}+\xi_{0,2}^{(3)})$.} 
The polygon $\M=[s_{\alpha_1}](\M_0)$ is a quadrangle with angles
$\pi/2,\,\pi/2,\,\pi/2,\,0$; the set
$$
P(\M)=[s_{\alpha_1}]\{\alpha_2,\,\alpha_3\}=
\{\alpha_2,\alpha_3,(1,4,1),(1,1,0)\}
$$
with the group of symmetries $[s_{\alpha_1}]$ which is $D_1$.
The generalized Cartan matrix is
$$
A_{2,II,\1o}=
\pmatrix
{2}&{0}&{-8}&{-2}\cr
{0}&{2}&{0}&{-1}\cr
{-1}&{0}&{2}&{0}\cr
{-2}&{-8}&{0}&{2}\cr
\endpmatrix.
$$
The Weyl vector
$\rho=({{5}\over{4}},{{1}\over{2}},{{1}\over{4}})$.
The automorphic form  $B_{\xi_{0,2}^{(1)}+\xi_{0,2}^{(3)}}$ coincides
with the automorphic form
$\Delta_{14}=\Delta_2 \Psi_{12}^{(2)}$ of the weight 14.
(We must correct the case $(2,II,\1o)$ in \cite{41, page 264}
in this way.) Fourier expansion of $\Delta_{14}$ is product of
the Fourier expansions of $\Delta_2$ and $\Psi_{12}^{(2)}$.

{\it The Algebra $\geg(\xi_{0,2}^{(2)}+\xi_{0,2}^{(3)})$.}
The polygon $\M=\M_0$ is the triangle with angles
$0,\,\pi/2,\pi/4$; the set
$$
P(\M)=P(\M_0),\ \ \  P(\M)_{\1o}=\{\alpha_2\}
$$
with the trivial group of symmetries and
with the generalized Cartan matrix
$$
A_{2,I,\1o}=
\pmatrix
2 & -1 & -2\\
-4&  2 &  0\\
-1&  0 &  2
\endpmatrix.
$$
The Weyl vector
$\rho=({{7}\over{4}},{{1}\over{2}},{{3}\over{4}})$.
The automorphic form $B_{\xi_{0,2}^{(2)}+\xi_{0,2}^{(3)}}$
coincides with the automorphic form $\Delta_9\Psi_{12}^{(2)}$ of
the weight $21$. Its Fourier expansion is product of the
Fourier expansions of $\Delta_9$ and $\Psi_{12}^{(2)}$.

{\it The Algebra $\geg(\xi_{0,2}^{(1)}+\xi_{0,2}^{(2)}+\xi_{0,2}^{(3)})$.}
The polygon $\M=\M_0$ is the triangle with angles
$0,\pi/2,\,\pi/4$ (the same as for the
$\geg(\xi_{0,2}^{(2)}+\xi_{0,2}^{(3)})$); the set
$$
P(\M)=\{\alpha_1\,,2\alpha_2\,,\alpha_3\}
$$
with the trivial group of symmetries and
with the generalized Cartan matrix
$$
A_{2,0,\1o}=
\pmatrix
2 & -2 & -2\\
-2&  2 &  0\\
-1&  0 &  2
\endpmatrix.
$$
The Weyl vector $\rho=(2,\,1,\,1)$.
The automorphic form
$B_{\xi_{0,2}^{(1)}+\xi_{0,2}^{(2)}+\xi_{0,2}^{(3)}}$
coincides with the automorphic form
$\Delta_2\Delta_9\Psi_{12}^{(2)}$ of the weight $23$.
Its Fourier expansion is given by products of the Fourier expansions
of $\Delta_2$, $\Delta_9$  and $\Psi_{12}^{(2)}$.
See \cite{41}.

\vskip20pt

\centerline{\bf Case $t=3$}

\vskip10pt

{\it The Algebra $\geg(\xi_{0,3}^{(1)})$.} The chamber
$\M=[s_{\alpha_1},\,s_{\alpha_3}](\M_0)$ is
the right hexagon with zero angles,
$$
P(\M)=[s_{\alpha_1},\,s_{\alpha_3}](\alpha_2)
=
\{(0,-1,0),\ (1,1,0),\ (2,5,1),\ (2,7,2),\ (1,5,2),\ (0,1,1)\}
$$
with the group of symmetries $[s_{\alpha_1},\,s_{\alpha_3}]$
which is $D_6$. The generalized Cartan matrix is
$$
A_{3,II}=
\pmatrix
\hphantom{-}{2}&{-2}&{-10}&{-14}&{-10}&{-2}\cr
{-2}&\hphantom{-}{2}&{-2}&{-10}&{-14}&{-10}\cr
{-10}&{-2}&\hphantom{-}{2}&{-2}&{-10}&{-14}\cr
{-14}&{-10}&{-2}&\hphantom{-}{2}&{-2}&{-10}\cr
{-10}&{-14}&{-10}&{-2}&\hphantom{-}{2}&{-2}\cr
{-2}&{-10}&{-14}&{-10}&{-2}&\hphantom{-}{2}\cr
\endpmatrix .
$$
The Weyl vector $\rho=({1\over 6},\,{1\over 2},\,{1\over 6})$.
The automorphic form  $B_{\xi_{0,3}^{(1)}}$ coincides with
the automorphic form $\Delta_1$ of the weight 1 introduced in
\cite{41}. Its Fourier expansion is
$$
\Delta_1=\sum_{M\ge 1}
\sum\Sb n,\,m >0,\,l\in \bz\\
\vspace{0.5\jot} n,\,m\equiv 1\mod 6\\
\vspace{0.5\jot} 4nm-3l^2=M^2\endSb
\hskip-4pt
\biggl(\dsize\frac{-4}{l}\biggr)
\biggl(\dsize\frac{12}{M}\biggr)
\sum\Sb a|(n,l,m)\endSb \biggl(\dsize\frac{6}{a}\biggr)
q^{n/6}r^{l/2}s^{m/6}\,.
$$
See \cite{39}, \cite{41}.

{\it The Algebra $\geg(\xi_{0,3}^{(2)})$.}
The chamber $\M=[s_{\alpha_3}](\M_0)$ is a triangle with
angles $0,\,0,\,\pi/3$. The sets
$$
P(\M)=\{\alpha_1,\,\alpha_2,\,s_{\alpha_3}(\alpha_1)=(0,2,1)\},\ \ \
P(\M)_{\1o}=\{\alpha_2\}
$$
with the group of symmetries $[s_{\alpha_3}]$ which is $D_1$.
The generalized Cartan matrix is
$$
A_{3,I,\0o}=
\pmatrix
2 & -1 & -1\\
-4&  2 & -4\\
-1& -1 &  2
\endpmatrix .
$$
The Weyl vector $\rho=({1\over 2},{1\over 2},{1\over 2})$.
The automorphic form $B_{\xi_{0,3}^{(2)}}$ coincides with
the automorphic form $D_6$ of the weight 6 introduced in
\cite{41}.
The form $D_6$ is the  lifting of its first Fourier-Jacobi
coefficient $\eta(z_1)^{11}\vartheta_{3/2}(z_1, z_2)$ where
$\vartheta_{3/2}(\tau, z)=
\eta(\tau)\vartheta(\tau,2z)/\vartheta(\tau,z)$.
Thus there is an exact formula for the Fourier coefficients
of $D_6$ in terms of the Fourier coefficients of this
Jacobi form (see \cite{41, Example 1.17}).

Fourier expansion of $\Delta_1D_{6}$ is also given by the
finite Hecke product of $\Delta_1$ (which is similar to the
Fourier expansion of $\Delta_{35}$ for $t=1$ above).
Let
$$
\align
[\Delta_1(z)]_{T(2)}&=
\prod_{a,b,c \mod 2}
\Delta_1 ({\tsize \frac{z_1+a}2 ,
\tsize\frac{z_2+b}2, \tsize\frac{z_3+c}2})\\
{}&\times 
\prod_{a \mod 2}
\Delta_1 ({\tsize\frac{z_1+a}2,z_2,2z_3})\,
\Delta_1 ({\tsize 2z_1, z_2, \frac{z_3+a}2})\\
\vspace{2\jot}
{}&\times\Delta_1 ({\tsize 2z_1, 2z_2, 2z_3})
\prod_{b\mod 2}
\Delta_1 ({\tsize 2z_1, -z_1+z_2, \frac{z_1-2z_2+z_3+b}2}).
\endalign
$$
In \cite{41}, it is shown that
$$
\Delta_1(z)D_{6}(z)={2^{22}[\Delta_1(z)]_{T(2)}\over \Delta_1(z)^8}\ .
$$
Thus,  Fourier expansion of $\Delta_1D_6$ is given by finite products
and quotients of known Fourier expansions.
See \cite{39}, \cite{41} for details.

{\it The Algebra $\geg(\xi_{0,3}^{(1)}+\xi_{0,3}^{(2)})$.} The
polygon $\M=[s_{\alpha_3}](\M_0)$
is a triangle with angles $0,\,0,\,\pi/3$
(the same as for $\xi_{0,3}^{(2)}$); the set
$$
P(\M)=\{\alpha_1,\,2\alpha_2,\,s_{\alpha_3}(\alpha_1)=(0,2,1)\}
$$
with the group of symmetries $[s_{\alpha_3}]$ which is $D_2$.
The generalized Cartan matrix is
$$
A_{3,0}=
\pmatrix
2 & -2& -1\\
-2&  2&-2\\
-1&	-2&	2
\endpmatrix.
$$
The Weyl vector $\rho=({2\over 3},1,{2\over 3})$.
The automorphic form $B_{\xi_{0,3}^{(1)}+\xi_{0,3}^{(2)}}$
is $\Delta_1 D_6$ of the weight 7.
We have discussed Fourier expansion of
$\Delta_1 D_6$ above. See \cite{39}, \cite{41}.

{\it The Algebra $\geg(\xi_{0,3}^{(3)})$.} 
The chamber $\M=[s_{\alpha_1},\,s_{\alpha_2}](\M_0)$ is an
infinite polygon with angles $\pi/3$ and
which is touching a horosphere with the
center at the Weyl vector $\rho=(1,0,0)$. The set
$$
P(\M)=[s_{\alpha_1},\,s_{\alpha_2}](\alpha_3)
$$
with the group of symmetries $[s_{\alpha_1},\,s_{\alpha_2}]$
which is $D_\infty$. The generalized
Cartan matrix is the symmetric matrix
$$
A_{3,\1o}=
\left({(\alpha,\,\alpha^\prime)\over 6}\right),\ \
\alpha,\alpha^\prime \in P(\M).
$$
The automorphic form
$B_{\xi_{0,3}^{(3)}}$ coincides with the automorphic form
$\Psi_{12}^{(3)}$ of the weight 12 from \cite{41}.
Like $\Psi_{12}^{(2)}$ above,
the automorphic form $\Psi_{12}^{(3)}$ and its Fourier expansion
can be obtained as the restriction of Borcherds automorphic form
$\Phi$ from \thetag{1.3.7}. See \cite{41, Remark 4.4}.

{\it The Algebra $\geg(\xi_{0,3}^{(1)}+\xi_{0,3}^{(3)})$.} 
The polygon $\M=[s_{\alpha_1}](\M_0)$ is a quadrangle with angles
$0,\pi/2,\pi/3,\,\pi/2$; the set
$$
P(\M)=[s_{\alpha_1}]\{\alpha_2,\,\alpha_3\}=
\{\alpha_2,\alpha_3,(2,6,1),(1,1,0)\}
$$
with the group of symmetries $[s_{\alpha_1}]$ which is $D_1$.
The generalized Cartan matrix is
$$
A_{3,II,\1o}=
\pmatrix
{2}&{0}&{-12}&{-2}\cr
{0}&{2}&{-1}&{-1}\cr
{-1}&{-1}&{2}&{0}\cr
{-2}&{-12}&{0}&{2}\cr
\endpmatrix .
$$
The Weyl vector
$\rho=({{7}\over{6}},{{1}\over{2}},{{1}\over{6}})$.
The automorphic form $B_{\xi_{0,3}^{(1)}+\xi_{0,3}^{(3)}}$ coincides with
$\Delta_1 \Psi_{12}^{(3)}$ of the weight 13.
(We must correct the case $(3,II,\1o)$ in \cite{41, page 264}
in this way.) Its Fourier expansion is product of the
Fourier expansions of $\Delta_1$ and $\Psi_{12}^{(3)}$.

{\it The Algebra $\geg(\xi_{0,3}^{(2)}+\xi_{0,3}^{(3)})$.} 
The polygon $\M=\M_0$ is the triangle with angles
$0,\,\pi/2,\pi/6$; the set
$$
P(\M)=P(\M_0),\ \ \  P(\M)_{\1o}=\{\alpha_2\}
$$
with the trivial group of symmetries and
with the generalized Cartan matrix
$$
A_{3,I,\1o}=
\pmatrix
2 & -1 & -3\\
-4&  2 &  0\\
-1&  0 &  2
\endpmatrix.
$$
The Weyl vector
$\rho=({{3}\over{2}},{{1}\over{2}},{{1}\over{2}})$.
The automorphic form $B_{\xi_{0,3}^{(2)}+\xi_{0,3}^{(3)}}$
coincides with $D_6 \Psi_{12}^{(3)}$ of the weight $18$.
See \cite{41}. Its Fourier expansion is product of the
Fourier expansions of $D_6$ and $\Psi_{12}^{(3)}$.

{\it The Algebra $\geg(\xi_{0,3}^{(1)}+\xi_{0,3}^{(2)}+\xi_{0,3}^{(3)})$.} 
The polygon $\M=\M_0$ is the triangle with angles
$0,\pi/2,\,\pi/6$ (the same as for the
$\geg(\xi_{0,3}^{(2)}+\xi_{0,3}^{(3)})$); the set
$$
P(\M)=\{\alpha_1,\,2\alpha_2,\,\alpha_3\}
$$
with the trivial group of symmetries and
with the generalized Cartan matrix
$$
A_{3,0,\1o}=
\pmatrix
2 & -2 & -3\\
-2&  2 &  0\\
-1&  0 &  2
\endpmatrix.
$$
The Weyl vector $\rho=({5\over 3},\,1,\,{2\over 3})$.
The automorphic form
$B_{\xi_{0,3}^{(1)}+\xi_{0,3}^{(2)}+\xi_{0,3}^{(3)}}$ is
$\Delta_1 D_6 \Psi_{12}^{(3)}$ of the weight $19$.
See \cite{41}. Its Fourier expansion is product of the Fourier
expansions of $\Delta_1$, $D_6$ and $\Psi_{12}^{(3)}$.

\vskip20pt

\centerline{\bf Case $t=4$}

\vskip10pt

{\it The Algebra $\geg(\xi_{0,4}^{(1)})$.} The chamber
$\M=[s_{\alpha_1},\,s_{\alpha_3}](\M_0)$ is
the infinite polygon with zero angles touching a horosphere with the
center at $\br_{++}\rho$,  where
$\rho=({1\over 8},{1\over 2},{1\over 8})$ is the Weyl vector; the set 
$$
P(\M)=[s_{\alpha_1},\,s_{\alpha_3}](\alpha_2)
$$
with the group of symmetries $[s_{\alpha_1},\,s_{\alpha_3}]$
which is $D_\infty$. The generalized Cartan matrix is
$$
A_{4,II,\0o}=\left(2(\alpha,\alpha^\prime)\right),\ \ \
\alpha,\,\alpha^\prime \in P(\M).
$$
The automorphic form $B_{\xi_{0,4}^{(1)}}$ coincides with
$\Delta_{1/2}$ of the weight 1/2 which is the theta-constant
of the genus 2. Its Fourier expansion is
$$
\Delta_{1/2}=\frac{1}2\sum\Sb n,m\in \bz \endSb
\,\biggl(\dsize\frac{-4}{n}\biggr)\biggl(\dsize\frac{-4}{m}\biggr)
q^{n^2/8}r^{nm/2}s^{m^2/8}\,.
$$
See \cite{41}.

{\it The Algebra $\geg(\xi_{0,4}^{(2)})$.} 
The chamber $\M=[s_{\alpha_3}](\M_0)$ is a triangle with
angles $0,\,0,\,0$. The sets
$$
P(\M)=\{\alpha_1,\,\alpha_2,\,s_{\alpha_3}(\alpha_1)=(0,2,1)\};\ \ \
P(\M)_{\1o}=\{\alpha_2\}
$$
with the group of symmetries $[s_{\alpha_3}]$ which is $D_1$.
The generalized Cartan matrix is
$$
A_{4,I,\0o}=
\pmatrix
2 & -1 &  -2\\
-4&  2 & -4\\
-2 & -1 &  2
\endpmatrix .
$$
The Weyl vector $\rho=({3\over 8},{1\over 2},{3\over 8})$.
The automorphic form $B_{\xi_{0,4}^{(2)}}$ coincides with 
\newline 
$\Delta_5^{(4)}/\Delta_{1/2}$ of the weight ${9\over 2}$ where
$\Delta_5^{(4)}=\Delta_5(z_1,2z_2,z_3)$ and $\Delta_5(z)$ was used
for $t=1$. Thus, $\Delta_5^{(4)}$ has Fourier expansion
$$
\Delta_5^{(4)}=
\sum\Sb n,l,m \equiv 1\operatorname{mod}2\\
\vspace{0.5\jot} n,m>0\endSb \hskip10pt
\sum_{d|(n,l,m)} (-1)^{\frac{l+d+2}2}\,d^{4}
\,\tau_9\left(\frac{4nm-l^2}{d^2}\right)\,
q^{n/2}\,r\,s^{m/2}\,
$$
and Fourier expansion of $\Delta_5^{(4)}/\Delta_{1/2}$
is quotient of the Fourier expansions of $\Delta_5^{(4)}$
and $\Delta_{1/2}$. See \cite{41}.

{\it The Algebra $\geg(\xi_{0,4}^{(1)}+\xi_{0,4}^{(2)})$.} The
polygon $\M=[s_{\alpha_3}](\M_0)$
is a triangle with angles $0,\,0,\,0$
(the same as for $\geg(\xi_{0,1}^{(1)})$ and $\geg(\xi_{0,4}^{(2)})$); 
the set 
$$
P(\M)=\{\alpha_1,\,2\alpha_2,\,s_{\alpha_3}(\alpha_1)=(0,2,1)\}
$$
with the group of symmetries $[s_{\alpha_3}]$ which is $D_1$.
The generalized Cartan matrix is
$$
A_{4,0,\0o}=A_{1,II}=
\pmatrix
2 & -2& -2\\
-2&  2&-2\\
-2& -2& 2
\endpmatrix
$$
(it is the same as for $\geg(\xi_{0,1}^{(1)})$).
The Weyl vector $\rho=({1\over 2},1,{1\over 2})$.
The automorphic form $B_{\xi_{0,4}^{(1)}+\xi_{0,4}^{(2)}}$
coincides with
$\Delta_5^{(4)}(z_1,z_2,z_3)=\Delta_5(z_1,2z_2,z_3)$
of the weight 5 with Fourier expansion above.
This case is equivalent to the case $\geg(\xi_{0,1}^{(1)})$ above.
See \cite{41}.

{\it The Algebra $\geg(\xi_{0,4}^{(3)})$.} 
The chamber $\M=[s_{\alpha_1},\,s_{\alpha_2}](\M_0)$ is an
infinite polygon with zero angles
which is touching a horosphere with the
center at $\br_{++}\rho$, where $\rho=(1,0,0)$ is the Weyl
vector. The set
$$
P(\M)=[s_{\alpha_1},\,s_{\alpha_2}](\alpha_3)
$$
with the group of symmetries $[s_{\alpha_1},\,s_{\alpha_2}]$
which is $D_\infty$. The generalized
Cartan matrix is the symmetric matrix
$$
A_{4,\1o}=
\left({(\alpha,\,\alpha^\prime)\over 8}\right),\ \
\alpha,\alpha^\prime \in P(\M).
$$
The automorphic form $B_{\xi_{0,4}^{(3)}}$ coincides with
the automorphic form
$\Psi_{12}^{(4)}$ of the weight 12 from \cite{41}.
Like $\Psi_{12}^{(2)}$ above,
the automorphic form $\Psi_{12}^{(4)}$ and its Fourier expansion
can be obtained as restriction of the Borcherds automorphic form
$\Phi$ from \thetag{1.3.7}. See \cite{41, Remark 4.4}.

{\it The Algebra $\geg(\xi_{0,4}^{(1)}+\xi_{0,4}^{(3)})$.} 
The polygon $\M=[s_{\alpha_1}](\M_0)$ is a quadrangle with angles
$\pi/2,\,0,\,\pi/2,\,0$; the set
$$
P(\M)=[s_{\alpha_1}]\{\alpha_2,\,\alpha_3\}=
\{\alpha_2,\alpha_3,(3,8,1),(1,1,0)\}
$$
with the group of symmetries $[s_{\alpha_1}]$ which is $D_1$.
The generalized Cartan matrix is
$$
A_{4,II,\1o}=
\pmatrix
{2}&{0}&{-16}&{-2}\cr
{0}&{2}&{-2}&{-1}\cr
{-1}&{-2}&{2}&{0}\cr
{-2}&{-16}&{0}&{2}\cr
\endpmatrix .
$$
The Weyl vector
$\rho=({9\over 8},{1\over 2},{1\over 8})$.
The automorphic form $B_{\xi_{0,4}^{(1)}+\xi_{0,4}^{(3)}}$
coincides with
$\Delta_{1/2} \Psi_{12}^{(4)}$ of the weight ${25\over 2}$.
(We must correct the case $(4,II,\1o)$ in \cite{41, page 264}
in this way.) Its Fourier expansion is product of the Fourier
expansions of $\Delta_{1/2}$ and $\Psi_{12}^{(4)}$.

{\it The Algebra $\geg(\xi_{0,4}^{(2)}+\xi_{0,4}^{(3)})$.} 
The polygon $\M=\M_0$ is the triangle with angles
$0,\,\pi/2,\,0$; the sets 
$$
P(\M)=P(\M_0),\ \ \  P(\M)_{\1o}=\{\alpha_2\}
$$
with the trivial group of symmetries and with the generalized
Cartan matrix
$$
A_{4,I,\1o}=
\pmatrix
2 & -1 & -4\\
-4&  2 &  0\\
-1&  0 &  2
\endpmatrix.
$$
The Weyl vector
$\rho=({{11}\over{8}},{{1}\over{2}},{{3}\over{8}})$.
The automorphic form $B_{\xi_{0,4}^{(2)}+\xi_{0,4}^{(3)}}$
coincides with
$\Psi_{12}^{(4)} \Delta_5^{(4)}/\Delta_{1/2}$ of the weight
${33\over 2}$. Its Fourier expansion is product and quotient of
known Fourier expansions of $\Psi_{12}^{(4)}$,
$\Delta_5^{(4)}$ and $\Delta_{1/2}$. See \cite{41}.

{\it The Algebra $\geg(\xi_{0,4}^{(1)}+\xi_{0,4}^{(2)}+\xi_{0,4}^{(3)})$.} 
The polygon $\M=\M_0$ is the triangle with angles
$0,\pi/2,\,0$ (the same as for the
$\geg(\xi_{0,4}^{(2)}+\xi_{0,4}^{(3)})$); the set
$$
P(\M)=\{\alpha_1,\,2\alpha_2,\alpha_3\}
$$
with the trivial group of symmetries and
with the generalized Cartan matrix
$$
A_{4,0,\1o}=
\pmatrix
2 & -2 & -4\\
-2&  2 &  0\\
-1&  0 &  2
\endpmatrix.
$$
The Weyl vector $\rho=({3\over 2},\,1,\,{1\over 2})$.
The automorphic form
$B_{\xi_{0,4}^{(1)}+\xi_{0,4}^{(2)}+\xi_{0,4}^{(3)}}$ coincides with
$\Delta_5^{(4)}\Psi_{12}^{(4)}$ of the weight $17$. Its Fourier
expansion is product of the Fourier expansions of
$\Delta_5^{(4)}$ and $\Psi_{12}^{(4)}$. See \cite{41}.

\vskip20pt

\centerline{\bf Case $t=8$}

\vskip10pt

{\it The Algebra $\geg(\xi_{0,8}^{(2)})$.} The chamber
$\M=[s_{\alpha_3}](\M_0)$ is the right quadrangle with
zero angles; the set
$$
P(\M)=[s_{\alpha_3}](\alpha_1,\,2\alpha_2,\,\alpha_4)=
\{\alpha_1,\,2\alpha_2,\,s_{\alpha_3}(\alpha_1)=(0,2,1),\,\alpha_4 \}
$$
with the group of symmetries $[s_{\alpha_3}]$ which is $D_1$.
The generalized Cartan matrix is
$A_{2,II}$ (the same as for $\geg(\xi_{0,2}^{(1)})$ for $t=2$).
The Weyl vector $\rho=({1\over 4},\,1,\,{1\over 4})$.
The automorphic form $B_{\xi_{0,8}^{(2)}}$ coincides with
$\Delta_2^{(8)}(z_1,z_2,z_3)=\Delta_2(z_1,2z_2,z_3)$ of the
weight 2 where $\Delta_2$ corresponds to $\geg(\xi_{0,2}^{(1)})$.
Thus, Fourier expansion of $\Delta_2^{(8)}$ is
$$
\Delta_2^{(8)}=\sum_{N\ge 1}\
\sum\Sb
 n,\,m >0,\,l\in \bz\\
\vspace{0.5\jot} n,\,m\equiv 1\mod 4\\
\vspace{0.5\jot} 2nm-l^2=N^2
\endSb
\hskip-4pt
N\biggl(\frac {-4}{Nl}\biggr)
\sum_{a\,|\,(n,l,m)} \biggl(\frac {-4}{a}\biggr)
\, q^{n/4} r\, s^{m/4}\,.
$$
This case is equivalent to $\geg(\xi_{0,2}^{(1)})$.

\vskip20pt

\centerline{\bf Case $t=9$}

\vskip10pt

{\it The algebra $\geg(\xi_{0,9}^{(2)})$.} The chamber
$\M=[s_{\alpha_3}](\M_0)$ is the pentagon with angles
$0$, $0$, $\pi/2$, $0$, $\pi/2$; the set
$$
\split
&P(\M)=[s_{\alpha_3}](\alpha_1,\alpha_2,\alpha_4)=
\{\alpha_1,\,\alpha_2,\,s_{\alpha_3}(\alpha_1)=(0,2,1),\,
s_{\alpha_3}(\alpha_4)=(1,9,2),\,
\alpha_4\},\\
&P(\M)_{\1o}=\{\alpha_2\}
\endsplit
$$
with the group of symmetries
$[s_{\alpha_3}]$ which is $D_1$.
The generalized Cartan matrix is
$$
\pmatrix
{2}&{-1}&{-7}&{-9}&{0}\cr
{-4}&{2}&{-4}&{-18}&{-18}\cr
{-7}&{-1}&{2}&{0}&{-9}\cr
{-4}&{-2}&{0}&{2}&{-2}\cr
{0}&{-2}&{-4}&{-2}&{2}\cr
\endpmatrix .
$$
The Weyl vector $\rho=({1\over 6},{1\over 2},{1\over 6})$.
The automorphic form $B_{\xi_{0,9}^{(2)}}$
coincides with the automorphic form $D_2$ of the weight 2 with Fourier
expansion
$$
D_2=\sum_{N\ge 1}
\sum\Sb m >0,\,l\in \bz\\
\vspace{0.5\jot} n,\,m\equiv 1\mod 6\\
\vspace{0.5\jot} 4nm-l^2=N^2\endSb
\hskip-4pt
N\biggl(\dsize\frac{-4}{N}\biggr)\biggl(\dsize\frac{12}{l}\biggr)
\sum\Sb a|(n,l,m)\endSb \biggl(\dsize\frac{\,6\,}{\,a\,}\biggr)
q^{n/6}r^{l/2}s^{m/6}\,.
$$
See \cite{41, \thetag{5.1.2}}.

\vskip20pt

\centerline{\bf Case $t=12$}

\vskip10pt

{\it The Algebra $\geg(\xi_{0,12}^{(2)}+\xi_{0,12}^{(3)})$.} 
The chamber $\M=[s_{\alpha_3},s_{\alpha_4}](\M_0)$ is the
right hexagon with zero angles; the set
$$
\split
P(\M)=&[s_{\alpha_3},s_{\alpha_4}](\alpha_1,\,2\alpha_2)=
\{\alpha_1,\,2\alpha_2,\,s_{\alpha_3}(\alpha_1)=(0,2,1),\\
&s_{\alpha_4}s_{\alpha_3}(\alpha_1)=(1,10,2),
s_{\alpha_4}(2\alpha_2)=(2,14,2),\,
s_{\alpha_4}(\alpha_1)=(2,10,1)\}
\endsplit
$$
with the group of symmetries $[s_{\alpha_3},\,s_{\alpha_4}]$ which
is $D_4$. The generalized Cartan matrix is
$A_{3,II}$ (the same as for $\geg(\xi_{0,3}^{(1)})$).
The Weyl vector is $({1\over 6},1,{1\over 6})$.
The automorphic form $B_{\xi_{0,12}^{(2)}+\xi_{0,12}^{(3)}}$
is $\Delta_1^{(12)}(z_1,z_2,z_3)=
\Delta_1(z_1,2z_2,z_3)$ of the weight 1 where $\Delta_1$
corresponds to $\geg(\xi_{0,3}^{(1)})$. Thus, Fourier expansion of
$\Delta_1^{(12)}$ is
$$
\Delta_1^{(12)}=\sum_{M\ge 1}
\sum\Sb n,\,m >0,\,l\in \bz\\
\vspace{0.5\jot} n,\,m\equiv 1\mod 6\\
\vspace{0.5\jot} 4nm-3l^2=M^2\endSb
\hskip-4pt
\biggl(\dsize\frac{-4}{l}\biggr)
\biggl(\dsize\frac{12}{M}\biggr)
\sum\Sb a|(n,l,m)\endSb \biggl(\dsize\frac{\,6\,}{a}\biggr)
q^{n/6}r\,s^{m/6}\,.
$$
This case is equivalent to $\geg(\xi_{0,3}^{(1)})$ above.

\vskip20pt

\centerline{\bf Case $t=16$}

\vskip10pt

{\it The Algebra $\geg(\xi_{0,16}^{(1)})$.} 
The chamber $\M=[s_{\alpha_3},s_{\alpha_4}](\M_0)$ is an infinite
polygon with zero angles touching a horosphere with
the center $\br_{++}\rho$ where $\rho=({1\over 8},1,{1\over 8})$
is the Weyl vector. The set
$$
P(\M)=[s_{\alpha_3},s_{\alpha_4}](\alpha_1,2\alpha_2,\alpha_5)
$$
with the group of symmetries $[s_{\alpha_3},s_{\alpha_4}]$ which is
$D_\infty$. The generalized Cartan matrix is
$$
\left({(\alpha,\,\alpha^\prime)\over 2}\right),\ \ \
\alpha,\ \alpha^\prime \in P(\M),
$$
which is the same as for $\geg(\xi_{0,4}^{(1)})$.
The automorphic form $B_{\xi_{0,16}^{(1)}}$ coincides with
$\Delta_{1/2}^{(16)}(z_1,z_2,z_3)=
\Delta_{1/2}(z_1,2z_2,z_3)$ of the weight 1/2
where $\Delta_{1/2}$ corresponds to $\geg(\xi_{0,4}^{(1)})$.
Thus, Fourier expansion of $\Delta_{1/2}^{(16)}$ is
$$
\Delta_{1/2}^{(16)}=
\frac{1}2\sum\Sb n,m\in \bz \endSb
\,\biggl(\dsize\frac{-4}{n}\biggr)\biggl(\dsize\frac{-4}{m}\biggr)
q^{n^2/8}r^{nm}s^{m^2/8}\,.
$$
This case is equivalent to $\geg(\xi_{0,4}^{(1)})$.

\vskip20pt

\centerline{\bf Case $t=36$}

\vskip10pt

{\it The Algebra $\geg(\xi_{0,36}^{(2)})$.} The chamber
$\M=[s_{\alpha_3},s_{\alpha_4}](\M_0)$
is the infinite periodic polygon
with angles
$\dots$, $0$, $\pi/2$, $0$, $0$, $0$, $0$, $\pi/2$, $0$, $\dots$,
with the center at the Weyl vector
$\rho=({1\over 24},{1\over 2},{1\over 24})$ at infinity.
The set
$$
P(\M)=[s_{\alpha_3},\,s_{\alpha_4}]
(\alpha_1,\,\alpha_2,\,\alpha_5,\,\alpha_6),\ \ \
P(\M)_{\1o}=[s_{\alpha_3},\,s_{\alpha_4}](\alpha_2)
$$
with the group of symmetries $[s_{\alpha_3},\,s_{\alpha_4}]$ which is
$D_\infty$. The generalized Cartan matrix is
$$
\left({2(\alpha,\alpha^\prime)\over (\alpha,\,\alpha)}\right),
\ \ \ \alpha,\ \alpha^\prime \in P(\M).
$$
The automorphic form $B_{\xi_{0,36}^{(2)}}$ coincides with
$D_{1/2}$ of the weight 1/2 with Fourier expansion
$$
D_{1/2}=\frac{1}2\sum_{m,n\in \bz}
\biggl(\frac{12}n\biggl)\biggl(\frac{12}m\biggl)
\,q^{{n^2}/{24}}r^{{nm}/2}s^{{m^2}/{24}}.
$$
See \cite{41; \thetag{5.1.3}}.

\vskip20pt

\head
3. Appendix: On Jacobi modular forms with integral Fourier coefficients
\endhead

Here we concentrate on calculational aspects of Jacobi modular forms
(or just Jacobi forms) which we use in this paper. We try to avoid
complicated subjects related with modular forms and try to be as
short as possible. We follow \cite{41}, \cite{33}, \cite{34}.

We shall consider Jacobi forms with respect to Jacobi group
$\Gamma^J$ which is a semi-direct product
$\Gamma^J=SL_2(\bz)\ltimes H(\bz)$ where $H(\bz )$ is the integral
Heisenberg group. The group $H(\bz )$ is the central extension of
$\bz^2$ by $\bz$ which is the center of $\Gamma^J$.
We denote by $(\lambda,\,\mu)$ an element of
the $\bz^2$ and by $\kappa$ an element of $\bz$.
The Jacobi group $\Gamma^J$ can be identified with
$\Gamma_\infty/\{\pm E_4\}$ where $\Gamma_\infty$ is a maximal
parabolic subgroup of $Sp_4(\bz)$, which consists of all elements
preserving a line.

The binary (to $\{\pm 1\}$) character $v_H$ on $H(\bz)$ which is
$$
v_H([\lambda,\mu;\,\kappa]):=(-1)^{\lambda+\mu+\lambda\mu+\kappa}
$$
can be extended to a binary character $v_J$ of the Jacobi group
if one puts $v_J|_{SL_2(\bz)}=1$.

\definition{Definition 3.1}Let $k\in \bz/2$,
$t\in \bz/2$, $t\ge 0$ and
$v:SL_2(\bz)\to \bc^*$ a character
(or a multiplier system) of finite order of $SL_2(\bz)$.
A holomorphic  function
$\phi(\tau ,z)$ on $\bh\times \bc$ (where
$\bh$ is the upper-half plane $\text{Im\,}\tau>0$) is called
a  {\it weak Jacobi form} of
weight $k$ and index $t$ with respect to $\Gamma^J$ with the character
(or the multiplier system) $v$ if
$$
\phi(\frac{a\tau+b}{c\tau+d},\,\frac z{c\tau+d})=
v(\left(\smallmatrix a&b\\c&d\endsmallmatrix\right))
(c\tau+d)^k e^{\frac {2\pi i t c z^2}{c\tau+d}}\phi(\tau,z)
\qquad (\left(\smallmatrix a&b\\c&d\endsmallmatrix\right) \in SL_2(\bz)),
\tag{3.1}
$$
and
$$
\phi(\tau, z+\lambda \tau+ \mu)=
(-1)^{2t(\lambda+\mu)}
e^{-2\pi i t(\lambda^2 \tau+2\lambda z)}\phi(\tau,z)
\qquad (\lambda, \mu\in \bz),
\tag{3.2}
$$
and it has Fourier expansion of type
$$
\phi(\tau ,z)=\sum\Sb n\ge 0,\ l \in t+\bz\\
\vspace{0.5\jot} \endSb
f(n,l) \exp{(2\pi i(n\tau +lz))}.
\tag{3.3}
$$
We denote the space of all weak Jacobi forms of weight $k$
and index $t$ with $v$ as $J_{k,t}(v)$.
\enddefinition

Further we denote $q=\exp(2\pi i\tau)$ and $r=\exp(2\pi i z)$,
thus $\exp{(2\pi i(n\tau +lz))}=q^nr^l$.

The Dedekind $\eta$-function is
$$
\eta(\tau)=
q^{\frac 1{24}}\,\prod_{n\ge 1} (1-q^n)=
\sum_{n\in \bn} \biggl(\frac{12}n\biggl)q^{{n^2}/{24}},
\tag3.4
$$
where
$$
\biggl(\frac{12}n\biggl)=
\cases
\hphantom{-}1\  \text{ if }\  n\equiv \pm 1 \operatorname{mod} 12,\\
-1\ \text{ if }\  n\equiv \pm 5 \operatorname{mod} 12,\\
\hphantom{-}0\  \text{ if }\  (n,\,12)\ne 1\,.
\endcases
$$
The $\eta(\tau)$ is $SL_2(\bz)$-modular form of
the weight $1/2$ with some
multiplier system $v_\eta$ which takes values in 24th roots of unity.
The function $\Delta_{12}=\eta^{24}$ (the Ramanujan function)
can be considered as a weak Jacobi form
$\Delta_{12}(\tau,\,z)=\Delta_{12}(\tau)$
of weight $12$ and index $0$, i.e. $\Delta_{12}\in J_{12,0}$.
It had been used in \thetag{1.3.5}. The function $\Delta_{12}$ is
equal to zero only at infinity $q=0$.

We remind the classical Eisenstein series
$$
E_4(\tau )=1+240\sum_{n=1}^{\infty}{\sigma_3(n)q^n},\hskip1cm
E_6(\tau )=1-504\sum_{n=1}^{\infty}{\sigma_5(n)q^n},
\tag{3.5}
$$
where
$\sigma_k(n)=\sum_{m|n}{m^k}$.
They are $SL_2(\bz)$-modular forms of the weight $4$ and $6$ respectively.
They give weak Jacobi forms $E_4\in J_{4,0}$ and $E_6\in J_{6,0}$
of index 0. We have $\Delta_{12}=(E_4^3-E_6^2)/1728$
which gives another formula for $\Delta_{12}$.
The modular forms $E_4$, $E_6$ and $\Delta_{12}$ generate
over $\bz$ the ring of holomorphic $SL_2(\bz)$-modular forms of
even weight with integral Fourier coefficients. Over $\bq$ similar
ring has two free generators $E_4$ and $E_6$.

A function $\phi(\tau,\,z)$ is called a {\it nearly holomorphic
Jacobi form} if $\Delta^N\phi(\tau,\,z)$ for some $N\ge 0$
is a weak Jacobi form ($\Delta^N\phi(\tau,\,z)\in J_{k,t}(v)$
for some $k$, $t$). Nearly holomorphic Jacobi forms is the most general
class of Jacobi forms which we consider here. They may have poles at
infinity $q=0$.  Their Fourier coefficients $f(n,l)$ depend only on
the norm $4tn-l^2$ and $\pm l\mod 2t$; the norm $4tn-l^2$
of non-zero Fourier coefficients $f(n,l)$ is bounded from below.
A nearly holomorphic Jacobi form is holomorphic at infinity iff
its non-zero Fourier coefficients have non-negative norm.
We denote by $J_{k,t}^{nh}$ the space of
all nearly holomorphic Jacobi forms of weight
$k\in \bz/2$ and index $t\in \bz/2$, $t\ge 0$ with
the trivial $SL_2(\bz)$-character. Respectively, $J_{k,t}$ denotes
its subspace of all weak Jacobi forms.

Below we shall describe generators of the ring
$J_{0,\ast}^{\bz}$ of weak Jacobi forms of zero weight,
integral index and with integral Fourier coefficients.

The Jacobi theta-series
$$
\align
\vartheta(\tau ,z)&=\hskip-2pt\sum\Sb n\equiv 1\, mod\, 2 \endSb
\,(-1)^{\frac{n-1}2}
\exp{(\frac{\pi i\, n^2}{4} \tau +\pi i\,n z)}=
\sum_{m\in \bz}\,\biggl(\frac{-4}{m}\biggr)\, q^{{m^2}/8}\,r^{{m}/2}\\
{}&=
-q^{1/8}r^{-1/2}\prod_{n\ge 1}\,(1-q^{n-1} r)(1-q^n r^{-1})(1-q^n)
\tag{3.6}
\endalign
$$
is a holomorphic Jacobi form of weight $1/2$ and index $1/2$
with the multiplier system $v_\eta^3$  where
$$
\biggl(\frac{-4}n\biggl)=
\cases
\pm 1\  &\text{ if }\  n\equiv \pm 1 \operatorname{mod} 4,\\
\ 0\      &\text{ if }\  (n,\,4)\ne 1\,.
\endcases 
$$
Using $\vartheta(\tau, z)$, we get weak Jacobi forms
$$
\split
&{}\phi_{0, \frac3{2}}(\tau, z) =
\frac{\vartheta(\tau, 2z)}{\vartheta(\tau, z)}\\
=& r^{-\frac 1{2}}
\prod_{n\ge 1}(1+q^{n-1}r)(1+q^{n}r^{-1})(1-q^{2n-1}r^2)
(1-q^{2n-1}r^{-2})\in J_{0,\,\frac 3{2}},
\endsplit
\tag{3.7}
$$
$$
{}\phi_{-1, \frac 1{2}}(\tau,z)=\frac{\vartheta(\tau ,z)}{\eta(\tau)^3}
=-r^{-1/2}\prod_{n\ge 1}\,(1-q^{n-1} r)(1-q^n r^{-1})(1-q^n)^{-2}
\in J_{-1,\frac 1{2}}
\tag{3.8}
$$
and weak Jacobi forms
$$
\phi_{0,3}(\tau, z)=\phi_{0,\,\frac3{2}}(\tau, z)^2=
\frac{\vartheta(\tau, 2z)^2}{\vartheta(\tau, z)^2}\in J_{0,3},
\tag{3.9}
$$
$$
\phi_{-2, 1}(\tau, z)=\phi_{-1, \frac 1{2}}(\tau, z)^2=
\frac{\vartheta(\tau, z)^2}{\eta(\tau)^6}
\in J_{-2,\,1}.
\tag{3.10}
$$
One can define two other weak Jacobi forms with integral Fourier
coefficients
$$
\align
\phi_{0,2}(\tau ,z)=
{\tsize\frac{1}2} \eta(\tau )^{-4}
\sum_{m\,,n\in \bz}
{(3m-n)}\biggl(\frac{-4}m\biggl)\biggl(\frac{12}n\biggl)
q^{\frac{3m^2+n^2}{24}}r^{\frac{m+n}2}\in J_{0,2}\,,
\tag{3.11}
\endalign
$$
$$
\split
\phi_{0,4}(\tau ,z)=\frac{\vartheta(\tau, 3z)}{\vartheta(\tau, z)}=&
r^{-1}\prod_{m\ge 1}(1+q^{m-1}r+q^{2m-2}r^{2})
(1+q^{m}r^{-1}+q^{2m}r^{-2})\\
\times &\prod\Sb n\equiv 1,2 \mod 3\\ n\ge 1\endSb
(1-q^{n}r^3) (1-q^{n}r^{-3})\ \in J_{0,4}.
\endsplit
\tag{3.12}
$$
A weak Jacobi form $\phi_{0,1}\in J_{0,1}$ with integral
Fourier coefficients is defined by the relation
$$
4\phi_{0,4}=\phi_{0,1}\phi_{0,3}-\phi_{0,2}^2.
\tag{3.13}
$$
The Jacobi forms $\phi_{0,1}$ and $\phi_{-2,1}$ were introduced
in \cite{24} using different definition. The Jacobi forms
$\phi_{0,\frac 1{2}}$, $\phi_{-1,\frac 1{2}}$,
$\phi_{0,2}$, $\phi_{0,3}$, $\phi_{0,4}$ were introduced in
\cite{41}.

By \cite{33}, the ring $J_{0,*}^\bz$ of all weak Jacobi
forms of weight $0$ with integral index $*$ and
with integral Fourier coefficients is generated over $\bz$
by the weak Jacobi forms $\phi_{0,1}$, $\phi_{0,2}$, $\phi_{0,3}$
and $\phi_{0,4}$ with the relation \thetag{3.13}.

Eisenstein-Jacobi series $E_{4,1}\in J_{4,1}$, $E_{4,2}\in J_{4,2}$,
$E_{6,1}\in J_{6,1}$, $E_{6,2}\in J_{6,2}$, $E_{6,3}\in J_{6,3}$
(see \cite{24} about general results on Eisenstein--Jacobi series)
which have integral Fourier coefficients (and the Fourier coefficient
$1$ for $q^0r^0$) can be found by relations
$$
E_4\phi_{0,1}-E_6\phi_{-2,1}=12E_{4,1},
\tag{3.14}
$$
$$
E_6\phi_{0,1}-E_4^2\phi_{-2,1}=12E_{6,1},
\tag{3.15}
$$
$$
E_{4,1}\phi_{0,1}-E_{6,1}\phi_{-2,1}=12E_{4,2},
\tag{3.16}
$$
$$
E_{6,1}\phi_{0,1}-E_4E_{4,1}\phi_{-2,1}=12E_{6,2},
\tag{3.17}
$$
$$
E_{4,1}\phi_{0,2}-E_4\phi_{0,3}=2E_{4,3},
\tag{3.18}
$$
$$
E_{6,1}\phi_{0,2}-E_6\phi_{0,3}=2E^\prime_{6,3},
\tag{3.19}
$$
$$
E^\prime_{6,3}=E_{6,3}+\frac{22}{61}\Delta_{12}\phi_{-2,1}^3.
\tag{3.20}
$$
By \cite{34}, the ring $J_{*,*}^\bz$ of weak Jacobi forms of integral
weight and integral index with integral Fourier coefficients
is generated over $\bz$ by the weak Jacobi forms
$E_4$, $E_6$, $\Delta_{12}$, $E_{4,1}$, $E_{4,2}$, $E_{4,3}$, $E_{6,1}$,
$E_{6,2}$, $E_{6,3}^\prime$, $\phi_{0,1}$, $\phi_{0,2}$, $\phi_{0,3}$,
$\phi_{0,4}$ and $\phi_{-2,1}$.
It was proved in \cite{24} that similar ring
$J_{*,*}^\bq$ over $\bq$ is generated by free generators
$E_4$, $E_6$, $\phi_{0,1}$ and $\phi_{-2,1}$.

Using $\phi_{0, \frac3{2}}(\tau,z)\in J_{0, \frac3{2}}$ and
$\phi_{-1, \frac 1{2}}(\tau,z) \in J_{-1, \frac 1{2}}$, one
can get similar results for half-integral index. Using
$\Delta_{12}$, one can generalize these results for nearly holomorphic
Jacobi forms with integral Fourier coefficients.
For example, any nearly holomorphic Jacobi form $\phi_{0,t}\in J_{0,t}^{nh}$,
$t\in \bn$, with integral Fourier coefficients needed in  Theorem 2.2.1
can be written as
$$
\phi_{0,t}={P(E_4,E_6,\Delta_{12},
E_{4,1},E_{4,2},E_{4,3}, E_{6,1}, E_{6,2},
E_{6,3}^\prime, \phi_{0,1}, \phi_{0,2}, \phi_{0,3}, \phi_{0,4},
\phi_{-2,1})\over \Delta_{12}^N}
$$
where $P$ is a polynomial with integral coefficients and $N\ge 0$.


\Refs
\widestnumber\key{70}

\ref
\key 1
\by W.L. Baily 
\paper Fourier--Jacobi series
\inbook Algebraic groups and discontinuous subgroups.
Proc. Symp. Pure Math. Vol. IX
\eds Borel A., Mostow  G.D.
\publ Amer. Math. Soc.
\publaddr Providence, Rhode Island
\yr 1966
\pages 296--300
\endref


\ref
\key 2
\by R. Borcherds
\paper Vertex algebras, Kac--Moody algebras, and the monster
\jour Proc. Natl. Acad. Sci. USA
\vol 83
\yr 1986
\pages 3068 -- 3071
\endref

\ref
\key 3
\by  R. Borcherds 
\paper Generalised Kac--Moody algebras
\jour J. of Algebra
\vol 115
\yr 1988
\pages 501--512
\endref

\ref
\key 4
\by  R. Borcherds
\paper The monster Lie algebra
\jour Adv. Math.
\vol 83
\yr 1990
\pages 30--47
\endref

\ref
\key 5
\by  R. Borcherds 
\paper The monstrous moonshine and monstrous Lie superalgebras
\jour Invent. Math.
\vol 109
\yr 1992
\pages 405--444
\endref

\ref
\key 6
\by  R. Borcherds
\paper Sporadic groups and string theory
\inbook Proc. European Congress of Mathem. 1992
\pages 411--421
\endref

\ref
\key 7
\by  R. Borcherds
\paper Automorphic forms on $O_{s+2,2}$ and
infinite products
\jour Invent. Math. \vol 120
\yr 1995
\pages 161--213
\endref

\ref
\key 8
\by  R. Borcherds
\paper The moduli space of Enriques surfaces and the fake monster Lie
superalgebra
\jour Topology
\vol 35
\issue 3
\yr 1996
\pages 699-710
\endref


\ref
\key 9
\by  R. Borcherds
\paper Automorphic forms with singularities on Grassmanians
\jour Invent. Math.
\vol 132
\issue 3
\yr 1998
\pages 491-562
\moreref alg-geom/9609022
\endref

\ref
\key 10
\by  R. Borcherds
\paper What is moonshine?
\inbook Proc. Int. Congr. Math. Berlin 1998
\vol 1
\pages 607--615
\moreref math.QA/9809110
\endref

\ref
\key 11
\by  R. Borcherds
\paper Vertex algebras
\inbook Topological field theory, primitive forms and related
topics (Kyoto, 1996)
\pages 35-77
\publ Progr. Math. 160.
\publaddr Birkh\"auser Boston, Boston, MA
\yr 1998
\moreref q-alg/9706008
\endref

\ref
\key 12
\by  R. Borcherds
\paper Reflection groups of Lorentzian lattices
\jour Duke Math. J.
\vol 104
\yr 2000
\issue 2
\pages 319--366
\moreref math.GR/9909123
\endref

\ref
\key 13
\by R. Borcherds, L. Katzarkov, T. Pantev T, N.I. Shepherd-Barron 
\paper Families of $K3$ surfaces
\jour J. Algebraic Geom.
\vol 7
\yr 1998
\issue 1
\pages 183-193
\moreref alg-geom/9701013
\endref


\ref
\key 14
\by J.H. Bruinier
\paper Borcherdsprodukte und Chernsche Klassen von
Hirzebruch-Zagier-\-Zykeln
\jour Dissertation, Universit\"at Heidelberg
\yr 1998
\endref

\ref
\key 15
\by J.H. Bruinier
\paper
Borcherds products and Chern classes of Hirzebruch-Zagier divisors
\jour Invent. Math.
\vol 138
\yr 1999
\issue 1
\pages 51 -- 83
\endref

\ref
\key 16
\by  J.H. Bruinier
\paper
Borcherds products on $O(2,\,l)$ and Chern classes of Heegner divisors
\jour Habilitationsschrift, Universit\"at Heidelberg
\yr 2000
\endref

\ref
\key 17
\by G.L. Cardoso
\paper Perturbative gravitational couplings and Siegel
modular forms in $D=4$, $N=2$ string models
\jour Nucl. Phys. Proc. Suppl.
\vol 56B
\yr 1997
\pages 94-101
\moreref hep-th/9612200
\endref


\ref
\key 18
\by G.L. Cardoso, G. Curio, D. Lust
\paper Perturbative coupling and modular forms in $N=2$ string models
with a Wilson line
\jour Nucl. Phys.
\vol B491
\yr 1997
\pages 147--183
\moreref hep-th/9608154
\endref


\ref
\key 19
\by J.H. Conway 
\paper
The automorphism group of the 26 dimensional even Lorentzian lattice
\jour J. Algebra
\yr 1983
\vol 80
\pages 159--163
\endref


\ref
\key 20
\by J.H.Conway, S. Norton 
\paper Monstrous moonshine
\jour Bull. London Math. Soc.
\vol 11
\yr 1979
\pages 308-339
\endref

\ref
\key 21
\by R. Dijkgraaf
\paper The mathematics of fivebranes
\inbook Proc. Int. Congr. Math. Berlin 1998
\vol 3
\pages 133--142
\moreref hep-th/9810157
\endref

\ref
\key 22
\by R. Dijkgraaf, G. Moore, E. Verlinde, H. Verlinde 
\paper Elliptic genera of symmetric products and second
quantized strings
\jour Commun. Math. Phys.
\yr 1997
\pages 197--209
\moreref hep-th/9608096
\endref

\ref
\key 23
\by R. Dijkgraaf, E. Verlinde, H. Verlinde 
\paper Counting dyons in $N=4$ string theory
\jour Nucl. Phys.
\yr 1997
\pages 543--561
\moreref hep-th/9607026
\endref

\ref
\key 24
\by M. Eichler,  D. Zagier 
\book The theory of Jacobi forms
\yr 1985
\publ Progress in Math. 55, Birkh\"auser
\endref

\ref
\key 25
\by I.B. Frenkel, J. Lepowsky, A. Meurmann 
\book Vertex operator algebras and the monster
\publ Academic Press
\yr 1988
\publaddr Boston, MA
\endref

\ref
\key 26
\by H. Garland, J. Lepowsky, 
\paper Lie algebra homology and the Macdonald--Kac formulas
\jour Invent. Math.
\vol 34
\yr 1976
\pages 37-76
\endref

\ref
\key 27
\by P. Goddard
\paper The work of Richard Ewen Borcherds
\inbook Proc. Int. Congr. Math. Berlin 1998
\vol 1
\pages 99--108
\endref

\ref
\key 28
\by P. Goddard, C.B. Thorn
\paper Compatibility of the dual Pomeron with unitarity and
the absence of ghosts in the dual resonance model
\jour Phys. Lett.
\vol B40
\issue 2
\pages 235 --238
\yr 1972
\endref

\ref
\key 29
\by V.A. Gritsenko 
\paper Jacobi functions of n-variables
\jour Zap. Nauk. Sem. LOMI
\vol 168
\yr 1988
\pages 32--45
\lang Russian
\transl\nofrills English transl. in
\jour J\. Soviet Math\.
\vol 53
\yr 1991
\pages 243--252
\endref

\ref
\key 30
\by V.A. Gritsenko
\paper Arithmetical lifting and its applications
\inbook Number Theory. Proceedings of Paris Seminar  1992--93
\eds S. David
\publ Cambridge Univ. Press
\yr 1995
\pages 103--126
\endref


\ref
\key 31
\by V.A. Gritsenko
\paper Modular forms and moduli spaces of Abelian and K3 surfaces
\jour Algebra i Analyz
\vol 6:6
\yr 1994
\pages 65--102
\lang Russian
\transl\nofrills  English transl. in
\jour St.Petersburg Math. Jour.
\vol 6:6
\yr 1995
\pages 1179--1208
\endref

\ref
\key 32
\by V.A. Gritsenko 
\paper Irrationality of the moduli spaces of polarized 
Abelian surfaces
\jour The International Mathematics Research Notices
\vol 6
\yr 1994
\pages  235--243,
In  full form  in
``{\it Abeli\-an varieties}'',  Proc. of the  Egloffstein conference
(1993)  de Gruyter, Berlin, 1995, pp. 63--81
\endref

\ref
\key 33
\by V.A. Gritsenko
\paper Elliptic genus of Calabi-Yau manifolds and Jacobi and
Siegel modular forms
\jour St. Petersburg Math. J.
\vol 11:5
\yr 1999
\pages 100--125
\moreref math.AG/9906190
\endref

\ref
\key 34
\by V.A. Gritsenko 
\paper Complex vector bundles and Jacobi forms
\jour Preprint MPI
\vol 76
\yr 1999
\moreref math.AG /9906191
\endref

\ref
\key 35
\by V.A. Gritsenko, K. Hulek 
\paper Minimal Siegel modular threefolds
\yr 1998
\jour Mathem. Proc. Cambridge Phil. Soc.
\vol 123
\pages 461--485
\moreref alg-geom/9506017
\endref

\ref
\key 36
\by V.A. Gritsenko, V.V. Nikulin 
\paper Siegel automorphic form correction of some Lorentzi\-an
Kac--Moody Lie algebras
\jour Amer. J. Math.
\yr 1997
\vol 119
\issue 1
\pages 181--224
\moreref  alg-geom/
\newline
9504006
\endref


\ref
\key 37
\by V.A. Gritsenko, V.V. Nikulin
\paper Siegel automorphic form correction of a Lorentzian
Kac--Moody algebra
\jour C. R. Acad. Sci. Paris S\'er. A--B
\vol 321
\yr 1995
\pages 1151--1156
\endref

\ref
\key 38
\by V.A. Gritsenko, V.V. Nikulin 
\paper K3 surfaces, Lorentzian Kac--Moody algebras and
mirror symmetry
\jour  Math. Res. Lett. \yr 1996 \vol 3 \issue 2 \pages 211--229;
\nofrills  alg-geom/9510008.
 \endref

\ref
\key 39
\by V.A. Gritsenko, V.V. Nikulin
\paper The Igusa modular forms and `the simplest'
Lorentzian Kac--Moody algebras
\jour Mat. Sb.
\yr 1996
\vol 187
\issue 11
\pages 27--66
\lang Russian
\transl\nofrills English transl. in
\jour Sb. Math.
\vol 187
\issue 11
\pages 1601--1641
\moreref alg-geom/9603010
\endref


\ref
\key 40
\by V.A. Gritsenko, V.V. Nikulin 
\paper Automorphic forms and Lorentzian Kac--Moody algebras. Part I
\jour Intern. J. Math. \yr 1998
\vol 9 \issue 2 \pages 153--199
\moreref alg-geom/9610022
\endref

\ref
\key 41
\by V.A. Gritsenko, V.V. Nikulin 
\paper Automorphic forms and Lorentzian Kac--Moody algebras. Part II
\jour Intern. J. Math. \yr 1998
\vol 9 \issue 2 \pages 201--275
\moreref alg-geom/9611028
\endref

\ref
\key 42
\by V.A. Gritsenko, V.V. Nikulin 
\paper The arithmetic mirror symmetry and Calabi--Yau manifolds
\jour Comm. Math. Phys.
\yr 2000
\pages 1--11
\vol 210
\moreref alg-geom/9612002
\endref


\ref
\key 43
\by V.A. Gritsenko, V.V. Nikulin 
\paper A lecture about classification of Lorentzian Kac--Moody
algebras of the rank three
\jour Preprint Newton Inst. Math. Sci.
\vol NI00036-SGT
\yr 2000
\pages 1--26
\moreref alg-geom/0010329
\endref

\ref
\key 44
\by V.A. Gritsenko, V.V. Nikulin 
\paper On the classification of meromorphic automorphic products
and related Lorentzian Kac--Moody algebras with respect to the
extended paramodular group
\toappear
\endref

\ref
\key 45
\by J. Harvey, G. Moore
\paper Algebras, BPS-states, and strings
\jour Nucl. Physics.
\vol B463
\yr 1996
\pages 315--368
\moreref hep-th/9510182
\endref

\ref
\key 46
\by J. Igusa
\paper On Siegel modular forms of genus two (II)
\jour Amer. J. Math.
\yr 1964
\vol 84
\issue 2
\pages 392--412
\endref

\ref
\key 47
\by V. Kac
\book Infinite dimensional Lie algebras
\yr 1990
\publ Cambridge Univ. Press
\endref


\ref
\key 48
\by V. Kac
\paper Infinite-dimensional algebras, Dedekind's $\eta$-function, classical
M\"obius function and the very strange formula
\jour Adv. Math.
\vol 30
\yr 1978
\pages 85--136
\endref

\ref
\key 49
\by V. Kac
\book Vertex algebras for beginners
\yr 1998
\publ Univ. Lect. Series (Providence, R.I.). Vol. 10
\publaddr Amer. Math. Soc.
\endref

\ref
\key 50
\by V. Kac, M. Wakimoto 
\paper Integrable highest weight modules over affine superalgebras
and number theory
\inbook Lie theory and geometry
\publ Progr. Math., 123, Birkh\"auser Boston
\publaddr Boston, MA
\yr 1994
\pages 415-456
\endref

\ref
\key 51
\by T. Kawai
\paper $N=2$ heterotic string threshold correction, K3
surfaces and generalized Kac--Moody superalgebra
\jour Phys. Lett.
\vol B371
\yr 1996
\page 59--64
\moreref hep-th/9512046
\endref

\ref
\key 52
\by T. Kawai 
\paper String duality and modular forms
\jour Phys. Lett.
\vol B397
\yr 1997
\pages 51--62
\moreref hep-th/ 9607078
\endref

\ref
\key 53
\by T. Kawai, K. Yoshioka
\paper String partition functions and infinite products
\jour  Adv. Theor. Math. Phys.
\vol 4
\yr 2000
\pages 397-485
\moreref hep-th/0002169
\endref

\ref
\key 54
\by Sh. Kond$\overline{\text{o}}$
\paper On the Kodaira dimension of the moduli space of $K3$ surfaces. II
\jour Compositio Math.
\vol 116
\yr 1999
\issue 2
\pages 111 -- 117
\endref

\ref
\key 55
\by H. Maass
\paper \"Uber ein Analogon zur Vermutung von Saito-Kurokawa
\jour Invent. math.
\yr 1980
\vol 60
\pages 85--104
\endref


\ref
\key 56
\by G. Moore
\paper
String duality, automorphic forms and generalised Kac--Moody algebras
\jour Nucl. Phys. Proc. Suppl.
\vol 67
\yr 1998
\pages 56--67
\moreref hep-th/9710198
\endref

\ref
\key 57
\by V.V. Nikulin
\paper Integral symmetric bilinear forms and some of
their geometric applications
\jour Izv. Akad. Nauk SSSR Ser. Mat.
\vol  43
\yr 1979
\pages 111--177
\lang Russian
\transl\nofrills English transl. in
\jour Math. USSR Izv.
\vol 14
\yr 1980
\endref

\ref
\key 58
\by V.V. Nikulin 
\paper On the quotient groups of the automorphism groups of
hyperbolic forms by the subgroups generated by 2-reflections,
Algebraic-geometric applications
\jour Current Problems in Math. Vsesoyuz. Inst. Nauchn. i
Tekhn. Informatsii, Moscow
\yr 1981
\pages 3--114
\lang Russian
\transl\nofrills English transl. in
\jour J. Soviet Math.
\yr 1983
\vol 22
\pages 1401--1476
\endref

\ref
\key 59
\by V.V. Nikulin
\paper On arithmetic groups generated by
reflections in Lobachevsky spaces
\jour Izv. Akad. Nauk SSSR Ser. Mat.
\vol  44   \yr 1980 \pages 637--669
\lang Russian
\transl\nofrills English transl. in
\jour Math. USSR Izv.
\vol 16 \yr 1981
\endref

\ref
\key 60
\by V.V. Nikulin 
\paper On the classification of arithmetic groups generated by
reflections in
Lobachev\-sky spaces
\jour Izv. Akad. Nauk SSSR Ser. Mat.
\vol  45
\issue 1
\yr 1981
\pages 113--142
\lang Russian
\transl\nofrills English transl. in
\jour Math. USSR Izv.
\vol 18
\yr 1982
\endref


\ref
\key 61
\by V.V. Nikulin
\paper Discrete reflection groups in Lobachevsky spaces and
algebraic surfaces
\inbook Proc. Int. Congr. Math. Berkeley 1986
\vol  1
\pages 654--669
\endref

\ref
\key 62
\by V.V. Nikulin
\paper A lecture on Kac--Moody Lie algebras of the arithmetic type
\jour Preprint Queen's University, Canada
\vol \#1994-16,
\yr 1994 \nofrills ; alg-geom/9412003.
\endref

\ref
\key 63
\by V.V. Nikulin
\paper Basis of the diagram method for generalized reflection
groups in Loba\-chevsky spaces and algebraic surfaces with nef
anticanonical class
\jour Int. J. Mathem.
\vol 7 \issue 1
\yr 1996
\pages 71 -- 108
\endref

\ref
\key 64
\by V.V. Nikulin
\paper Reflection groups in Lobachevsky spaces and
the denominator identity for Lorentzian Kac--Moody algebras
\jour Izv. Akad. Nauk of Russia. Ser. Mat.
\vol  60
\issue 2
\yr 1996
\pages 73--106
\lang Russian
\transl\nofrills English transl. in
\jour Russian Acad. Sci. Izv. Math.
\nofrills ; alg-geom/9503003.
\endref

\ref
\key 65
\by V.V. Nikulin
\paper The remark on discriminants of K3 surfaces moduli as sets
of zeros of automorphic forms
\jour J. Math. Sci., New York
\vol 81
\issue 3
\yr 1996
\pages 2738--2743
\moreref alg-geom/9512018
\endref

\ref
\key 66
\by V.V. Nikulin
\paper $K3$ surfaces with interesting groups of automorphisms
\jour J. Math. Sci., New York
\vol 95
\issue 1
\pages 2028--2048
\yr 1999
\moreref  alg-geom/9701011
\endref

\ref
\key 67
\by V.V. Nikulin
\paper A theory of Lorentzian Kac--Moody algebras
\inbook Trudy Mezhdunar. konf. posvya\-shch. 90-letiyu so dnya
rozhdeniya L.S. Pontryagina, T.8: Algebra (Proc. Int. Conf.
Devoted to the 90-th Anniversary of L.S. Pontryagin, vol. 8: Algebra)
\publ VINITI (Itogi nauki i tekhniki.
Sovremennaya matematika i ee prilozheniya. Tematicheskie obzory,
vol. 69) \yr 1999 \pages 148-167
\publaddr Moscow
\lang Russian
\moreref math.AG/9810001
\endref


\ref
\key 68
\by V.V. Nikulin
\book On the classification of hyperbolic root systems of the
rank three
\publ Trudy Matem. Instit. V.A. Steklov,  T. 230
\publaddr Moscow
\yr 2000
\pages 255
\lang Russian
\transl\nofrills English translation in
\publ Proc. Steklov Math. Institute. Vol. 230
\publaddr Moscow
\yr 2000
\moreref alg-geom/9711032; alg-geom/9712033; math.AG/9905150
\endref

\ref
\key 69
\by U. Ray
\paper A character formula for generalised Kac--Moody superalgebras
\jour J. of Algebra
\vol 177
\yr 1995
\pages 154--163
\endref


\ref
\key 70
\by U. Ray 
\paper Generalized Kac--Moody algebras and some related topics
\jour Bull. Amer. Math. Soc.
\vol 38 \issue 1
\yr 2001
\pages 1--42
\endref

\ref
\key 71
\by \'E.B. Vinberg 
\paper The absence of crystallographic reflection groups in Lobachevsky
spaces of large dimension
\jour Trudy Moscow. Mat. Obshch.
\vol  47 \yr 1984  \pages 68 -- 102
\lang Russian
\transl\nofrills English transl. in
\jour Trans. Moscow Math. Soc.
\vol 47 \yr 1985
\endref


\ref
\key 72
\by  \'E.B. Vinberg 
\paper Hyperbolic reflection groups
\jour Uspekhi Mat. Nauk
\vol 40
\yr 1985
\pages 29--66
\lang Russian
\transl\nofrills English transl. in
\jour Russian Math. Surveys
\vol 40
\yr 1985
\endref

\ref
\key 73
\by  \'E.B. Vinberg 
\paper Discrete reflection groups in Lobachevsky spaces
\inbook Proc. Int. Congr. Math. Warsaw, 1983
\vol  1, 2
\pages 593--601
\endref

\endRefs
\enddocument

\end